\documentclass[10pt]{amsart}%
\usepackage{amsmath,amssymb,amscd}%
\usepackage{array}
\usepackage[new]{old-arrows}
\usepackage{xcolor}%
\usepackage{mathrsfs}%
\usepackage{bbm}
\usepackage[all]{xy}%
\usepackage{stmaryrd}
\usepackage{hyperref}
\usepackage[french,british,ngerman]{babel}
\usepackage{cancel}
\usepackage[normalem]{ulem}
\usepackage{colonequals}
\usepackage{comment}
\usepackage{csquotes}
\usepackage{cleveref}
\usepackage{tikz-cd}
\usepackage[cal=boondox]{mathalpha}
\usetikzlibrary{babel}
\usepackage[T1]{fontenc}

\theoremstyle{plain}
    \newtheorem{theorem}{Theorem}[section]
    \newtheorem{proposition}[theorem]{Proposition}
    \newtheorem{lemma}[theorem]{Lemma}
    \newtheorem{corollary}[theorem]{Corollary}

    \newtheorem{proposition-definition}[theorem]{Proposition-Definition}
    \newtheorem{theorem-definition}[theorem]{Theorem-Definition}
    \newtheorem{lemma-definition}[theorem]{Lemma-Definition}
    \newtheorem*{theorem*}{Theorem}
    \newtheorem*{conjecture*}{Conjecture}
    \newtheorem*{corollary*}{Corollary}
\theoremstyle{definition}
    \newtheorem{definition}[theorem]{Definition}
    \newtheorem{notation}[theorem]{Notation}
    \newtheorem{assumption}[theorem]{Assumption}
    \newtheorem{example}[theorem]{Example}
    \newtheorem{remark}[theorem]{Remark}

    \newtheorem*{remark*}{Remark}

\crefname{theorem}{Theorem}{Theorems}
\crefname{lemma}{Lemma}{Lemmata}
\crefname{corollary}{Corollary}{Corollaries}
\crefname{proposition}{Proposition}{Propositions}
\crefname{proposition-definition}{Proposition-Definition}{Proposition-Definitions}
\crefname{claim}{Claim}{Claims}
\crefname{theorem-definition}{Theorem-Definition}{Theorem-Definitions}
\crefname{lemma-definition}{Lemma-Definition}{Lemma-Definitions}
\crefname{definition}{Definition}{Definitions}
\crefname{conjecture}{Conjecture}{Conjectures}
\crefname{question}{Question}{Questions}
\crefname{example}{Example}{Examples}
\crefname{algorithm}{Algorithm}{Algorithms}
\crefname{remark}{Remark}{Remarks}
\crefname{assumption}{Assumption}{Assumptions}
    
\def\Alphabet{A,B,C,D,E,F,G,H,I,J,K,L,M,N,O,P,Q,R,S,T,U,V,W,X,Y,Z}
\def\alphabet{a,b,c,d,e,f,g,h,i,j,k,l,m,n,o,p,q,r,s,t,u,v,w,x,y,z}
\def\endpiece{xxx}
\def\makeAlphabet[#1]{\expandafter\makeA#1,xxx,}
\def\makealphabet[#1]{\expandafter\makea#1,xxx,}
\def\makeA#1,{\def\temp{#1}\ifx\temp\endpiece\else%
\mkbb{#1}\mkbbm{#1}\mkfrak{#1}\mkbf{#1}\mkcal{#1}\mkscr{#1}\mkbs{#1}\expandafter\makeA\fi}%
\def\makea#1,{\def\temp{#1}\ifx\temp\endpiece\else\mkbbm{#1}\mkfrak{#1}\mkbf{#1}\mkbs{#1}\expandafter\makea\fi}%
\def\mkbb#1{\expandafter\def\csname bb#1\endcsname{\mathbb{#1}}}
\def\mkbbm#1{\expandafter\def\csname bm#1\endcsname{\mathbbm{#1}}}
\def\mkfrak#1{\expandafter\def\csname fr#1\endcsname{\mathfrak{#1}}}
\def\mkbf#1{\expandafter\def\csname b#1\endcsname{\mathbf{#1}}}
\def\mkcal#1{\expandafter\def\csname c#1\endcsname{\mathcal{#1}}}
\def\mkscr#1{\expandafter\def\csname s#1\endcsname{\mathscr{#1}}}
\def\mkbs#1{\expandafter\def\csname bs#1\endcsname{{\boldsymbol{#1}}}}
\def\makeop[#1]{\xmakeop#1,xxx,}
\def\mkop#1{\expandafter\def\csname #1\endcsname{{\operatorname{#1}}}} %
\def\xmakeop#1,{\def\temp{#1}\ifx\temp\endpiece\else\mkop{#1}\expandafter\xmakeop\fi}%
\def\makeup[#1]{\xmakeup#1,xxx,}
\def\mkup#1{\expandafter\def\csname #1\endcsname{{\mathrm{#1}\,}}} %
\def\xmakeup#1,{\def\temp{#1}\ifx\temp\endpiece\else\mkup{#1}\expandafter\xmakeup\fi}%
\makeAlphabet[\Alphabet]
\makealphabet[\alphabet]
\makeop[Hom,Aut,Ext,Pic,Jac,Div,NS,Br,H,R,ord,BSD,T,E,Res,Frob,rk,nr,codim,div,triv,Gal,free,Tr,res,cores,id,CH,tor,sep,sh,h,Spec,Spf,Spa,coker,im,ker]
\makeup[Zar,fppf,rig,cris,conv,an,alg,red,Isoc,Val,reg,sing,rank]
  \makeatletter
    
    \@addtoreset{equation}{section}
  \makeatother

\makeatletter
\newcommand{\labitem}[2]{%
\def\@itemlabel{\textbf{#1}}
\item
\def\@currentlabel{#1}\label{#2}}
\makeatother

\makeatletter
\newenvironment{abstracts}{%
  \ifx\maketitle\relax
    \ClassWarning{\@classname}{Abstract should precede
      \protect\maketitle\space in AMS document classes; reported}%
  \fi
  \global\setbox\abstractbox=\vtop \bgroup
    \normalfont\Small
    \list{}{\labelwidth\z@
      \leftmargin3pc \rightmargin\leftmargin
      \listparindent\normalparindent \itemindent\z@
      \parsep\z@ \@plus\p@
      
      \itemsep\medskipamount
    }%
}{%
  \endlist\egroup
  \ifx\@setabstract\relax \@setabstracta \fi
}

\newcommand{\abstractin}[1]{%
  \otherlanguage{#1}%
  \item[\hskip\labelsep\scshape\abstractname.]%
}
\makeatother

\newcommand{\et}{\mathrm{\acute{e}t}}
\newcommand*{\sheafhom}{\sH\kern-.5pt \mathrm{om}}
\newcommand*{\sheafext}{\sE\kern-.5pt \mathrm{xt}}

\DeclareFontFamily{U}{wncy}{}
\DeclareFontShape{U}{wncy}{m}{n}{<->wncyr10}{}
\DeclareSymbolFont{mcy}{U}{wncy}{m}{n}
\DeclareMathSymbol{\Sha}{\mathord}{mcy}{"58}

\newcommand{\Char}{{\mathrm{char}}}

\definecolor{darkgreen}{rgb}{0,0.5,0}
\definecolor{darkblue}{rgb}{0,0,0.8}
\definecolor{darkred}{rgb}{0.8,0,0}

\begin{document}
\selectlanguage{british}
\title{Comparison of different Tate Conjectures}

\author{}
\address{}
\email{}
\urladdr{}
\thanks{}

\author{Veronika Ertl}
\address{(V.\,Ertl) \newline
\parbox{\linewidth}{Universität Regensburg, Fakultät für Mathematik, Universitätsstraße 31, 93053 Regensburg, Germany}
\parbox{\linewidth}{Instytut Matematyczny Polskiej Akademii Nauk, ulica Śniadeckich 8, 00-656 Warszawa, Poland}
\parbox{\linewidth}{Laboratoire de Mathématiques Nicolas Oresme,  Université de Caen Normandie, 6 boulevard maréchal Juin, 14032 Caen Cedex, France} 
}
\email{veronika.ertl@unicaen.fr}
\urladdr{\url{https://ertlvroni.github.io/}}

\author{Timo Keller}
\address{(T.\,Keller)\newline Institut für Mathematik, Universität Würzburg, Emil-Fischer-Strasse 30, 97074,
	Würzburg, Germany \newline
     Rijksuniveriteit Groningen, Bernoulli Institute, Bernoulliborg, Nijenborgh 9, 9747 AG Groningen, The Netherlands\newline
     Leibniz Universität Hannover, Institut für Algebra, Zahlentheorie und Diskrete Mathematik, Welfengarten 1, 30167 Hannover, Germany\newline
     Department of Mathematics, Chair of Computer Algebra, Universität Bay\-reuth, Germany}
\email{math@kellertimo.de}
\urladdr{\url{https://www.timo-keller.de}}

\author{Yanshuai Qin}
\address{(Y.\,Qin)\newline Universität Regensburg, Fakultät für Mathematik, Universitätsstraße 31, 93053 Regensburg, Germany}
\email{yanshuai.qin@ur.de}
\urladdr{\url{https://sites.google.com/view/ysqin/}}

\thanks{TK was supported by the Deutsche Forschungsgemeinschaft (DFG),
	Geschäftszeichen STO~299/18-1, AOBJ: 667349 and by the 2021 MSCA Postdoctoral Fellowship
	01064790 -- Ex\-pli\-cit\-Rat\-Points while working on this article.
YQ was supported by the DFG through CRC 1085 ``Higher Invariants''.
VE was supported by the DFG through CRC 1085 ``Higher Invariants'' and by the ERC-grant KAPIBARA under the European
Union's Horizon 2020 research and innovation programme (grant
agreement No.\,802787). }

\subjclass[2020]{11G40 (Primary) 11G05, 11G10, 14G10 (Secondary)}
\keywords{Birch--Swinnerton-Dyer conjecture, positive characteristic, higher-dimensional function fields, étale cohomology, rigid cohomology, Brauer groups, abelian varieties}

\date{\today}

\begin{abstracts}
\abstractin{british}
	For an abelian variety $A$ over a finitely generated field $K$ of characteristic $p > 0$, we prove that the algebraic rank of $A$ is  is less or equal to a suitably defined analytic rank. Moreover, we prove that equality, i.e., the BSD rank conjecture, holds for $A/K$ if and only if a suitably defined Tate--Shafarevich group of $A/K$ (1) has finite $\ell$-primary component for some/all $\ell \neq p$, or (2) finite prime-to-$p$ part, or (3) has $p$-primary part of finite exponent, or (4) is of finite exponent. We give an algorithm to verify those conditions for concretely given $A/K$.

\abstractin{french}
Pour une variété $A$ sur un corps de génération finie $K$ de caractéristique $p>0$, 
on montre que le rang algébrique de $A$ est au plus égal à un certain rang analytique. 
En outre,  on montre l'égalité, c'est-à-dire que la conjecture de BSD concernant le rang est vraie pour A/K, si et seulement si un certain groupe de Tate--Shafarevich de A/K satisfait une des conditions suivantes:
(1) sa composante $\ell$-primaire est finie pour un/tout $\ell\neq p$, (2) la partie qui  est premier à $p$ est finie, (3) la partie $p$-primaire est d'exposant fini, (4) il est d'exposant fini.
On donne un algorithme pour vérifier ces conditions pour la donné d'un $A/K$ concrète. 
\end{abstracts}

\maketitle
\selectlanguage{british}

\tableofcontents

\section*{Introduction}\label{sec: introduction}

\subsection*{Several conjectures by Tate}

In  \cite{Tate_1965}, Tate formulated several conjectures concerning
 algebraic cycles (\cite[Conj.\,1]{Tate_1965} and \cite[Conj.\,B+S-D, Conj.\,2, Conj.\,B+S-D+2]{Tate_1965}, cf.\,Conjectures~\hyperlink{conj1}{$\T^i(X,\ell)$},~\hyperlink{bsd}{B+S-D},~\hyperlink{conj: E}{$\E^i(X,\ell)$},~\hyperlink{conj2}{2} and~\hyperlink{BS3}{B+S-D+2} below) 
in terms of $L$-functions. He asked about the relations between these conjectures (cf.\,\cite[Conj. B+SD+2]{Tate_1965}) and proved some equivalences between them for smooth projective varieties over finite fields (cf.\,\cite[Thm. 2.9]{Tate_1994} and~\cite[Thm.\,3.1]{TateBSD1965}). In this paper, we consider his questions on the relations between these conjectures for smooth projective varieties over finitely generated fields. 
Let us
first recall Tate's formulations of these conjectures. (We will use Tate's notation throughout.)

Let $X$ be a smooth projective geometrically connected variety over a finitely generated field $K$ (i.e. finitely generated over a finite field or $\bbQ$).  
 It extends
to a projective and smooth morphism $\rho\colon \sX\longrightarrow \sY$ of schemes of finite type over $\bbZ$ with $\sY$ integral and regular, whose generic fibre is $X/K$ (compare \cite{Poonen_2017}). 
(We often call such a morphism a smooth model of $X/K$.) 
Let $|\sY|$ denote the set of closed points of $\sY$. For $y\in |\sY|$, let $q_y$ denote the cardinality of the residue field $k(y)$ of $y$.
Set
$$ 
P_{y,i}(T)\colonequals\det(1-\sigma_{y}^{-1}T| \H^i(\sX_{\bar{y}},\bbQ_\ell)),
$$
and the $L$-function 
$$
\Phi_i(s)\colonequals\prod_{y\in |\sY|}\frac{1}{P_{y,i}(q_y^{-s})},
$$
where $\sigma_y\in \Gal(\overline{k(y)}/k(y))$ is the arithmetic Frobenius and $\ell\neq \Char(k(y))$ is a prime number. 
(If we want to emphasise that $P_{y,i}$ or $\Phi_i(s)$ depend on   $\sX/\sY$ we add this to the notation.)  The product converges absolutely for $\Re(s)>\dim\sY+i/2$ (cf.\,\cite[Thm.\,1]{Serre_1965}). Replacing $\sY$ by an open dense subscheme will divide $\Phi_i(s)$ by a product which converges absolutely for $\Re(s)>\dim \sY+i/2-1$. So the zeros and poles of $\Phi_i$ in the strip
$$\dim\sY+i/2-1<\Re(s)\leqslant \dim\sY+i/2$$
depend only on $X/K$ and not on the choice of $\sX/\sY$.
The zeta function of $\sX$ is defined as
$$ 
\zeta(\sX,s)\colonequals\prod_{x\in |\sX|}\frac{1}{1-q_x^{-s}}.
$$
By definition and Grothendieck's trace formula,
$$ 
\zeta(\sX,s)=\prod_{y\in |\sY|}\zeta(\sX_y,s)=\prod_{i=0}^{2\dim X}\Phi_i( s)^{(-1)^i}.
$$
Let $\ell\neq \Char(K)$ be a prime number. Let $A^i(X)$ denote the group of classes of algebraic cycles of codimension $i$ on $X$, with coefficients in $\bbQ$, for $\ell$-adic homological equivalence. Let $N^i(X)\subseteq A^i(X)$ denote the group of classes of cycles that are numerically equivalent to zero. Let $d$ denote $\dim(X)$. Tate made the following conjectures:

\begin{conjecture*}[\hypertarget{conj1}{$\T^i(X,\ell)$}]
The cycle class map
$$
A^i(X)\otimes_\mathbb{Q}\mathbb{Q}_\ell\longrightarrow \H^{2i}(X_{K^\sep},\mathbb{Q}_\ell(i))^{G_K}
$$
is surjective. 
\end{conjecture*} 

\begin{conjecture*}[\hypertarget{conj: E}{$\E^i(X,\ell)$}]
$N^i(X)=0$, i.e., numerical equivalence is equal to $\ell$-adic homological equivalence for algebraic cycles of codimension $i$ on $X$.
\end{conjecture*}

\begin{conjecture*}[\hypertarget{bsd}{B+S-D}]
The rank of $\Pic^0_{X/K}(K)$ is equal to the order of the zero of $\Phi_1(s)$ at $s=\dim(\sY)$ (and of $\Phi_{2d-1}(s)$ at $s=\dim\sX-1$ by duality).
\end{conjecture*}

\begin{conjecture*}[\hypertarget{conj2}{2}]
The dimension of $A^i(X)/ N^i(X)$ is equal to the order of the poles of $\Phi_{2i}(s)$ at $s=\dim\sY+i$ (and of $\Phi_{2d-2i}(s)$ at $s=\dim\sX-i$ by duality).
\end{conjecture*}

Let $\sX$ be a regular scheme of finite type over $\bbZ$. Assume that its zeta function $\zeta(\sX,s)$ can be meromorphically continued to the point $s=\dim\sX-1$. Let $e(\sX)$ be the order of $\zeta(\sX,s)$ at that point, and put
\begin{equation} \label{eq:def z}\tag{$\ast$}
z(\sX)=\rk\, \H^0(\sX,\sO_{\sX}^*)-\rk\, \H^1(\sX,\sO_{\sX}^*)-e(\sX).
\end{equation}
Tate showed that $z(\sX)$ is a birational invariant, and any two of the following statements imply the third:
\begin{itemize}
\item[(i)] Conjecture~\hyperlink{bsd}{B+S-D} for $X/K$,
\item[(ii)]  Conjecture~\hyperlink{conj2}{2},  for $i=1$, for $X/K$ and
\item[(iii)]  The equality $z(\sX)=z(\sY)$.
\end{itemize}
Since $z(\sY)=0$ if $\sY$ is of dimension $1$, Tate conjectured the following:
\begin{conjecture*}[\hypertarget{BS3}{B+S-D+2}]
If $\sX$ is a regular scheme of finite type over $\bbZ$, then the order of $\zeta(\sX,s)$ at the point $s=\dim \sX-1$ is equal to $\rk\, \H^0(\sX,\sO_{\sX}^*)-\rk\, \H^1(\sX,\sO_{\sX}^*)$.
\end{conjecture*}
For a morphism $f\colon \sX\rightarrow \sY$ as above with $\dim \sY=1$, Tate wrote that Conjecture~\hyperlink{BS3}{B+S-D+2} seems equivalent to Conjecture~\hyperlink{bsd}{B+S-D} for $X/K$ and the Conjecture~\hyperlink{conj2}{2},  for $i=1$, for $X/K$. This was proved to be true in positive characteristic by Geisser \cite{Geisser_2021}. 

Motivated by this question of Tate, we define a Tate--Shafarevich group for an abelian variety over arbitrary finitely generated fields and study its relation with Brauer groups (the obstruction of Conjecture~\hyperlink{conj1}{$\T^1(X,\ell)$}). We prove a relation between Tate--Shafarevich groups and Brauer groups for $\rho\colon \sX\rightarrow \sY$. As a consequence of this relation, we obtain some equivalences involving the above conjectures of Tate. Then, we use this equivalence to prove that the finiteness of the exponent of the Tate--Shafarevich group for an abelian variety over a finitely generated field of positive characteristic is equivalent to Conjecture~\hyperlink{bsd}{B+S-D}.
This generalises the rank part of the theorems of Schneider \cite{Schneider_1982} and Kato--Trihan \cite{katotrihan2003} over a global function field to finitely generated fields with higher transcendence degrees.

\begin{remark*}
Conjecture~\hyperlink{conj1}{$\T^i(X,\ell)$} is often simply referred to as ``the Tate Conjecture'', 
whereas Conjecture~\hyperlink{BS3}{B+S-D+2} is sometimes also called ``the analytic Tate Conjecture''.

The two invariants appearing in Conjecture~\hyperlink{bsd}{B+S-D} are called the algebraic and the analytic rank of $X$, and denoted in this paper by $\rk_\alg$ and $\rk_\an$, respectively. 
The conjecture itself is often referred to as ``the rank part of the Birch--Swinnerton-Dyer Conjecture'' or simply ``the BSD rank Conjecture''. 
\end{remark*}

\subsection*{Main Theorems}

\subsubsection*{Tate--Shafarevich group and the BSD rank conjecture}
Let $A$ be an abelian variety over a finitely generated field $K$. Let $disc\Val(K)$ denote the set of discrete valuations on $K$. Define the Tate--Shafarevich group of $A/K$ as
$$
\Sha(A/K)\colonequals\ker\Big(\H^1(K,A)\rightarrow\prod_{v\in disc\Val(K)}\H^1(K_v^{\sh}, A)\Big),
$$
where $K_v^{\sh}$ denotes the Henselian field of $v$. 
In the case that $K$ is a global field, the usual Tate--Shafarevich group is a subgroup with a finite index of the one we defined. For an abelian group $M$ and a prime number $\ell$, denote by $M(\ell)$ the subgroup of elements with order equal to powers of $\ell$.
Combining Theorems~\ref{thm: bsdconj} and \ref{BSD in char p} we obtain:
\begin{theorem*}[\hypertarget{bsdconj}{A}]
Let $A$ be an abelian variety over a finitely generated field $K$ of characteristic $p>0$. Let $m$ denote the transcendence degree of $K$ over $\bbF_p$. Then,
$$ \rk\, A(K) \leqslant \ord_{s=m}\Phi_{1}(s),$$
and the following statements are equivalent.
\begin{itemize}
\item[(1)] $\Sha(A/K)(\ell)$ is finite for some prime $\ell \neq p.$

\item[(2)] $\Sha(A/K)(p)$ is of finite exponent.

\item[(3)] $\Sha(A/K)$ is of finite exponent.

\item[(4)] The Conjecture~\hyperlink{bsd}{B+S-D} for $A$ holds, i.e., the rank  of $A(K)$  is equal to $\ord_{s=m}\Phi_{1}(s)$.
\end{itemize}
\end{theorem*}
\begin{remark*}
The above theorem was proved for the case that $m=1$ and $A$ is a Jacobian of a curve by Artin--Tate \cite{TateBSD1965} and Milne \cite{milne1975conjecture}, for the case that $m=1$ by Schneider \cite{Schneider_1982} (prime to $p$-part) and Kato--Trihan \cite{katotrihan2003} ($p$-primary part) and for the case that $A$ extends to an abelian scheme over a smooth projective variety by the second author 
\cite{Keller_2019} (prime to $p$-part). Note that our Tate--Shafarevich group agrees with the one defined in~\cite{Keller_2016} by the second author 
if $A/K$ has good reduction everywhere, i.e., if it extends to an abelian scheme over a smooth projective model of $K$. Our proof follows Artin--Tate's approach. All of the previous work mentioned above also included the full BSD formula. Our work only concerns the rank part of the BSD conjecture. It seems interesting to formulate a full BSD formula in the general case.
\end{remark*}

\subsubsection*{Brauer groups and the Tate conjecture for divisors}
For any Noetherian scheme $X$,  the \emph{cohomological Brauer group}
$$
\Br(X)\colonequals\H^2(X,\mathbb{G}_m)_{\tor}
$$
is defined to be the torsion part of the \'etale cohomology group $\H^2(X,\bbG_m)$. In \cref{nsrank} we obtain the following statement: 
\begin{theorem*}[\hypertarget{firstthm}{B}]
Let $X$ be a smooth projective variety over a finitely generated field $K$ of characteristic $p>0$. Let $m$ denote the transcendence degree of $K$ over its prime field. Then 
$$
\dim_{\bbQ}(A^i(X)/N^i(X))\leqslant\ord_{s=m+i}\Phi_{2i}(s)
$$
with equality if and only if $\E^i(X,\ell)$ and $\T^i(X,\ell)$ hold.
\end{theorem*}
This was proved by Tate for the case that $K$ is a finite field (cf.\,\cite[Thm.\,2.9]{Tate_1994}). The general case was known to experts and its proof was sketched by P\'al \cite[Prop.\,6.9]{Pal_2022}. For completeness, we write down a proof of it following P\'al's proof of the rigid analogue. In the case that $i=1$, it is well-known that the Conjecture~\hyperlink{conj1}{$\T^1(X,\ell)$} is equivalent to the finiteness of the geometric Brauer groups of $X$ (cf.\,Proposition~\ref{prop2.1}).
\begin{corollary*}[\hypertarget{cor0}{C}]
 The following statements are equivalent: 
\begin{itemize}
\item[(1)] $\Br(X_{K^\sep})^{G_K}(\ell)$ is finite for some prime $\ell \neq p$.

\item[(2)] The prime-to-$p$ part of $\Br(X_{K^\sep})^{G_K}$ is finite. 

\item[(3)] Conjecture~\hyperlink{conj1}{$\T^1(X,\ell)$} holds for all primes $\ell \neq p$.

\item[(4)] Conjecture~\hyperlink{conj2}{2} holds for $i=1$, i.e.\ the rank of $\NS(X)$
is equal to $-\ord_{s=m+1}\Phi_{2}(s)$.
\end{itemize}
\end{corollary*}
\begin{remark*}
In the case that $K$ is finite, this was proved by Tate \cite{Tate_1994} and Lichtenbaum \cite{lichtenbaum1983zeta}.
The equivalence between (1), (3) and (4) is a consequence of the theorem above which is due to P\'al \cite[Prop.\,6.9]{Pal_2022}. The equivalence between (2) and (3) is due to Cadoret--Hui--Tamagawa~\cite{CadoretHuiTamagawa_2023}.  We will use this result to deduce the prime-to-$p$ part of our Theorem~\hyperlink{bsdconj}{A}. For this, we prove a relation between Brauer groups and Tate--Shafarevich groups.
\end{remark*}
Let $K$ be a finitely generated field. Its \emph{unramified Brauer group} $\Br_{\nr}(K/k)$ is defined to be the intersection of $\Br(\sO_{K,v})$ for all discrete valuation rings $\sO_{K,v}$ with function field $K$. For a torsion abelian group and a prime number $\ell$, define $V_\ell M=\lim \limits_{\substack{\leftarrow \\ n}}M[\ell^n]\otimes_{\mathbb{Z}_\ell}\mathbb{Q}_\ell.$
\begin{theorem*}[\hypertarget{bigtheorem}{D}]
Let $X$ be a smooth projective geometrically connected variety over a finitely generated field $K$. Then there exists a canonical exact sequence
$$
0\longrightarrow V_\ell\Br_{\nr}(K/k)\longrightarrow 
\ker\big(V_\ell\Br_{\nr}(K(X))\rightarrow V_\ell\Br(X_{K^{\sep}})^{G_K}\big) 
\longrightarrow V_\ell\Sha(\Pic^0_{X/  K,\red})\longrightarrow 0.
$$
Moreover, if $\ell \neq \Char(K)$, the natural map
$$\Br_{\nr}(K(X))(\ell)\rightarrow \Br(X_{K^{\sep}})^{G_K}(\ell)$$
has a finite cokernel and is surjective for sufficiently large $\ell$.
\end{theorem*}

In the case that $K$ is a global field, the above theorem was proved in \cite{qinbrauer0,QinBrauer2} by the third named author. We also obtain a generalisation of it relative to a finitely generated field.

\begin{theorem*}[\hypertarget{bigthm}{E}]
Let $\rho\colon\sX\rightarrow\sY$ be a dominant $k$-morphism between smooth geometrically connected varieties over a finitely generated field $k$. Let $K$ be the function field of $\sY$. Assume that the generic fibre $X$ of $\rho$ is smooth projective geometrically connected over $K$. Set $K^\prime\colonequals Kk^\sep$. Let $\ell\neq \Char(k)$ be a prime, and set
$$
\sK\colonequals\ker (V_\ell\Br(\sX_{k^\sep})^{G_k}\rightarrow V_\ell \Br(X_{K^\sep})^{G_K}) .
$$
Then we have canonical exact sequences
$$
0\rightarrow V_\ell\Br(\sY_{k^\sep})^{G_k}\rightarrow \sK \rightarrow V_\ell\Sha_{K^\prime}(\Pic^0_{X/K, \red})^{G_k}\rightarrow 0 ,
$$
$$0\rightarrow \sK\rightarrow V_\ell\Br(\sX_{k^\sep})^{G_k}\rightarrow V_\ell \Br(X_{K^\sep})^{G_K}\rightarrow 0, 
$$
where $V_\ell\Sha_{K^\prime}(\Pic^0_{X/K,\red})^{G_k}$ is defined in \cref{prop: sha} and only depends on $A/K$.
Moreover, let $B$ denotes the $K/k$-trace of $\Pic^0_{X/K, \red}$, then there is a canonical exact sequence
$$
0\rightarrow V_\ell\Sha(B) \rightarrow V_\ell\Sha(\Pic^0_{X/K, \red}) \rightarrow V_\ell\Sha_{K^\prime}(\Pic^0_{X/K, \red})^{G_k}\rightarrow 0 .
$$
In particular, if $k$ is finite, $V_\ell\Sha(\Pic^0_{X/K,\red})=V_\ell\Sha_{K^\prime}(\Pic^0_{X/K,\red})^{G_k}$.
\end{theorem*}
\begin{remark*}
The relation between the Tate--Shafarevich group and the Brauer group has been studied in the following cases:
\begin{itemize}
\item for fibrations of relative dimension $1$ over curves: by Artin and Grothendieck (cf.\,\cite[$\S$\,4]{Grothendieck_1968-III} or \cite[Prop.\,5.3]{Ulm}), Milne~\cite{Mil2}, Gonzales-Aviles \cite{Goa} and Geisser  \cite{Geisser_2021},
\item in the case of good reduction everywhere for relative dimension $1$ over higher-dimensional bases by the second author 
\cite[Thm.\,4.27]{Keller_2016},
\item for higher relative dimension over curves by Geisser~\cite{Geisser_2021} and the third author \cite{qinbrauer0}.
\end{itemize}
For fibrations over arbitrary finitely generated fields, Ulmer \cite[\S\,7.3]{Ulm} formulated some questions on relations between Tate conjecture and BSD conjecture in terms of $L$-functions. The previous theorem
may be considered as a partial answer to his question. The Tate--Shafarevich group $V_\ell\Sha_{K^\prime}(A)^{G_k}$ there
should be regarded as the obstruction of Nagao's conjecture \cite[Conj.\,1.1]{rosen1998rank} when $A$ is a non-isotrivial elliptic curve and Hindry--Pacheco's generalisation \cite[Thm.\,1.3]{hindry2005rang} of Nagao's conjecture for Jacobian of curves. It seems plausible to formulate Nagao's conjecture for arbitrary abelian varieties relative to a finitely generated field and verify that the conjecture in positive characteristic is equivalent to the vanishing $V_\ell\Sha_{K^\prime}(A)^{G_k}=0$.
\end{remark*}
\begin{corollary*}[\hypertarget{cor1intro}{F}]
Let $\sX\rightarrow\sY$  be as in the above theorem, and assume that $k$ is a finite field. Then  Conjecture~\hyperlink{BS3}{B+S-D+2} for  $\sX$  is equivalent to 
$$\emph{Conjecture} \ \hyperlink{BS3}{\text{B+S-D+2}}  \ \emph{for}\ \sY \ + \emph{Conjecture} \ \hyperlink{bsd}{\text{B+S-D}} \ \emph{for} \ \Pic^0_{X/K} + \ \emph{Conjecture} \ \hyperlink{conj2}{2} \ \emph{for}\  X  \ \emph{when}\ i=1.
$$
\end{corollary*}
\begin{remark*}
For the case that $\sY$ is a smooth projective curve and $\pi$ is proper, this result was firstly proved by Geisser (cf.\,\cite[Thm.\,1.1]{Geisser_2021}).
\end{remark*}

\subsection*{Verification in concrete cases}

	Using Theorem~\hyperlink{bsdconj}{A}, it is now possible to verify the Conjecture~\hyperlink{bsd}{B+S-D} in some concrete cases, under the assumption
	that $A$ is given by equations in some projective space or as the Jacobian of a curve.
	This then also implies that the Tate--Shafarevich group is of finite exponent (and similar results for the Brauer group), 
	which implies in particular that its prime-to-$p$ part is finite by~\cite[Thm.\,1.3]{QinBrauer1} and all the conjectures in~\cref{cor1}, \cref{cor: conjectures follow from Tate}, Theorem~\hyperlink{bsdconj}{A} and~\cref{thm: obstruction analytic Tate}. 
	Let $K$ be a finitely generated field of characteristic $p > 0$, and $A/K$ an abelian variety.
	To show Conjecture~\hyperlink{bsd}{B+S-D} for $A/K$, one has to verify the equality $r_\alg(A/K) = r_\an(A/K)$. 
	
	In the case at hand, one can compute $r_\an(A/K)$ via the étale cohomology of a model of $A$, 
	as this allows to compute the $L$-function and the analytic rank of the model as described in~\cite{MadoreOrgogozo_2015}. 
	(Note that here one takes advantage of the fact that the $L$-function is a rational function in $q^{-s}$.
	Over number fields, one cannot do this---the only tool available is the Gross--Zagier formula and the knowledge of the sign of the functional equation.  
	With this one can only show that $r_\an(A/K)$ is contained in a certain finite subset of natural numbers, more precisely $r_\an(A/K) \in \{0,1,2,3\}$.) 
	
	Next, recall that we already have the inequality $r_\alg(A/K) \leqslant r_\an(A/K)$ from~\cref{prop: rBSD inequality}.
	Thus, to prove equality, it suffices to construct $r_\an(A/K)$ independent points in $A(K) \otimes \bbQ$. 
	To do so, one enumerates points of bounded height in $A(K)$: if one has equations of $A/K$, one has a projective embedding of $A$ into some projective space $\bbP^n_K$. One can list all points in $\bbP^n(K)$ by increasing height (basically by increasing multidegree if $K = k(X_1,\ldots,X_d)$---in general, $K$ is a finite extension of a rational function field, and there is a finite-to-one map from $K$ to the function field) and test whether they lie in $A$.
	(For example, if $A = E$ is an elliptic curve given by a Weierstraß equation, or more generally if $A = \Jac(X)$ is the Jacobian of a curve $X/K$, one can easily enumerate divisors.
	One assumes that one can decide about linear equivalence of the divisors on $A$, i.e., equality in the Picard group.
	For example for hyperelliptic Jacobians this is possible using the Mumford representation, or in general using Riemann--Roch for the curve $X/K$.)
	To check linear independence of the points of bounded height, which one has obtained, one computes the rank of the height pairing matrix of pairs out of these points. 
	(Recall that the height pairing in positive characteristic is the intersection pairing and non-degenerate, see~\cite[Thm.\,9.15]{Conrad_2006} for the N\'eron--Tate height pairing.)

\subsection*{Outline of proofs}
In Section~\ref{sec: L-functions}, we first recall the definitions and basic properties of $L$-functions. Then we prove the equivalence of the above Conjectures for divisors, i.e., $i = 1$. For an abelian variety $A/K$, in~\cref{prop: rBSD inequality} we prove that its
algebraic rank is bounded above by 
the analytic rank, and that the latter equals the dimension of the generalized $1$-eigenspace of Frobenius on the $\ell$-adic Selmer group $\H^1(\sY_{k^\alg},V_\ell \sA)$.
The proofs are based on methods from étale cohomology, including Deligne's theory of weights \cite{Deligne_1980}, Poincaré duality and the trace formula.
 We also use Matsusaka's Theorem to move between algebraic, $\ell$-adic
 homological, and numerical equivalence for divisors.

In Section~\ref{subsec: obstruction}, we prove in Subsection~\ref{sec:Tate ell adic} the $\ell$-part of Corollary~\hyperlink{cor0}{C}, namely that the finiteness of the $\ell$-primary part of the Brauer group is the obstruction of the Tate Conjecture~\hyperlink{conj1}{$T^1(X,\ell)$ and~\hyperlink{BS3}{B+S-D+2},}
and an analogous result with the unramified Brauer group in Subsection~\ref{sec:Tate ell adic}. We use the cycle class map in \'etale, crystalline and rigid cohomology, and as a key property that it is Galois equivariantly split (see~\cref{ssplit} in the $\ell$-adic and~\cref{splitnesscor} in the $p$-adic case) To overcome the problem that we do not have smooth projective models in general, we use de Jong's alterations. In the $p$-adic case, lacking resolutions of singularities, one key idea is to use the unramified Brauer group because it only depends on the function field and equals the cohomological Brauer group for varieties with a smooth compactification. We give a proof of Theorem~\hyperlink{firstthm}{B} and its analogue for abelian varieties.

In Section~\ref{sec:Sha and Brauer}, we relate the $\ell$-adic Tate modules of the Tate--Shafarevich group to the Brauer group in~\cref{big} to get the prime-to-$p$ part of Theorem~\hyperlink{bigtheorem}{D} and Theorem~\hyperlink{bigthm}{E}, and analogously the $p$-adic Tate module of the unramified Brauer group in~\cref{thm: fibrations p-torsion}. For the latter, we need to construct suitable models of the base whose Brauer group does not differ too much from its unramified one. Again we use de Jong's alterations.

In the final two Sections~\ref{sec:elladic obstruction of BSD rank} and~\ref{sec:padic obstruction to BSD rank}, we prove Theorem~\hyperlink{bsdconj}{A} on the equivalence of the BSD rank conjecture and the finiteness (finite exponent for the $p$-part) of the Tate--Shafarevich group for an abelian variety $A/K$. Following Schneider's approach \cite{Schneider_1982}, we give a direct proof of the $\ell$-primary part of Theorem~\hyperlink{bsdconj}{A} without relating it to Tate conjecture for divisors. We also include one section on the N\'eron--Tate height pairing and Yoneda pairing. The $p$-adic case uses D'Addezio's theorem \cite{dAddezio}.

\subsection*{Notation and Terminology}

\subsubsection*{Fields}
By a \emph{finitely generated field}, we mean a field which is finitely generated over a prime field.\\
For any field $k$, denote by $k^\sep$ the separable closure. Denote by $G_k=\Gal(k^\sep/k)$ the absolute Galois group of $k$.

\subsubsection*{Henselisation}
Let $R$ be a Noetherian local ring, denote by $R^\h$ (resp. $R^{\sh}$) the Henselisation (resp.\ strict Henselisation) of $R$ at the maximal ideal. If $R$
is a domain, denote by $K^\h$ (resp. $K^{\sh}$) the fraction field of $R^{\h}$ (resp. $R^{\sh}$).

\subsubsection*{Schemes}
All schemes are assumed to be separated over their bases. For a Noetherian scheme $S$, denote by $|S|$ the set of closed points and by $S^1$ the set of points of codimension $1$.\\
By a \emph{variety} over a field $k$, we mean a scheme which is separated and of finite type over $k$.\\
For a smooth proper geometrically connected variety $X$ over a field $k$, we use $\Pic^0_{X/k}$ to denote the underlying reduced closed subscheme of the identity component of the Picard scheme $\underline{\Pic}_{X/k}$.

\subsubsection*{Cohomology}
The default sheaves and cohomology over schemes are with respect to the
small \'etale site. So $\H^i$ is usually the abbreviation of $\H_{\et}^i$ unless otherwise
stated. We use $\H^i_{\fppf}$ to denote flat cohomology for sheaves on fppf site.

\subsubsection*{Abelian groups}
For any abelian group $M$, integer $m$ and prime $\ell$, we set\\
$$M[m]=\{x\in M| mx=0\},\quad M_{\tor}=\bigcup\limits_{m\geqslant 1}M[m],\quad  M(\ell)=\bigcup\limits_{n\geqslant 1}M[\ell^n], $$
$$T_\ell M=\Hom_\mathbb{Z}(\mathbb{Q}_\ell/\mathbb{Z}_\ell, M)=\lim \limits_{\substack{\leftarrow \\ n}}M[\ell^n],\quad V_\ell M= T_\ell(M)\otimes_{\mathbb{Z}_\ell}\mathbb{Q}_\ell.$$

\subsection*{Acknowledgements} 

The last named author would like to thank Xinyi Yuan for important communications and Marco D'Addezio for pointing out that Theorem (B) was proved by P\'al \cite{Pal_2022}.

\section{$L$-functions, the BSD-rank conjecture and Tate conjectures}\label{sec: L-functions}

This section is devoted to the definition and basic properties of $L$-functions. 
In \cref{cor: nsrank}, we relate the Tate conjecture \hyperlink{conj1}{$\T^1(X/K,\ell)$} for divisors and $\ell$ invertible on the base to the pole order of certain $L$-functions. 
This appears also as part of a result due to Pál which we recall in  \cref{thm: equivalence Tate over function fields}. 
(However, as he only gave the prove for the $p$-adic case, we provide the argument in the $\ell$-adic for the benefit of the reader.)
In particular, Pál's result implies that the Conjecture~\hyperlink{conj1}{$\T^1(X/K,\ell)$} is independent of $\ell$. 
In~\cref{prop: rBSD inequality}, we prove that in our setting the algebraic rank of an abelian variety $A/K$ less than or equal to its analytic rank and relate the analytic rank to the dimension of the generalised $1$-eigenspace of Frobenius of the $\ell$-adic Selmer group of $A/K$.

\begin{notation} 
Let $X/K$ be a smooth projective variety. 
For a smooth model $\rho\colon\sX\rightarrow \sY$ of $X/K$, $n\in\bbN_0$
and a prime $\ell$ invertible on $\sY$, we will denote the lisse $\bbQ_\ell$-sheaves $\R^i\rho_\ast\bbQ_\ell$ by $\sR^i_\ell(\sX/\sY)$ and its $L$-function by $\Phi_i(\sX/\sY,s)$.
\end{notation}

We will use that the order of poles and zeros of $L$-functions at certain points is independent of the choice of a model. 
This is shown in the following proposition using Deligne's theory of weights \cite{Deligne_1980}.

\begin{proposition}\label{conv}
In the set-up of the introduction $\Phi_i(\sX/\sY,s)$ converges absolutely for $\Re(s)>\dim\sY+i/2$.
Moreover, the zeros and poles of $\Phi_i(\sX/\sY,s)$ in the strip
$$\dim\sY+i/2-1<\Re(s)\leqslant \dim\sY+i/2$$
are independent of the choice of a smooth projective model  $\sX\rightarrow \sY$.
\end{proposition}
\begin{proof}
By Deligne's theorem, the $P_{y,i}( T)$ are of fixed degree with reciprocal roots of absolute value $q_y^{i/2}$.
Therefore, by \cite[Thm.\,1]{Serre_1965}, $\Phi_i(\sX/\sY,s)$ converges absolutely for $\Re(s)>\dim\sY+i/2$.
If we replace $\sY$ by a non-empty open subscheme, and let $\sZ$  denote the complement, we have to divide $\Phi_i(\sX/\sY,s)$ by 
$$ 
\prod_{z\in |\sZ|} P_{z,i}(q_z^{-s}),
$$
which converges absolutely for $\Re(s)>\dim \sY+i/2-1$ since $\dim\sZ\leqslant \dim \sY -1$. It follows that the zeros and poles of $\Phi_i$ in the strip
$$\dim\sY+i/2-1<\Re(s)\leqslant \dim\sY+i/2$$
depend only on $X/K$ and not on the choice of $\sX/\sY$.
\end{proof}

The following propositions were sketched by Tate in~\cite{Tate_1965}.

\begin{proposition-definition}\label[proposition-definition]{bira}
Let $\sX$ be a regular scheme of finite type over $\bbZ$. Assume that its zeta function $\zeta(\sX,s)$ can be meromorphically continued to the point $s=\dim\sX-1$. Let $e(\sX)$ be the order of $\zeta(\sX,s)$ at this
point, and put
$$
z(\sX)=\rk \, \H^0(\sX,\sO_{\sX}^*)-\rk\, \H^1(\sX,\sO_{\sX}^*)-e(\sX).
$$
Then $z(\sX)$ is a birational invariant.
\end{proposition-definition}

\begin{proof}
	If one removes from $\sX$ a closed irreducible subscheme $\sZ$, then we have 
	$$\zeta(\sX-\sZ,s)=\zeta(\sX,s)/\zeta(\sZ,s).$$
	Since $\zeta(\sZ,s)$ converges absolutely for $\Re(s)>\dim \sZ$ and has a simple pole at $s=\dim \sZ$ (cf.\,\cite[Thm. 1, 2, 3]{Serre_1965}), the order of $\zeta(\sZ,s)$ at $\dim\sX-1$ is $0$ if the codimension of $\sZ$ is at least $2$ and is $1$ if the codimension of $\sZ$ is $1$. If $\sZ$ is a divisor, set $U=\sX-\sZ$, there is an exact sequence\\
	$$
	0\rightarrow \H^0(\sX,\sO_{\sX}^*)\rightarrow \H^0(U,\sO_{U}^*)\rightarrow \bbZ\rightarrow\Pic(\sX)\rightarrow\Pic(U)\rightarrow 0.
	$$
	It follows that $e(\sX)=e(U)-1$, so $z(\sX)=z(U)$.
\end{proof}

\begin{proposition}\label[proposition]{bira2}
Let $\rho\colon \sX\rightarrow\sY$ be a smooth projective morphism of schemes of finite type over $\bbZ$ whose generic fibre is a smooth projective geometrically connected variety $X$ over a finitely generated field $K$. 
Then any two of the following statements imply the third:
\begin{itemize}
\item[(i)] Conjecture~\hyperlink{bsd}{B+S-D} for $X/K$.
\item[(ii)] Conjecture~\hyperlink{conj2}{2},  for $i=1$, for $X/K$.
\item[(iii)] The equality $z(\sX)=z(\sY)$.
\end{itemize}
\end{proposition}

\begin{proof}
Since $\rho$ is proper and the generic fibre is geometrically irreducible and $\sY$ is normal, we have $\rho_*\sO_{\sX}=\sO_{\sY}$. Thus $\H^0(\sX,\sO_{\sX}^*)=\H^0(\sY,\sO_{\sY}^*)$. By \cite[Chap.\,III, Cor.\,11.3]{Hartshorne}, all fibres of $ \rho$ are geometrically connected. Let $D$ be a vertical prime divisor on $\sX$. Since $\rho$ is smooth proper and all fibres are connected, one can show that $D= \rho^{-1}( \rho(D))$ as a Weil divisor. Thus, any vertical divisor on $\sX$ is the pullback of some divisor on $\sY$. By the Leray spectral sequence
$$\H^p(\sY,\R^q \rho_*\bbG_m)\Rightarrow \H^{p+q}(\sX,\bbG_m),$$
and $ \rho_*\bbG_m=\bbG_m$, we get an injection
$$\H^1(\sY,\bbG_m)\rightarrow  \H^1(\sX,\bbG_m).$$
Thus, we get an exact sequence
$$ 0\rightarrow \Pic(\sY)\rightarrow \Pic(\sX)\rightarrow \Pic(X)\rightarrow  0.$$
Hence $z(\sX)=z(\sY)\Leftrightarrow e(\sX)-e(\sY)=-\rk \,\Pic(X)$. Since $\Phi_{2d}(\sX/\sY,s)=\zeta(\sY,s-d)$, it suffices to show that $\Phi_i(\sX/\sY,s)$ has no pole or zero at $s=\dim \sX-1$ for $i\neq 2d-1, 2d-2, 2d.$ Since all reciprocal roots of $P_{y,i}(T)$ have absolute values $q_y^{\frac{i}{2}}$, $\Phi_i(\sX/\sY,s)$ converges absolutely for $\Re(s)>\dim \sY+i/2$, the claim follows.
\end{proof}

The next proposition gives a bound for
the dimension of the group of cycles of codimension $i$ with respect to $\ell$-adic homological equivalence by the dimension of a $2i$-th \'etale cohomology group. 
And this in turn is bounded by the pole order of $\Phi_{2i}$.

\begin{proposition}\label[proposition]{nsrank}
Let $X$ be a smooth projective geometrically connected variety over a finitely generated field $K$ of characteristic $p>0$ with a smooth model $\rho\colon \sX\rightarrow \sY$.
For $\ell\neq p,$ 
	$$ 
	\dim_\bbQ (A^i(X)/N^i(X))\leqslant \dim_{\bbQ_\ell} 
	\H^{2i}(X_{  K^{\sep}},\bbQ_\ell(i))^{G_{  K}}\leqslant -\ord_{s=\dim(S)+i}\Phi_{2i}(\sX/\sY,s).
	$$
Moreover the equation $\dim_\bbQ (A^i(X)/N^i(X))=-\ord_{s=\dim(S)+i}\Phi_{2i}(\sX/\sY,s)$ holds if and only if Conjecture~\hyperlink{conj1}{$\T^i(X/K,\ell)$} and  Conjecture~\hyperlink{conj: E}{$\E^i(X,\ell)$}  is true for $\ell\neq p$. 
\end{proposition}

\begin{proof}
	The cycle classes of $A^{\dim(X)-i}(X)$ generate a $\bbQ_\ell$-subspace in $\H^{2\dim(X)-2i}(X_{  K^{\sep}},\bbQ_\ell(\dim(X)-i))$. 
	Let $b_1,\ldots,b_s \in A^{\dim(X)-i}(X)$ such that their cycle classes form a basis of this subspace. 
	Let $B \subseteq A^{\dim(X)-i}(X)$ be the subspace spanned by the $b_j$, $j = 1,\ldots,s$. 
	Then the intersection pairing
	$$A^i(X)/ N^i(X)\times B\longrightarrow \bbQ$$
	is left non-degenerate. 
	Thus, the induced map $A^i(X)/ N^i(X)\rightarrow B^*$ is injective. 
	It follows that $A^i(X)/ N^i(X)$ has finite dimension. 
	Let $a_1,..,a_t\in A^i(X)$ such that they form a basis in $A^i(X)/ N^i(X)$. 
	Then there exist $a^*_1,\ldots, a^*_t\in B$ such that $a_i.a^*_j=\delta_{ij}$. 
	Hence $\gamma^\ell_X(a_i).\gamma^\ell_X(a_j^*)=\delta_{ij}$. 
	Thus $\gamma^\ell_X(a_1),\ldots, \gamma^\ell_X(a_t)$ are $\bbQ_\ell$-linearly independent in $\H^{2i}(X_{  K^{\sep}},\bbQ_\ell(i))^{G_{  K}}$. 
	This proves the first inequality.
	
	Let $a\in  N^i(X)$ with $\gamma^\ell_X(a)\neq 0$, then $\gamma^\ell_X(a_1),\ldots, \gamma^\ell_X(a_t), \gamma^\ell_X(a)$ are $\bbQ_\ell$-linear independent. 
	Let $\lambda_i\in \bbQ_\ell$ such that 
	$$\lambda_1\gamma^\ell_X(a_1)+\ldots+\lambda_t\gamma^\ell_X(a_t)+\lambda_{t+1}\gamma^\ell_X(a)=0$$
	Intersecting with $a_j^*$, we get $\lambda_j=0$ for $1\leqslant j\leqslant t$. Since $\gamma^\ell_X(a)\neq 0$, so $\lambda_{t+1}=0$. 
	This implies that
	$$\dim_\bbQ (A^i(X)/ N^i(X))< \dim_{\bbQ_\ell} \H^{2i}(X_{  K^{\sep}},\bbQ_\ell(i))^{G_{  K}}$$
	if $ N^i(X)\neq 0$. 
	Assuming the second inequality, we proved the $``\Rightarrow"$-direction of the second claim.
	
	Next we will prove the second inequality. 
	Let  $K^\prime$ denote $  Kk^{\sep}$. 
	By Deligne's theorem, the lisse sheaf
	$\sR^{2\dim(X)-2i}_\ell(\sX/\sY)$ is pure of weight $2\dim(X)-2i$. 
	By Grothendieck's trace formula, we have 
	\begin{align*}
	\Phi_{2\dim(X)-2i}(\sX/\sY,s)&=L(\sY,\sR^{2\dim(X)-2i}_\ell(\sX/\sY),s) \\
	&= 
	\prod_{j=0}^{2\dim(S)}
	\det(1-q^{-s}\Frob^{-1}|_{\H^j_{\et,c}(\sY_{\bar{k}}, \sR^{2\dim(X)-2i}_\ell(\sX/\sY))})^{(-1)^{j+1}},
	\end{align*}
	where $\Frob$ denotes  the arithmetic Frobenius element in $G_k$. 
	Since $\H^j_{\et,c}(\sY_{\bar{k}}, \sR^{2\dim(X)-2i}_\ell(\sX/\sY))$ is mixed of weight $\leqslant j+2\dim(X)-2i$,   $-\ord_{s=\dim(X)+\dim(S)-i} \Phi_{2\dim(X)-2i}(\sX/\sY,s) $ is equal to the dimension of the generalised $q^{\dim(X)+\dim(S)-i}$-eigenspace of $\Frob^{-1}$ on $\H^{2\dim(S)}_{\et,c}(\sY_{\bar{k}}, \sR^{2\dim(X)-2i}(\sX/\sY))$. 
	By Poincar\'e duality, 
	$$
	\H^{2\dim(S)}_{\et,c}(\sY_{\bar{k}}, \sR^{2\dim(X)-2i}(\sX/\sY))\cong \H^0(\sY_{\bar{k}}, \sR^{2\dim(X)-2i}(\sX/\sY)^\vee)^\vee(-\dim(S)).
	$$ 
	Therefore $-\ord_{s=\dim(X)+\dim(S)-i} \Phi_{2\dim(X)-2i}(\sX/\sY,s) $ is equal to the dimension of the generalised $1$-eigenspace of  $ \H^0(\sY_{\bar{k}}, \sR^{2\dim(X)-2i}(\sX/\sY)^\vee)(n-\dim(X))$.  $\sR^{2\dim(X)-2n}(\sX/\sY)^\vee$ corresponds to the $G_{  K}$-representation on $\H^{2d-2i}(X_{  K^{\sep}},\bbQ_\ell)^\vee$ which is isomorphic to $\H^{2i}(X_{  K^{\sep}},\bbQ_\ell)(\dim(X))$. 
	
	It follows that $-\ord_{s=\dim(X)+\dim(S)-i} \Phi_{2\dim(X)-2i}(\sX/\sY,s) $ is equal to the dimension of the $1$-eigenspace of  $\H^{2i}(X_{  K^{\sep}},\bbQ_\ell(i))^{G_{K^\prime}}$. 
	By Poincar\'e duality, we have 
	$$\Phi_{2\dim(X)-2i}(\sX/\sY,s)=\Phi_{2i}(s-\dim(X)+2i).$$ 
	Thus  $-\ord_{s=\dim(S)+i} \Phi_{2i}(\sX/\sY,s)$ is equal to the dimension of the generalised $1$-eigenspace of $\H^{2i}(X_{  K^{\sep}},\bbQ_\ell(i))^{G_{K^\prime}}$ which contains $\H^{2i}(X_{  K^{\sep}},\bbQ_\ell(i))^{G_{  K}}$. This proves the second inequality.
	
	For the ``$\Leftarrow$''-direction of the second statement, assuming $\T^i(X/  K,\ell)$ and $\E^i(X/  K,\ell)$,  
	we have 
	$$A^i(X)_{\bbQ_\ell}=\H^{2i}(X_{  K^{\sep}},\bbQ_\ell(i))^{G_{  K}}.$$
	Consider the $G_k$-equivariant pairing 
	\begin{equation}\label{directsum}
		\H^{2i}(X_{  K^{\sep}},\bbQ_\ell(i))^{G_{K^\prime}}\times A^{\dim(X)-i}(X)_{\bbQ_\ell}\longrightarrow \bbQ_\ell.
	\end{equation}
	$\E^i(X/  K,\ell)$ implies that 
	$$A^i(X)_{\bbQ_\ell}\times A^{\dim(X)-i}(X)_{\bbQ_\ell} \longrightarrow \bbQ_\ell$$
	is left nondegenerate. 
     Thus, the restriction of (\ref{directsum}) 
	$$\H^{2i}(X_{  K^{\sep}},\bbQ_\ell(i))^{G_  K}\times A^{\dim(X)-i}(X)_{\bbQ_\ell}\longrightarrow \bbQ_\ell$$
	is also left non-degenerate. 
	By the subsequent Lemma \ref{splitrep},  
	$\H^{2i}(X_{  K^{\sep}},\bbQ_\ell(i))^{G_  K}$ is a direct summand of $\H^{2i}(X_{  K^{\sep}},\bbQ_\ell(i))^{G_{K^\prime}}$ 
	as a $G_k$-representation. 
	This implies that the $1$-eigenspace of $\Frob^{-1}$ on $\H^{2i}(X_{  K^{\sep}},\bbQ_\ell(i))^{G_{K^\prime}}$ is a direct summand, 
	and therefore, is equal to the generalised $1$-eigenspace of $\Frob^{-1}$. 
	It follows that
	$$\dim_{\bbQ_\ell} \H^{2i}(X_{  K^{\sep}},\bbQ_\ell(i))^{G_  K}= -\ord_{s=\dim(S)+i} \Phi_{2i}(\sX/\sY,s).$$
	This completes the proof.
\end{proof}

\begin{remark}
In the above proof, we used the following well known equalities. 
Let $X/K$ be a smooth projective geometrically connected variety with a smooth model $\sX/\sY$.
Note that by Poincaré duality, we have for any $n\in\bbN_0$
$$
\Phi_{2\dim(X)-2n}(\sX/\sY,s)=\Phi_{2n}(\sX/\sY,s-\dim(X)+2n).
$$
So 
\begin{align*}
\ord_{s=\dim(\sY)+1}\Phi_2(\sX/\sY,s) &= \ord_{s=\dim(\sY)+\dim(X)-1} \Phi_{2\dim(X)-2}(\sX/\sY,s)\\
\ord_{s=\dim(\sY)+n}\Phi_{2n}(\sX/\sY,s) &= \ord_{s=\dim(\sY)+\dim(X)-n} \Phi_{2\dim(X)-2n}(\sX/\sY,s).
\end{align*}
\end{remark}

\begin{lemma}\label{splitrep}
	Let $k$ be a field and $V$ be a finite-dimensional $\bbQ_\ell$-linear continuous representation of $G_k$. Let $W\subseteq V$ be a sub-representation. Assuming that there exists a $G_k$-equivariant pairing 
	$$V\times V^\prime \longrightarrow \bbQ_\ell$$
	where $V^\prime$ is a $\bbQ_\ell$-linear representation of a finite quotient of
	$G_k$ such that the restriction of the pairing to $W\times V^\prime$ is left non-degenerate, then $W$ is a direct summand of $V$ as a $G_k$-representation.
\end{lemma}

\begin{proof}
	Since the pairing $W\times V^\prime$ is left non-degenerate, it induces a surjection
	$$V^\prime\longrightarrow W^*,$$
	where $W^*$ is the dual representation of $W$. Since the $G_k$-action on $V^\prime$ factors through a finite quotient, so $V^\prime$ is semisimple. Thus, there exists a subrepresentation $W^\prime$ of $V^\prime$ such that the pairing 
	$$W\times W^\prime$$
	is perfect. Define
	$$(W^\prime)^\bot\colonequals \{v\in V| (v, W^\prime)=0\}.$$
	It follows that $V=W\oplus(W^\prime)^\bot$. This proves the claim.
\end{proof}

\begin{corollary}\label[corollary]{cor: nsrank}
Let $X$ be a smooth projective geometrically irreducible variety over $K$. 
Then Conjecture~\hyperlink{conj1}{$\T^1(X/K,\ell)$} is equivalent to the Conjecture~\hyperlink{conj2}{2} for $i=1$.  
\end{corollary}
\begin{proof}
This follows immediately from \cref{nsrank}
using 
the fact that Conjecture~\hyperlink{conj: E}{$\E^1(X,\ell)$} is true by Matsusaka's Theorem~\cite[Thm.\,4]{Matsusaka_1957}
(which can be shown to hold for general fields).
\end{proof}

More generally, Pál obtained the following theorem:
\begin{theorem}[\cite{Pal_2022}]\label[theorem]{thm: equivalence Tate over function fields}
Let $X$ be smooth projective geometrically connected variety over $  K$ of characteristic $p>0$.
Then the following statements are equivalent:
\begin{enumerate}
\item Conjecture~\hyperlink{conj1}{$\T^1(X/K,\ell)$} is true for some prime number $\ell$.
\item Conjecture~\hyperlink{conj1}{$\T^1(X/K,\ell)$} is true for all prime numbers $\ell$. 
\item Conjecture~\hyperlink{conj2}{2} for $i=1$ holds.
\end{enumerate}
\end{theorem}

\begin{example}
As an application, we compute the Picard number of an abelian variety or a $K3$ surface. Let $K$ be a finitely generated field of characteristic $p>0$.
We consider the following two situations:
\begin{enumerate}
\item $X$ is an abelian variety over $  K$.
\item $k$ is of characteristic $p>2$ and $X$ is a $K3$-surface over $K$.
\end{enumerate}
Then 
$$
\rk_{\bbZ}(\NS(X))= -\ord_{s=\dim(S)+1}\Phi_2(\sX/\sY,s)
$$
for any smooth model $\sX/\sY$ of $X/K$.

Indeed, $\T^1(X/  K,\ell)$ holds for abelian varieties over function 
fields of positive characteristic \cite{Zarhin_1974a,Zarhin_1974b} and for $K3$-surfaces over function fields of characteristic $p>2$ \cite{Madapusi_2015}.
Thus the formula follows from \cref{thm: equivalence Tate over function fields} above.
\end{example}

For an abelian variety $A$ over a finitely generated field $K$ of characteristic $p>0$. Denote the transcendence degree of $K$ by $m$.

\noindent
\begin{proposition}\label[proposition]{prop: rBSD inequality}
Let $K$ be a finitely generated field over a finite field $k$. 
Let $A/K$ be an abelian variety and 
assume that  it extends to an abelian scheme $\rho\colon \sA\rightarrow\sY$ over a smooth integral variety $\sY$ over $k$.
\begin{enumerate}
\item The inequality  $\rk_{\alg}(A/ K) \leqslant \rk_{\an}(A/  K)$ 
for the algebraic and analytic rank holds.
\item The analytic rank $\rk_{\an}(A/  K)$ is equal to the dimension of 
the generalised $1$-eigenspace of the Frobenius action on 
$\H^1(\sY_{k^{\alg}}, V_\ell\sA)$. 
\end{enumerate}
\end{proposition}

\begin{proof}
By definition,
\begin{align*}
\Phi_1(\sA/\sY,s) =& L(\sY,\sR^1_{\ell}(\sA/\sY),q^{-s})\\
=& \prod_{j=0}^{2\dim(\sY)}
	\det(1-q^{-s}\Frob^{-1}|_{\H^j_{\et,c}(\sY_{\bar{k}}, \sR^{1}_\ell(\sA/\sY))})^{(-1)^{j+1}},
\end{align*}
where $\Frob$ is the arithmetic Frobenius in $G_k$. 
By Deligne's theorem, $\sR^1_{\ell}(\sA/\sY)$ is pure of weight $1$ and hence 
$\H^j_{\et,c}(\sY_{\bar{k}}, \sR^{1}_\ell(\sA/\sY))$ 
is mixed of weight $\leqslant j+1$, 
and, in particular, $H^{2\dim(\sY)}_{\et,c}(\sY_{\bar{k}}, \sR^{1}_\ell(\sA/ \sY))$ is pure of weight $2m+1$ by Poincaré duality. 
Therefore, only the factor contributing to $\rk_{\an}(A/ K)=\ord_{s=\dim(\sY)} \Phi_1(\sA/ \sY,s)$ is $\det(1-q^{-s}\Frob^{-1}|_{\H^{2m-1}_{\et,c}( \sY_{\bar{k}}, \sR^{1}_\ell(\sA/ \sY))})$. 
By Poincaré duality,
$$
\H^{2m-1}_{\et,c}( \sY_{\bar{k}}, \sR^{1}_\ell(\sA/ \sY)) \cong 
\H^1_{\et,c}( \sY_{\bar{k}}, \sR^{1}_\ell(\sA/ \sY)^{\vee})^{\vee}(-\dim(\sY)),
$$
and hence $\ord_{s=\dim(\sY)} \Phi_1(\sA/ \sY,s)$ equals to the dimension of the 
generalised $1$-eigenspace of 
$\Frob^{-1}$ on $\H^1_{\et,c}( \sY_{\bar{k}}, \sR^{1}_\ell(\sA/ \sY)^{\vee})$ 
and therefore
\begin{equation}\label{equ: ineq1}
\dim(\H^1_{\et,c}( \sY_{\bar{k}}, \sR^{1}_\ell(\sA/ \sY)^{\vee})^{G_k}) \leqslant \ord_{s=\dim(\sY)} \Phi_1(\sA/ \sY,s).
\end{equation}
As the lisse sheaf $\sR^{1}_\ell(\sA/ \sY)$ corresponds to the representation of $G_k$ on $\H^1(A_{\overline{ K}},\bbQ_{\ell})$, $ \sR^{1}_\ell(\sA/ \sY)^{\vee}$ corresponds to the $G_k$-representation $\H^1(A_{\overline{ K}},\bbQ_{\ell})^\vee$, 
which can be identified with the rational Tate module $V_\ell A$. 
Moreover, since by assumption $\sA$ is an abelian scheme over $ \sY$, 
there is for any $n$ a Kummer exact sequence of étale sheaves on $ \sY$
$$
0 \longrightarrow \sA[\ell^n] \longrightarrow \sA \xrightarrow{\cdot\ell^n} \sA \longrightarrow 0,
$$
where the finite étale $ \sY$-group scheme $\sA[\ell^n]$ represents a locally constant sheaf on $ \sY$. 
The $\bbQ_\ell$-sheaf  $V_\ell\sA$ associated to the inverse system of $\sA[\ell^n]$ corresponds to the $G_k$-representation $V_\ell A$, 
and thus can be identified with $ \sR^{1}_\ell(\sA/ \sY)^{\vee}$. 
This shows the second statement. 

The Kummer exact sequence induces a short exact sequence of cohomology groups
$$
0\longrightarrow \H^0( \sY_{\bar{k}}, \sA)\otimes_{\bbZ} \bbQ_\ell 
\longrightarrow \H^1( \sY_{\bar{k}}, V_\ell\sA)
$$
According to \cite[Thm.\,3.3]{Keller_2016}, $\sA$ satisfies the Néron mapping property
$$
\sA \cong j_\ast j^\ast \sA
$$
with the generic point $j\colon \Spec( K)\rightarrow   \sY$. 
Writing $K'=  K\bar{k}$, 
we can identify
$$
\H^0( \sY_{\bar{k}}, \sA)= A(K'),
$$
and therefore $\H^0( \sY_{\bar{k}}, \sA)^{G_k}=A( K)$.
Thus taking $G_k$-invariants of the above short exact sequence, 
we obtain
$$
0\longrightarrow A(  K) \otimes_{\bbZ} \bbQ_\ell 
\longrightarrow \H^1( \sY_{\bar{k}}, V_\ell\sA)^{G_k}
$$
and hence the inequality
$$
\rk\, A(K) \leqslant \dim (\H^1( \sY_{\bar{k}}, V_\ell\sA)^{G_k}).
$$
Together with \eqref{equ: ineq1}, this implies the first statement. 
\end{proof}


\section{The obstruction of the Tate Conjecture}\label{subsec: obstruction}

In this section we prove that the vanishing of the $\ell$-adic (respectively $p$-adic) 
Tate module of the Brauer group is the obstruction of the analytic Tate Conjecture~\hyperlink{BS3}{B+S-D+2} for divisors. This is proved in~\cref{mlem} and~\cref{thm: obstruction analytic Tate}, respectively.

\subsection{Definition of the Brauer group and some basic properties}

\begin{definition} 
Let $X$ be a Noetherian scheme.
\begin{enumerate}
\item
The \textit{cohomological Brauer group} of $X$ is given by
$$
\Br(X)\colonequals \H^2(X,\sO^\times_{X})_{\tor}.
$$

\item Let $K/k$ be a finitely generated field extension. The \emph{unramified Brauer group} of $K/k$ is defined as the subgroup
$$
\Br_{\nr}(K/k)\colonequals \bigcap_A \im\big(\Br(A)\rightarrow  \Br(K)\big) \subset \Br(K),
$$
where $A$ runs through the discrete valuation rings containing $k$ with field of fractions $K$. (See~\cite[Def.\,6.2.1]{ColliotTheleneSkorobogatov_2021}.)

\item Let $k$ be a finitely generated field. The \emph{unramified Brauer group} of $k$ is defined as the subgroup
$$
\Br_{\nr}(k)\colonequals \bigcap_A \im\big(\Br(A)\rightarrow  \Br(k)\big) \subset \Br(k),
$$
where $A$ runs through the valuation rings of all discrete valuations of $k$.

\item
Let $K/k$ be a finitely generated field extension. 
We also use the notation
$$
\Br^{\nr}(K/k)\colonequals \bigcap_{\nu\in div\Val(K)} \im\big(\Br(A_\nu)\rightarrow  \Br(K)\big) \subset \Br(K),
$$
where $\nu$ runs through the divisorial valuations with fraction field $K$ such that the associated discrete valuation ring $A_\nu$ contains $k$.
\end{enumerate}
\end{definition}

We collect some facts on the behaviour 
of the unramified Brauer group under field extensions.

\begin{lemma}\label[lemma]{lem: unram Br}
Let $k\subset K\subset L$ be extensions of fields. 
\begin{enumerate}
\item
The unramified Brauer group is functorial in the sense that there is a restriction map $\Br_{\nr}(K/k)\rightarrow \Br_{\nr}(L/k)$ compatible with field extensions.
\item 
If $L/K$ is finite separable, then there is also a corestriction 
$\Br_{\nr}(L/k)\rightarrow \Br_{\nr}(K/k)$.
\item
Assume that $K$ is finitely generated over a finite field $k=\bbF_{p^r}$ of characteristic $p>0$, and that $L/K$ is finite and purely inseparable.
Then $\Br_{\nr}(K/k)(p)$ is of finite exponent if and only if $\Br_{\nr}(L/k)(p)$ is of finite exponent.
\end{enumerate}
\end{lemma}

\begin{proof}
The first two items are 
\cite[Props.\,6.2.3 and 6.2.4]{ColliotTheleneSkorobogatov_2021}. 

Let us prove the third point. 
Assume first that $L/K$ is finite. 
Therefore, the restriction map $\Br_{\nr}(K/k)\rightarrow \Br_{\nr}(L/k)$ has a kernel killed by $[L:K]$. 
Hence, if $\Br_\nr(L/k)(p)$ is of finite exponent, 
$\Br_{\nr}(K/k)(p)$ is as well. 

For the reverse implication: Since $L/K$ is purely inseparably, 
there is an $n$ such that $\Frob^n(L) \subseteq K$, 
where $\Frob$ denotes the arithmetic Frobenius of $G_k$.  
Since $\Frob^n$ is the identity on $k$ (in other words a $k$-morphism), 
we have $\Br_\nr(\Frob^n(L)/k)(p) = \Br_\nr(L/k)(p)$. 
But since $K/\Frob^n(L)$ is finite, 
we see, by the same argument as above, that if $\Br_{\nr}(K/k)(p)$ is of finite exponent, then 
$\Br_{\nr}(\Frob^n(L)/k)(p) = \Br_\nr(L/k)(p)$ so, too. 
\end{proof}

\begin{remark}\label[remark]{rem: brauer}
The unramified and the usual Brauer groups are related by the following facts:
\begin{enumerate}
\item
For a regular, proper variety $X/k$ with function field $  K(X)$, 
the natural inclusion $\Br(X)\subseteq \Br(  K(X))$ induces an isomorphism
$
\Br(X) \xrightarrow{\sim} \Br_{\nr}(K(X)/k)
$
\cite[Prop.\,6.2.7]{ColliotTheleneSkorobogatov_2021}. 
\item
Let $X/k$ be a smooth variety with a smooth compactification $X \subseteq\overline{X}/k$. 
Then $ \Br_{\nr}(  K(X)/k) = \Br_{\nr}(K(\overline{X})/k)\cong \Br(\overline{X})$. 
\item
For a regular variety $X/k$ with function field $  K(X)$, 
there is a natural inclusion $\Br_{\nr}(K(X)/k)\subseteq \Br(X)$.
Indeed, we have
\begin{align*}
\Br_{\nr}(K(X){/k})= &
\bigcap_{v\in disc\Val(  K(X))} \im\big(\Br(\sO_{X,v}) \rightarrow  \Br( K(X))\big)\\
&\subseteq \bigcap_{x\in X^1} \im\big(\Br(\sO_{X,x})\rightarrow  \Br( K(X))\big) = \Br(X)
\end{align*}
where the last equality holds because for for regular schemes, the Brauer classes are those unramified at all codimension $1$ places.
\end{enumerate}
\end{remark}

\subsection{The obstruction of the Tate Conjecture in the $\ell$-adic case} \label{sec:Tate ell adic}

It will be important to understand the obstruction of the Tate conjecture in different situations.
We are interested in statements over finitely generated fields $K$ over a finite field $k=\bbF_q$, $q=p^r$, (for example, $K$ could be the function field of a smooth projective geometrically irreducible variety $S/k$).

\begin{proposition}\label[proposition]{prop2.1}
Let $X$ be a smooth projective geometrically connected variety over a field $K$  
and $\ell\neq \mathrm{char}(K)$ a prime number.
There is a canonical exact sequence
$$
0\longrightarrow \NS(X)_{\bbQ_\ell} \longrightarrow \H^2(X_{  K^{\sep}},\mathbb{Q}_\ell(1))^{G_{  K}} \longrightarrow V_\ell \Br(X_{  K^{\sep}})^{G_{  K}}\longrightarrow 0.
$$
\end{proposition}
\begin{proof}
This is well-known to experts (cf. \cite[Thm.\,1.2]{Yua2}).
\end{proof}

Tate showed in  \cite[Thm.\,5.2.b)]{Tate_1994} that $\T^1(X/ K,\ell)$, 
$\ell\neq \mathrm{char}(K)$,
is birationally invariant. His proof also gives the following:
\begin{proposition}\label[proposition]{bir}
	Let $X$ be a smooth  geometrically connected variety over  a finitely generated field $K$ of characteristic $p>0$, 
	and $U\subseteq X$ an open dense subvariety. Let $\ell \neq p$ be prime.
	Then the induced map
	$$
	V_\ell\Br(X_{  K^{\sep}})^{G_{  K}}\longrightarrow V_\ell\Br(U_{  K^{\sep}})^{G_{  K}}
	$$
	is an isomorphism.
\end{proposition}

\begin{proof}
	This follows from~\cite[Prop.~2.1]{qinbrauer0}.
\end{proof}

The following result generalises \cref{prop2.1} to noncompact varieties $X$ over finite fields $k$. 
We will use it to link the finiteness of 
$\Br(X_{\bar{k}})^{G_{k}}(\ell)$ 
to the analytic Tate Conjecture~\hyperlink{BS3}{B+S-D+2}.

\begin{proposition}\label[proposition]{ssplit}
Let $X/k$ be a smooth geometrically connected variety over a finite field of characteristic $p>0$ and $\ell\neq p$ a prime. 
Then there is a split exact sequence of $G_k$-modules
$$0\longrightarrow \Pic(X_{\bar{k}})_{\bbQ_\ell} \longrightarrow \H^2(X_{\bar{k}},\bbQ_\ell(1)) \longrightarrow V_\ell\Br(X_{\bar{k}})\longrightarrow 0.$$
\end{proposition}

\begin{proof}
Firstly, if the claim holds for $G_{k_1}$-actions where $k_1/k$ is a finite Galois extension, then it also holds for $G_k$-actions. 
Namely, let $V$ be a sub-representation of some $G_k$-representation $W$ over $\bbQ_\ell$. 
Let $p\colon W\rightarrow V$ be a $G_{k_1}$-equivariant projection. 
Then
$$\frac{1}{|G_{k}/G_{k_1}|}\sum\limits_{\sigma \in G_{k}/G_{k_1}}p^\sigma$$ 
is a $G_{k}$-equivariant projection.
	
Secondly, we assume that $X$ has a smooth projective compactification $X\subseteq \overline{X}$ over $k$, 
and show that the claim holds in this case. 
There is a commutative diagram with exact rows
\begin{displaymath}
		\xymatrix{
			0\ar[r] &\Pic(\overline{X}_{\bar{k}}) \otimes_\bbZ \bbQ_\ell \ar[r] \ar[d]^a& \H^2(\overline{X}_{\bar{k}},\bbQ_\ell(1))\ar[d]^b 
			\ar[r] & V_\ell\Br(\overline{X}_{\bar{k}}) \ar[d]^c \ar[r] & 0
			\\
			0\ar[r] & \Pic(X_{\bar{k}}) \widehat{\otimes}_\bbZ{\bbQ_\ell} \ar[r] & \H^2(X_{\bar{k}},\bbQ_\ell(1)) 
			\ar[r] & V_\ell\Br(X_{\bar{k}}) \ar[r] & 0
		}
\end{displaymath}
where the last vertical map $c$ is injective, hence $\ker(b)\cong\ker(a)$,
and by the theorem of cohomological purity \cite{Fujiwara_2002}, 
$\ker(b)\cong\ker(a)$ is generated by cycle classes of divisors contained 
in $\overline{X}\setminus X$. (Note that $\Pic(\overline{X}_{\bar{k}}) \otimes_\bbZ \bbQ_\ell = \Pic(\overline{X}_{\bar{k}}) \widehat{\otimes}_\bbZ \bbQ_\ell \colonequals \varprojlim_{n}\Pic(\overline{X}_{\bar{k}})/\ell^n \otimes_{\bbZ_\ell} \bbQ_\ell$ because $\Pic^0(\overline{X}_{\bar{k}})$ is $\ell$-divisible.)
Therefore, the fact that the first row splits together with the snake lemma implies that
$$
0\longrightarrow \Pic(X_{\bar{k}})\widehat{\otimes}_\bbZ \bbQ_\ell \longrightarrow\im(b) 
$$
is a split monomorphism.
Together with the above, it follows that $\Pic(X_{\bar{k}})$ is finitely generated modulo torsion, and we have in fact 
a split monomorphism
$$
0\longrightarrow \Pic(X_{\bar{k}})_{\bbQ_\ell} \longrightarrow\im(b) .
$$
By extending $k$, we may assume  that the $G_k$-action on 
$\Pic(X_{\bar{k}})_{\bbQ_\ell}$ is trivial. 
Then $\Pic(X_{\bar{k}})_{\bbQ_\ell}$ is contained in the 
generalised $1$-eigenspace of the Frobenius action on 
$\H^2(X_{\bar{k}},\bbQ_\ell(1))$. 
This eigenspace is a direct summand and is also contained in $\im(b)$ by the proof of \cref{bir}. 
It follows that $\Pic(X_{\bar{k}})_{\bbQ_\ell}$ is a direct summand of the generalised $1$-eigenspace, 
and therefore, it is also a direct summand of 
$\H^2(X_{\bar{k}},\bbQ_\ell(1))$. 	
	
Thirdly, by de Jong's theorem, 
there exists an alteration $f\colon X'\rightarrow X$ 
such that $X'$ admits a smooth projective compactification. 
It induces a commutative diagram 
$$
\xymatrix{
& \Pic(X_{\bar{k}})_{\bbQ_\ell} \ar[r] \ar[d]^{f^\ast} 
& \H^2(X_{\bar{k}},\bbQ_\ell(1)) \ar[r] \ar[d]^{f^\ast} 
& V_\ell\Br(X_{\bar{k}}) \ar[r]\ar[d]^{f^\ast} & 0 \\
0 \ar[r] & \Pic(X'_{\bar{k}})_{\bbQ_\ell} \ar[r] 
& \H^2(X'_{\bar{k}},\bbQ_\ell(1)) \ar[r] 
& V_\ell\Br(X'_{\bar{k}}) \ar[r] & 0
}
$$
with the lower roe being split exact by the second reduction step.
Since $f$ is proper and generically finite, 
there are pushforward maps
\begin{align*}
f_*\colon \H^2(X'_{\bar{k}},\bbQ_\ell(1))\longrightarrow \H^2(X_{\bar{k}},\bbQ_\ell(1))
&&\text{and}&& 
f_\ast\colon \Pic(X'_{\bar{k}}) \longrightarrow \Pic(X_{\bar{k}})
\end{align*}
such that $f_*f^*=\deg(f)$. It follows that 
the first and second vertical map in the above diagram are injective, 
and so is the map 
$\Pic(X_{\bar{k}})_{\bbQ_\ell} \rightarrow \H^2(X_{\bar{k}},\bbQ_\ell(1))$ in the first row. 
Through $f^*$, 
$\Pic(X_{\bar{k}})_{\bbQ_\ell}$ is a direct summand of 
$\Pic(X'_{\bar{k}})_{\bbQ_\ell}$, 
and therefore, is also a direct summand of $\H^2(X'_{\bar{k}},\bbQ_\ell(1))$. 
It follows that $\Pic(X_{\bar{k}})_{\bbQ_\ell}$ is a direct summand of $\H^2(X_{\bar{k}},\bbQ_\ell(1))$ and this finishes the proof.
\end{proof}

We can now prove the first main result of this section.
\begin{theorem}\label[theorem]{mlem}
Let $X/k$ be a smooth and irreducible  variety over a finite field of characteristic $p>0$, and $\ell\neq p$ prime. Then the analytic Tate Conjecture~\hyperlink{BS3}{B+S-D+2} for $X$ is equivalent to the finiteness of $\Br(X_{k^{\sep}})^{G_k}(\ell)$.
\end{theorem}
\begin{proof}
Using Grothendieck's trace formula, we can write the $\zeta$-function of $X$ as
$$
	\zeta(X,s)=\prod_{j=0}^{2\dim(X)}\det(1-q^{-s}\Frob^{-1} |_{\H^{j}_{\et,c}(X_{\bar{k}},\bbQ_\ell)})^{(-1)^{j+1}},
	$$
where $\Frob$ denotes the arithmetic Frobenius element in $G_k$. 
By Deligne's theorem, $\H^{j}_{\et,c}(X_{\bar{k}},\bbQ_\ell)$
is mixed of weight 
$\leq j$
and $\H^{2{\dim(X)}}_{\et,c}(X_{\bar{k}},\bbQ_\ell)$
is pure of weight $2\dim(X)$. 
Therefore, only the factors 
$\det(1-q^{-s}\Frob^{-1} |_{\H^{2\dim(X)-2}_{\et,c}(X_{\bar{k}},\bbQ_\ell)})$ and 
$\det(1-q^{-s}\Frob^{-1} |_{\H^{2\dim(X)-1}_{\et,c}(X_{\bar{k}},\bbQ_\ell)})$ 
contribute to the pole order of $\zeta(X,s)$ at $s=\dim(X)-1$, 
meaning that we have to study the generalised $q^{{\dim(X)}-1}$-eigenspaces of $\Frob^{-1}$ on 
$\H^{2\dim(X)-2}_{\et,c}(X_{\bar{k}},\bbQ_\ell)$ and 
$\H^{2\dim(X)-1}_{\et,c}(X_{\bar{k}},\bbQ_\ell)$. 
For this, we will make use of Poincaré duality, which tells us that 
\begin{align*}
\H^{2\dim(X)-1}_{\et,c}(X_{\bar{k}},\bbQ_\ell)\cong \H^1(X_{\bar{k}},\bbQ_\ell)^\vee(-{\dim(X)}),\\
\H^{2\dim(X)-2}_{\et,c}(X_{\bar{k}},\bbQ_\ell)\cong \H^2(X_{\bar{k}},\bbQ_\ell)^\vee(-{\dim(X)}).
\end{align*}

First, we show that the dimension of the generalised $q^{\dim(X)-1}$-eigenspace 
of $\Frob^{-1}$ on
$\H^1(X_{\bar{k}},\bbQ_\ell)^\vee(-{\dim(X)})$ is equal to 
$\rk\, \bbG_m(X)=\rk\, \H^0(X,\cO_{X}^{*})$.
By Kummer theory, we obtain an exact sequence 
$$
0\longrightarrow \H^0(X_{\bar{k}},\bbG_m)\otimes_\bbZ \bbQ_\ell 
\longrightarrow \H^1(X_{\bar{k}},\bbQ_\ell(1)) 
\longrightarrow V_\ell \H^1(X_{\bar{k}},\bbG_m).
$$
By de Jong's alteration theorem, there exists an alteration 
$f\colon X^\prime \rightarrow X$ such that $X^\prime$ admits a smooth projective compactification $X'\subset\overline{X}'$ over $k$. 
Since $z(X)$ from~\eqref{eq:def z} is a birational invariant, we might shrink $X$ 
to an open dense subset such that $f$ is finite flat. 
Due to the fact that $f$ is proper and generically finite, 
there is a canonical norm map $f_*\bbG_m \rightarrow  \bbG_m$ 
whose composition with the canonical map $\bbG_m\rightarrow  f_*\bbG_m$ 
is equal to the multiplication by $\deg(f)$. 
Since $\H^j(X^\prime,\bbG_m)=\H^j(X, f_*\bbG_m)$ because $f$ is finite, 
it follows that the composition
$$ \H^j(X,\bbG_m)\longrightarrow \H^j(X^\prime,\bbG_m)=\H^j(X, f_*\bbG_m) \longrightarrow \H^j(X,\bbG_m)
$$
is equal to the multiplication by $\deg(f)$.
Taking $j=1$ we deduce that 
$$
V_\ell \H^1(X_{\bar{k}},\bbG_m)\longrightarrow V_\ell \H^1(X^\prime_{\bar{k}},\bbG_m)
$$
is injective.
Now we take advantage of the fact that $X'$ has a smooth compactification $X'\subset\overline{X}'$, and hence the natural map
$$
\Pic(\overline{X}'_{\bar{k}})\longrightarrow \Pic(X'_{\bar{k}})
$$
is surjective and has a finitely generated kernel (generated by the finitely many irreducible components of codimension $1$ of $\overline{X}'_{\bar{k}} \setminus X'_{\bar{k}}$). 
But $\Pic(\overline{X}'_{\bar{k}})$ is an extension 
of a finitely generated abelian group 
by the divisible torsion group  $\Pic^0_{\overline{X}'/k}(\bar{k})$,
and therefore the induced map
$$
V_\ell\Pic(\overline{X}'_{\bar{k}})\longrightarrow V_\ell\Pic(X^\prime_{\bar{k}})
$$
is an isomorphism.
Since 
$V_\ell\Pic(\overline{X}'_{\bar{k}})\cong \H^1(\overline{X}_{\bar{k}},\bbQ_\ell(1))$,
$V_\ell \H^1(X_{\bar{k}},\bbG_m)$ is pure of weight $-1$.
It follows that the generalised $1$-eigenspace of 
$\H^1(X_{\bar{k}},\bbQ_\ell(1))$ 
can be identified with 
$\H^0(X_{\bar{k}},\bbG_m)^{G_k}\otimes_\bbZ \bbQ_\ell= \H^0(X,\bbG_m)\otimes_\bbZ \bbQ_\ell$. 
This implies that the dimension of the generalised $q^{{\dim(X)}-1}$-eigenspace of 
$\H^1(X_{\bar{k}},\bbQ_\ell)^\vee(-{\dim(X)})$ is equal to 
$\rk\,\bbG_m(X)=\rk\, \H^0(X,\cO_{X}^{*})$ as desired.
	
Considering that 
$z(X) \colonequals \rk\, \bbG_m(X) - \rk\,\Pic(X)-\ord_{s=\dim(X)-1}\zeta(X,s)$ 
as defined in~\eqref{eq:def z}, 
to show that the analytic Tate Conjecture~\hyperlink{BS3}{B+S-D+2} for $X$ is equivalent to the finiteness of 
$\Br(X_{k^{\sep}})^{G_k}(\ell)$,  
it suffices to show that 
the dimension of the generalised $q^{\dim(X)-1}$-eigenspace 
of $\Frob^{-1}$ on $\H^2(X_{\bar{k}},\bbQ_\ell)^\vee(-{\dim(X)})$
(which equals the dimension of the generalised $1$-eigenspace of 
$\H^2(X_{\bar{k}},\bbQ_\ell(1))$) 
is equal to $\rk\, \Pic(X)$ if and only if $V_\ell\Br(X_{\bar{k}})^{G_k}$ vanishes. 
By Lemma \ref{ssplit}, there is an exact sequence\\
	$$ 0\longrightarrow \Pic(X)\otimes_\bbZ \bbQ_\ell \longrightarrow \H^2(X_{\bar{k}},\bbQ_\ell(1))^{G_k}\longrightarrow V_\ell\Br(X_{\bar{k}})^{G_k}\longrightarrow 0.
	$$
On the one hand, if the dimension of the generalised $1$-eigenspace of $\H^2(X_{\bar{k}},\bbQ_\ell(1))$ is equal to $\rk\, \Pic(X)$, 
then
$$ 
\Pic(X)\otimes_\bbZ \bbQ_\ell\longrightarrow \H^2(X_{\bar{k}},\bbQ_\ell(1))^{G_k}
$$
is an isomorphism, 
and it follows that $V_\ell \Br(X_{\bar{k}})^{G_k}=0$. 
On the other hand if $V_\ell \Br(X_{\bar{k}})^{G_k}=0$, 
then
$$
\Pic(X)\otimes_\bbZ \bbQ_\ell \longrightarrow \H^2 (X_{\bar{k}},\bbQ_\ell(1))^{G_k}
$$
is an isomorphism. 
Considering that the $G_k$-action on $\Pic(X_{\bar{k}})\otimes_\bbZ \bbQ_\ell$ is semisimple (as the $G_k$-action factors through a finite quotient), \cref{ssplit} implies that $\H^2(X_{\bar{k}},\bbQ_\ell(1))^{G_k}$ is a direct summand of $\H^2(X_{\bar{k}},\bbQ_\ell(1))$ as a $G_k$-representation. 
It follows that
the generalised $1$-eigenspace of $\Frob^{-1}$ on 
$\H^2(X_{\bar{k}},\bbQ_\ell(1))$ equals the $1$-eigenspace,
which coincides with $\Pic(X)\otimes_\bbZ \bbQ_\ell$. 
This completes the proof.
\end{proof}

\begin{corollary}
Assuming that the Tate Conjecture $\T^1(X/k,\ell)$ holds 
for all smooth projective varieties over a finite field $k$ of characteristic $p>0$, 
the analytic Tate 
Conjecture~\hyperlink{BS3}{B+S-D+2} holds for all smooth varieties over $k$.
\end{corollary}
\begin{proof}
Let $X$ be a smooth variety over $k$. 
According to de Jong's alteration theorem, there is a dense open subset
$U\subset X$, and a finite flat morphism $f\colon U'\rightarrow  U$, 
such that $U'$ admits a smooth projective compactification $\overline{U}'$. 
By assumption, $\T^1(\overline{U}'/k,\ell)$ holds, 
and hence $V_\ell\Br(\overline{U}'_{k^{\sep}})^{G_k}=0$. 
With \cref{bir} it follows that 
$V_\ell \Br(U'_{k^{\sep}})^{G_k} =0$ as well. 
But the same is true for $V_\ell \Br(U_{k^{\sep}})^{G_k}$ because the map
$V_\ell\Br(U_{k^{\sep}})^{G_k}\rightarrow V_\ell\Br(U'_{k^{\sep}})^{G_k}$ 
induced by $f$ is injective. 
Now we apply again \cref{bir} to obtain 
$V_\ell\Br(X_{k^{\sep}})^{G_k}=V_\ell\Br(\sU_{k^{\sep}})^{G_k}=0$. 
Finally the statement follows by the above \cref{mlem}.
\end{proof}

\subsection{The obstruction of the Tate Conjecture in the $p$-adic case} \label{sec:Tate p adic}

We work again over the finite field $k=\bbF_q$, $q=p^r$. 
{Denote by $K_0$ the fraction field of the Witt vectors $W(k)$.}
In this subsection we are interested in the obstruction of the Tate Conjecture in the $p$-adic case.
We start with a rather well-known result.

\begin{proposition} \label[proposition]{Br is obstruction to Tp}
Let $X/k$ be a smooth projective geometrically connected variety. Then there is a short exact sequence
$$
	0 \to \NS(X)_{\mathbb{Q}_{p}}\longrightarrow  \H^2_\cris(X/W(k))^{F=p}_{\mathbb{Q}_{p}} \longrightarrow  V_p\Br(X)\longrightarrow  0.
$$
\end{proposition}

\begin{proof}
As rational crystalline cohomology is a Weil cohomology theory
and as homological equivalence lies between algebraic and numerical equivalence,
we can again invoke Matsusaka's Theorem \cite[Thm.\,4]{Matsusaka_1957} 
to deduce that 
algebraic equivalence for divisors equals numerical equivalence for divisors 
modulo torsion.
Hence the cycle class map 
$\Pic(X)_{\mathbb{Q}_{p}} \rightarrow  \H^2_\cris(X/W(k))^{F=p}_{\mathbb{Q}_{p}}$ factors through  the map 
$$
\NS(X)_{\mathbb{Q}_{p}}  \longrightarrow  \H^2_\cris(X/W(k))^{F=p}_{\mathbb{Q}_{p}},
$$
which is moreover injective.
Next according to \cite[II.\,Thm.\,5.5]{Illusie_1979} 
$\H^2_\cris(X/W(k))^{F=p}_{\mathbb{Q}_{p}} \cong \H^2_{\fppf}(X,\mathbb{Q}_p(1))$
and this is compatible with cycle class maps.
But as explained in \cite[II.\,(5.8.5)]{Illusie_1979}  we have an exact sequence
$$
0\longrightarrow  \NS(X)_{\mathbb{Q}_{p}} \longrightarrow  \H^2_{\fppf}(X,\mathbb{Q}_p(1)) \longrightarrow  V_p\Br(X)\longrightarrow  0,
$$
and by what we just said this shows the (existence and) exactness of the desired sequence.
\end{proof}

In the case of a variety over a function field, 
the situation is a bit more subtle. 
Let again $S$ be a smooth projective geometrically irreducible variety over the finite field $k=\bbF_q$, $q=p^r$, with function field $K$.

\begin{proposition}
Let $X$ be a smooth projective geometrically connected variety over $  K$.
	The following are equivalent:
	\begin{enumerate}
		\item The map $\NS(X)_{\bbQ_p}\rightarrow   \H^2_\rig(X/  K,K_0)^{F=p}$ is an isomorphism.
		\item The map $\NS(X)_{K_0}\rightarrow   \H^2_\rig(X/  K,K_0)^{F^r=q}$ is an isomorphism.
	\end{enumerate}
In other words, both statements express the $p$-adic Tate Conjecture~\hyperlink{conj1}{$\T^1(X/K,\ell)$}.
\end{proposition}
\begin{proof}
Since the map in (ii) is injective,
the same holds for (i).
Since the implication $\text{(ii)} \Rightarrow \text{(i)}$  is obvious, it suffices to show the inverse $\text{(i)} \Rightarrow \text{(ii)}$. 
Thus assume that the map in (i) is an isomorphism.
By tensoring it with $K_0$, we get 
$$\NS(X)_{K_0}\cong \H^2_\rig(X/  K,K_0)^{F=p}\otimes_{\bbQ_p} K_0.$$
Then, the claim follows from the subsequent lemma for $V=\H^2_\rig(X/  K,K_0)$.
\end{proof}

\begin{lemma}
 For an $F$-isocrystal $V$ with coefficients in $K_0=W(k)\left[\frac{1}{p}\right]$, 
 there is an isomorphism $V^{F=p} \otimes_{\bbQ_p} K_0 \cong V^{F^r=q}$.
\end{lemma}
\begin{proof}
Let $\sigma$ be the arithmetic Frobenius element in $\mathrm{Gal}(k/\mathbb{F}_p)$. Then, $\sigma\mapsto p^{-1}F$ defines an action of $\mathrm{Gal}(k/\mathbb{F}_p)$ on $V$. Since $V^{F=p}=V^{G}$, the claims follows from Galois descent.
\end{proof}

\begin{remark}\label[remark]{remark: when Tate Conjecture for function fields}
The assumption $\coker\big(\NS(X)_{K_0}\rightarrow   \H^2_\rig(X/  K,K_0)^{F^r=q}\big)=0$ holds for example in the case that $X$ is an abelian variety
because the Tate Conjecture is known in this case.
More precisely, Zarhin in \cite{Zarhin_1974a,Zarhin_1974b} extended Tate's result in \cite{Tate_1966}
concerning the Tate Conjecture~\hyperlink{conj1}{$\T^n(X/K,\ell)$} 
from abelian varieties over finite fields (of characteristic $p$) to abelian varieties over fields which are finitely generated over $\bbF_p$. 
Hence by Proposition \ref{thm: equivalence Tate over function fields} $\T^1(X/  K,p)$ is also known to be true 
which means that  $\coker\big(\NS(X)_{K_0}\rightarrow   \H^2_\rig(X/  K,K_0)^{F^r=q}\big)$ vanishes. 
\end{remark}

\begin{assumption}
In the remainder of this {sub}section, we assume that the {base field}
$k$ is $\bbF_p$.
\end{assumption}

We now want to address the $p$-adic equivalent of \cref{mlem}
which concerns the obstruction of the analytic Tate Conjecture~\hyperlink{BS3}{B+S-D+2}. 
By \cite[\S2]{Petrequin_2003}, rigid cohomology is a Bloch--Ogus cohomology.

\begin{lemma} \label[lemma]{proper surjective implies injective of Hrig}
	Let $f\colon X \to Y$ be a proper surjective morphism between smooth and integral varieties over $\bbF_p$. Assume that $\dim{X} = \dim{Y} \equalscolon m$. Then $f^*\colon \H_{\rig}^i(Y) \to \H_{\rig}^i(X)$ is injective for any $i \geqslant 0$.
\end{lemma}

\begin{proof}
Since $f$ is proper, it induces a morphism $f^*\colon \H_{\rig,c}^i(Y) \to \H_{\rig,c}^i(X)$ for any $i\geqslant 0$.
	
It is straight forward to show injectivity of $f^\ast$ for $i=2m$. 
Indeed, if both $X,Y$ are geometrically connected both groups are $1$-dimensional vector spaces in characteristic $0$ and $f^*$ is multiplication by the degree. In general, one uses the fact that  each connected component of $Y$ is dominated by a connected component of $X$.
	
From this, one can deduce injectivity in general as follows. 
Let $0\neq\alpha\in \H_{\rig,c}^i(Y)$. 
By Poincaré duality, there is $\beta\in \H_{\rig,c}^{2m-i}(Y)$ such that $0\neq\alpha \cup \beta \in \H_{\rig,c}^{2m}(Y)$. 
By what we showed above, $0\neq f^\ast(\alpha\cup\beta)= f^*\alpha \cup f^*\beta$, where the equality holds because of compatibility of $f^\ast$ with the cup product. 
In particular, $f^\ast\alpha\neq 0$ and hence $f^\ast$ is injective.
\end{proof}

\begin{remark}
	In fact, there exists $f_*\colon \H_{\rig}^i(X) \to \H_{\rig}^i(Y)$ with $f_* \circ f^* = \deg(f)$ (and one has the adjunction $f_*\alpha \cup \beta= \alpha \cup f^*\beta$).
\end{remark}

\begin{lemma} \label[lemma]{pullback purely inseparable morphism is iso}
	If $f\colon X \to Y$ is proper {over $\bbF_p$} {with purely inseparable extension of the function fields $K(X)/K(Y)$},
	then there exists {an open dense subvariety} $U \subseteq Y$
	such that $f^*\colon \H_{\rig}^i(U) \to \H_{\rig}^i(V)$ is an isomorphism where $V \colonequals f^{-1}(U)$.
\end{lemma}

\begin{proof}
	By the previous lemma, $f^*$ is injective. 
	There exists $U \subseteq Y$ open dense {contained in the smooth locus} such that {a power of the absolute Frobenius} $F_U^n\colon U \to U$ factors  as $U \rightarrow V \xrightarrow{f} U$ for some $n > 0$ and the first arrow is surjective. Hence again by the previous lemma $\dim \H_{\rig}^i(V) \leq \dim \H_{\rig}^i(U)$, hence $f^*$ is an isomorphism.
\end{proof}

\begin{corollary}\label[corollary]{splitnesscor}
	Let $X/\bbF_p$ be a smooth integral variety. Then $\Pic(X)_{\bbQ_p} \to \H_{\rig}^2(X)$ is injective and split with $F$-action. (Here $F$ acts on $\Pic(X)_{\bbQ_p}$ by multiplication by $p$.)
\end{corollary}

\begin{proof}
	Let $\overline{X}$ be a projective compactification of $X$. Let $f\colon Y \to \overline{X}$ be an alteration such that $Y$ is smooth and projective over $\bbF_p$. Let $V = f^{-1}(X)$, so $f\colon V \to X$ is proper. By the lemma above, $f^*\colon \H_{\rig}^2(X) \to \H_{\rig}^2(V)$ is injective. Similarly, $\Pic(X)_{\bbQ_p}\rightarrow \Pic(V)_{\bbQ_p}$ is also injective. (By \cite[Ex.\,8.1.7, p.\,135]{fulton2013intersection}, the pull-back $f^*$ and $f_*$ on the Chow group $\CH^1$ satisfy $f_*\circ f^*=\deg(f)$.) Hence it suffices to show that $\Pic(V)_{\bbQ_p} \to \H_{\rig}^2(V)$ is injective and split.
	
	Consider {the exact sequence} $\H_{\rig}^2(Y) \to \H_{\rig}^2(V) \to \H_{\rig,Y \setminus V}^3(Y)$, where $\H_{\rig,Y \setminus V}^3(Y)$ has weight different from the generalised $(F = p)$-eigenspace $\H_{\rig}^2(V)^{F=p,g} \colonequals \bigcup_{n \geqslant 1}\ker((F-p)^n\colon \H_{\rig}^2(V) \to \H_{\rig}^2(V))$. 
	Indeed, by the Gysin isomorphism \cite[Thm.\,2.10]{Petrequin_2003} $\H_{\rig,Y \setminus V}^3(Y)\cong \H_{\rig,c}^{3-2\codim Y\setminus V}(Y \setminus V)$, which has weight $\leq 3-2\codim Y\setminus V$. 
	On the other hand $\H_{\rig}^2(V) $ has weight $\geqslant 2$.
 Hence $\H_{\rig}^2(Y)^{F=p,g} \twoheadrightarrow \H_{\rig}^2(V)^{F=p,g}$ is surjective. Let $W \subseteq \H_{\rig}^2(Y)^{F = p,g}$ be a direct summand complement of $\Pic(Y)_{\bbQ_p}$.
 Thus $W \cap \Pic(Y)_{\bbQ_p} = 0$, and $\H_{\rig}^2(Y) \to \H_{\rig}^2(V)$ restricted to $W$ is injective.
 Indeed, the intersection of $W$ with $\ker(\H_{\rig}^2(Y) \to \H_{\rig}^2(V))$ is contained in $\ker(\H_{\rig}^2(Y) \to \H_{\rig}^2(V))\cap\H_{\rig}^2(Y)^{F = p,g}$, which is equal to $\ker(\H_{\rig}^2(Y) \to \H_{\rig}^2(V))\cap\H_{\rig}^2(Y)^{F = p} \subseteq \Pic(Y)_{\bbQ_p}$. 
 Its image gives a direct summand complement  of $\Pic(V)_{\bbQ_p}$.
\end{proof}

\begin{remark}
Since $X/\bbF_p$, $\H_{\rig}^2(X)$ is a $W(\bbF_p)[p^{-1}] = \bbQ_p$-vector space, {the Frobenius} $F$ is $\bbQ_p$-linear on $\H_{\rig}^2(X)$. 

By \cite[Thm.\,2.10]{Petrequin_2003}, the group $\ker(\H_{\rig}^2(Y) \to \H_{\rig}^2(V))$ is generated by cycle classes. 
{Hence for sufficiently large $m$, we have  so $F^m=p^m$ on  $\ker(\H_{\rig}^2(Y) \to \H_{\rig}^2(V))$.}
\end{remark}

\begin{lemma}\label[lemma]{lemma coker invariant under shrinking}
Let $X/\bbF_p$ be a smooth integral variety and $U\subseteq X$ be an open dense subvariety. Then $\coker(\Pic(X)_{\bbQ_p} \to \H_{\rig}^2(X)^{F=p})\cong \coker(\Pic(U)_{\bbQ_p} \to \H_{\rig}^2(U)^{F=p})$
\end{lemma}
\begin{proof}
The kernel $\ker(\H_{\rig}^2(X)^{F=p,g}\rightarrow \H_{\rig}^2(U)^{F=p,g})\subseteq \Pic(X)_{\bbQ_p}$ is {a} direct summand of $\H_{\rig}^2(X)^{F=p,g}$ since $\Pic(X)_{\bbQ_p}$ is a direct summand of $\H_{\rig}^2(X)^{F=p,g}$ by~\cref{splitnesscor}. Since $\H_{\rig}^2(X)^{F=p,g}\rightarrow \H_{\rig}^2(U)^{F=p,g}$ is surjective, taking $F=p$, we get a surjective map
$$\H_{\rig}^2(X)^{F=p}\twoheadrightarrow \H_{\rig}^2(U)^{F=p}.$$
Then the claim follows from the snake lemma, noting that $\Pic(X) \twoheadrightarrow \Pic(U)$ is surjective.
\end{proof}
\begin{remark}
By \cref{splitnesscor}, $\coker(\Pic(X)_{\bbQ_p} \to \H_{\rig}^2(X)^{F=p})=\coker(\Pic(X)_{\bbQ_p} \to \H_{\rig}^2(X))^{F=p}
.$
\end{remark}

\begin{corollary}\label[corollary]{Gdescent}
Let $f\colon X \rightarrow Y$ be a proper surjective morphism between smooth integral varieties over $k=\bbF_p$. Assume that there is a finite group $G$ acting on $X \rightarrow Y$ with a trivial action on $Y$ such that the field extension $K(X)^G/K(Y)$ is purely inseparable.

Then the pullback induces an isomorphism
$$\coker(\Pic(Y)_{\bbQ_p} \to \H_{\rig}^2(Y))^{F=p}\xrightarrow{\sim}(\coker(\Pic(X)_{\bbQ_p} \to \H_{\rig}^2(X))^{F=p})^G.
$$
\end{corollary}

\begin{proof}
Since $\coker(\Pic(Y)_{\bbQ_p} \to \H_{\rig}^2(Y))^{F=p}$ does not change when shrinking $Y$, we may assume that $f$ is a composition of a finite flat pure inseparable morphism with a finite Galois \'etale morphism. It suffices to prove the claim separately for these two kind morphisms.

First, assume that $f$ is a finite \'etale Galois cover.
By the spectral sequence 
$$E_2^{p,q}=\H^p(G, \H^q(X,\bbG_m))\Rightarrow \H^{p+q}(Y,\bbG_m),$$
$$\Pic(Y)\rightarrow \Pic(X)^G$$
has a kernel and a cokernel of finite exponent since $G$ is finite. From \'etale Galois descent for rigid cohomology, one also has that
$$\H_{\rig}^2(Y)\rightarrow\H_{\rig}^2(X)^G$$
is an isomorphism since rigid cohomology groups are $K_0$-vector spaces. Since $G$ is {a} finite group, $\H^0(G,-)$ is exact for finite dimensional $\bbQ_p$-representations. Hence, the pullback induces an isomorphism
$$\coker(\Pic(Y)_{\bbQ_p} \to \H_{\rig}^2(Y))\rightarrow\coker(\Pic(X)_{\bbQ_p} \to \H_{\rig}^2(X))^G
$$
Taking $F=p$, this proves the claim for $f$ being finite \'etale and Galois.

Second, assume that $f$ is finite flat and purely inseparable. Using~\cref{lemma coker invariant under shrinking}, {and possibly} shrinking $Y$, we may assume that the absolute Frobenius $F_X^n\colon X\rightarrow X$ can be written as a composition
$$Y\rightarrow X\stackrel{f}{\rightarrow} Y.$$
Thus the pullbacks $\Pic(Y)_{\bbQ_p}\rightarrow\Pic(X)_{\bbQ_p}$ and $\H_{\rig}^2(Y) \rightarrow \H_{\rig}^2(X)$ are isomorphisms: For the Picard group it follows from a pullback-pushforward argument and for the rigid cohomology by~\cref{proper surjective implies injective of Hrig}. This proves the claim for $f$ purely inseparable. 
\end{proof}

\begin{proposition}\label[proposition]{prop: nr Br Tate-obstruction for open}
Let $k = \bbF_p$ and $X/k$ be an irreducible (hence integral) smooth variety.
Then there is an exact sequence
$$
\Pic(X)_{\mathbb{Q}_{p}}\longrightarrow  \H^2_\rig(X/k,K_0)^{F=p} \longrightarrow  V_p\Br_{\nr}( K(X)/k)\longrightarrow  0.
$$
In other words, there exists a canonical isomorphism
$$
V_p\Br_{\nr}( K(X)/k)\cong \coker(\Pic(X)_{\bbQ_p}\rightarrow \H^2_\rig(X/k,K_0))^{F=p}.
$$
\end{proposition}

\begin{proof}
Let $\overline{X}$ be a Nagata compactification of $X$ over $k$. 
By de Jong's alteration theorem \cite[Thm.\,7.3]{deJong_1996}, there exists a finite extension $k'/k$ and an alteration
$$
f\colon \overline{X}_1\rightarrow \overline{X}_{k'},
$$
where $\overline{X}_1/k'$ now is smooth projective,
and a subgroup $G\subset \Aut_{\overline{X}_{k'}} (\overline{X}_1)$
such that $ K(\overline{X}_{k'})\subset  K(\overline{X}_1)^G$ is purely inseparable. 
Denote by $K_0$ the fraction field of $W(k)$. 
By \cref{Gdescent}, we get a natural map
\begin{align*}
\coker\big(\Pic(X_{k'})_{\bbQ_p}\rightarrow \H^2_\rig(X_{k'},K_0))^{F=p}\big)
\xrightarrow{\sim} &(\coker\big(\Pic(\overline{X}_1)_{\bbQ_p}\rightarrow \H^2_\rig(\overline{X}_1,K_0))^{F=p}\big)^G\\
\subset &\coker\big(\Pic(\overline{X}_1)_{\bbQ_p}\rightarrow \H^2_\rig(\overline{X}_1,K_0)\big)^{F=p}.
\end{align*}
According to \cref{lem: unram Br}, 
we also get a natural map 
$$
V_p\Br_{\nr}(K(X_{k'})/k) \xrightarrow{\sim} V_p\Br_{\nr}(K(\overline{X}_1)^G/k) \cong V_p\Br_{\nr}(K(\overline{X}_1))^G \subset V_p\Br_{\nr}(K(\overline{X}_1)/k).
$$
In other words, 
\begin{align*}
f^\ast \coker\big(\Pic(X_{k'})_{\bbQ_p}\rightarrow \H^2_\rig(X_{k'},K_0)\big)^{F=p} &\cong  \big(\coker\big(\Pic(\overline{X}_1)_{\bbQ_p}\rightarrow \H^2_\rig(\overline{X}_1,K_0)\big)^{F=p}\big)^G\\
f^\ast V_p\Br_{\nr}(K(X_{k'})/k) &\cong V_p\Br_{\nr}(K(\overline{X}_1)/k)^G
\end{align*}
Since $\overline{X}_1$ is smooth and projective over $k'$, 
its rigid and rational crystalline cohomology coincide and by \cref{Br is obstruction to Tp}
there exists a canonical isomorphism
$$
g_{X_1}\colon V_p\Br_{\nr}( K(\overline{X}_1)/k)\xrightarrow{\sim}
\coker\big(\Pic(\overline{X}_1)_{\bbQ_p}\rightarrow \H^2_\rig(\overline{X}_1,K_0)^{F=p}\big), 
$$
which is $G$-invariant. 
Thus this induces a morphism
$$
g_{X_{k'}}\colon V_p\Br_{\nr}( K(X_{k'})/k) \rightarrow \coker\big(\Pic(X_{k'})_{\bbQ_p}\rightarrow \H^2_\rig(X_{k'},K_0)^{F=p}\big)
$$
which is determined by the diagram
$$
\xymatrix{
V_p\Br_{\nr}( K(X_{k'})/k) \ar[r]^{f^\ast} \ar[d]_{g_{X_{k'}}} & V_p\Br_{\nr}( K(X_1)/k) \ar[d]_{g_{\overline{X}_1}} \\
\coker\big(\Pic(X_{k'})_{\bbQ_p}\rightarrow \H^2_\rig(X_{k'},K_0)^{F=p}\big) \ar[r]^{f^\ast} & \coker\big(\Pic(\overline{X}_1)_{\bbQ_p}\rightarrow \H^2_\rig(\overline{X}_1,K_0)^{F=p}\big).
}
$$
Thus it suffices to show that $g_{X_{k'}}$ is compatible with the action by $\Gal(k'/k)$. 

To show this, let $\sigma\colon \Spec(k')\rightarrow \Spec(k')$ be induced by an element $\sigma\in\Gal(k'/k)$
and consider the diagram
$$
\xymatrix{
\overline{X}_1^\sigma \ar[r]^{f^\sigma} \ar[d]^{\sigma} & \overline{X}_{k'} \ar[d]^{\sigma} \ar[r] & \Spec(k') \ar[d]^\sigma\\
\overline{X}_1 \ar[r]^f & \overline{X}_{k'} \ar[r] & \Spec(k'),
}
$$
where $\overline{X}_1^\sigma= \overline{X}_1\times^\sigma_{k'}\Spec(k')$. 
Set $\tau=\sigma^{-1}$. 

We observe first that
\begin{equation}\label{claim2}
g_{\overline{X}_1}\circ \tau^\ast = \tau^\ast\circ g_{\overline{X}_1^\sigma}.
\end{equation}
Indeed, $\tau\colon \overline{X}_1 \rightarrow \overline{X}_1^\sigma$ is a $k$-morphism between smooth projective varieties. 
In this case, rigid and rational crystalline cohomology coincide, and as explained in \cref{Br is obstruction to Tp} 
the subspace where $F=p$ compares to fppf-cohomology.
What is more, this comparison is compatible with the action by $\Gal(k'/k)$, so that we obtain a commutative diagram
$$
\xymatrix{
\H^2_{\rig}(\overline{X}^\sigma_1,K_0)^{F=p} \ar[r]^{\tau^\ast} & \H^2_{\rig}(\overline{X}_1,K_0)^{F=p}\\
\H^2_{\fppf}(\overline{X}_1^\sigma,\bbQ_p(1)) \ar[u] \ar[r]^{\tau^\ast} & \H^2_{\fppf}(\overline{X}_1,\bbQ_p(1)) \ar[u]
}
$$
where the vertical maps induce $g_{\overline{X}_1}$ and $g_{\overline{X}_1^\sigma}$, respectively. 
This shows the claim.

Moreover, there is a commutative diagram
$$
\xymatrix{
\coker\big(\Pic(X_{k'})_{\bbQ_p}\rightarrow \H^2_\rig(X_{k'},K_0)^{F=p}\big) \ar[r]^{\tau^\ast} \ar[d]^{(f^\sigma)^\ast} & \coker\big(\Pic(X_{k'})_{\bbQ_p}\rightarrow \H^2_\rig(X_{k'},K_0)^{F=p}\big) \ar[d]^{f^\ast}\\
\coker\big(\Pic(\overline{X}_1^\sigma)_{\bbQ_p}\rightarrow \H^2_\rig(\overline{X}_1^\sigma,K_0)^{F=p}\big) 
\ar[r]^{\tau^\ast} & \coker\big(\Pic(\overline{X}_1)_{\bbQ_p}\rightarrow \H^2_\rig(\overline{X}_1,K_0)^{F=p}\big)
}
$$

Let $a\in V_p\Br_{\nr}(K(X_{k'})/k)$ and $b\colonequals g_{\overline{X}_1}(a)\in \coker\big(\Pic(X_{k'})_{\bbQ_p}\rightarrow \H^2_\rig(X_{k'},K_0)^{F=p}\big)$. 
Then 
\begin{equation}\label{claim1}
g_{\overline{X}_1^\sigma}\big((f^\sigma)^\ast(a)\big) = (f^\sigma)^\ast(b).
\end{equation}
Indeed, $f^\sigma\colon \overline{X}_1^\sigma\rightarrow \overline{X}_k'$ 
satisfies the same conditions as $f$.
Let now $h\colon Y\rightarrow \overline{X}_1\times_{\overline{X}_k'} \overline{X}_1^\sigma$ be an alteration, such that $Y$ is regular. 
By the universal property of the fibre product there is a commutative diagram
$$
\xymatrix{
Y \ar[r]^{p_1} \ar[d]^{p_2} & \overline{X}_1^\sigma \ar[d]^{f^\sigma}\\
\overline{X}_1 \ar[r]^f & \overline{X}_{k'}.
}
$$
For any $a\in V_p\Br_{\nr}(K(\overline{X}_k')/k)$ and 
$b=g_{X_{k'}}(a)$ we have
$f^\ast(b)=g_{X_1}f^\ast(a)$, and hence
$$
p_1^\ast\circ f^\ast(b)= p_1^\ast\circ g_{X_1}(f^\ast(a))=
g_Y(p_1^\ast f^\ast(a)).
$$
But by (\ref{claim2}) we have $p_1^\ast\circ g_{X_1}= g_Y\circ p_1^\ast$, so that we compute
$$
(f\circ p_1)^\ast(b) = g_Y((f\circ p_1)^\ast y) = p_2^\ast\circ g_{X_1^\sigma}((f^\sigma)^\ast(a)).
$$
From $(f\circ p_1)^\ast(b)=(f^\sigma\circ p_2)^\ast(b)$ we obtain
$$
p_2^\ast((f^\sigma)^\ast(b))=p_2^\ast\circ g_{X_1}\circ((f^\sigma)^\ast(a)).
$$
But because of the injectivity of $p_2^\ast$, 
we obtain $g_{X_1^\sigma}(f^\sigma)^\ast(a))= (f^\sigma)^\ast(b)$ as desired.

Since  $f^\ast(\tau^\ast(a)) = \tau^\ast(f^\sigma)^\ast(a)$,
it follows from (\ref{claim1}) that 
$$
g_{\overline{X}_1}\circ f^\ast(\tau^\ast(a)) = 
(g_{\overline{X}_1}\circ \tau^\ast)\big((f^\sigma)^\ast(a)\big).
$$

Using (\ref{claim2}) and (\ref{claim1}), we obtain the following equalities:
\begin{align*}
(g_{\overline{X}_1}\circ\tau^\ast)\big((f^\sigma)^\ast(a)\big) &=
\tau^\ast\circ g_{\overline{X}^\sigma_1}\circ(f^\sigma)^\ast(a)\\
&=\tau^\ast\circ(f^\sigma)^\ast(b)\\
&=f^\ast\circ\tau^\ast(b).
\end{align*}
Thus
\begin{align*}
f^\ast\circ\tau^\ast(b) &= g_{\overline{X}_1}\circ f^\ast(\tau^\ast(a))\\
&=f^\ast\circ g_{X_{k'}}(\tau^\ast(a)),
\end{align*}
and therefore $g_{X_{k'}}(\tau^\ast(a))= \tau^\ast(b)$, 
which shows that $g_{X_{k'}}$ is compatible with Galois descent.
This finishes the proof of the proposition.
\end{proof}

We conclude this section with the following main theorem.

\begin{theorem}\label[theorem]{thm: obstruction analytic Tate}
Let $X/k$ be an irreducible (hence in particular integral) smooth variety. 
Then the analytic Tate Conjecture~\hyperlink{BS3}{B+S-D+2} holds if and only if 
the $p$-torsion subgroup  $\Br_{\nr}(K(X)/k)(p)$ of its unramified Brauer group is of finite exponent.
\end{theorem}

\begin{proof}
By~\cref{prop: nr Br Tate-obstruction for open}, it suffices to show that 
{ the analytic Tate Conjecture~\hyperlink{BS3}{B+S-D+2} is equivalent to the equation $\coker(\Pic(X)_{\bbQ_p}\rightarrow \H^2_\rig(X)^{F=p})=0$.}
 Replacing $\ell$-adic cohomology by rigid cohomology in the proof of~\cref{mlem} using the theory of weights for rigid cohomology~\cite[Thm.\,5.3.2]{Kedlaya_2006} and the Etesse--Le Stum trace formula~\cite[{Thm.\,6.3}]{EtesseLestum_1993}, we are reduced to show that
$$\Pic(X)_{\bbQ_p}=\H^2_\rig(X)^{F=p}\iff
\Pic(X)_{\bbQ_p}=\H^2_\rig(X)^{F=p,g}.$$
By \cref{splitnesscor}, there exists a direct summand $W$ of $\H^2_\rig(X)^{F=p,g}$ such that
$$\Pic(X)_{\bbQ_p}\oplus W=\H^2_\rig(X)^{F=p,g}.$$
So $\Pic(X)_{\bbQ_p}\oplus W^{F=p}=\H^2_\rig(X)^{F=p}$. If $W^{F=p}=0$, then $W^{F=p,g}=0$. Since $W=W^{F=p,g}$, $W=0$.
\end{proof}

\section{Tate--Shafarevich and Brauer groups for abelian varieties over function fields} \label{sec:Sha and Brauer}

The purpose of this section is twofold: 
firstly to give a definition of the Tate--Shafarevich group
suited for the case of an abelian variety over a higher-dimensional function field;
and secondly to study its relation with the Brauer group in both the $\ell$-adic and the $p$-adic case.
Our main result is the proof of Theorem~\hyperlink{bigtheorem}{D}.
The technical problems are the lack of resolution of singularities in positive characteristic and that N\'eron models need not exist over higher-dimensional bases if the abelian variety does not have good reduction everywhere. (For a criterion when they exist, see~\cite[Thm.~1.2]{Holmes2019}.)

\subsection{Definition of the Tate--Shafarevich group and some basic properties}

Recall that to an abelian variety $A$ over a global field $K$, one may associate the Tate--Shafarevich group $\Sha(A/K)$. 
It classifies everywhere locally trivial $A/K$ torsors
and the definition can be given in terms of Galois cohomology as
$$\Sha(A/K)\colonequals\ker\Big(\H^1(G_{K},A) \rightarrow  \prod_\nu \H^1(G_{K_\nu},A)\Big)$$
where $\nu$ runs over all places of $K$ and $K_\nu$ is the completion of $K$ with respect to $\nu$.

For a finitely generated field $K$ over a field $k$, 
and $A/K$ which extends to an abelian scheme $\sA$ over a smooth proper base $S$ with function field $K$, 
a definition is given in \cite{Keller_2016}, where 
$$\Sha(A/K) \colonequals \H^1(S,\sA).$$
We call this the ``good reduction case''.
Without this assumption, we still have the following general definition, which we will use in this section:

\begin{definition}[Tate--Shafarevich group]\label[definition]{def: Shafarevich-Tate for abelian varieties}
Let $K$ be a finitely generated field over $k$ and $A/K$ an abelian variety.
\begin{enumerate}
\item
Define the \textit{Tate--Shafarevich group} of $A/ K$ by
$$
\Sha'(A/K) \colonequals \ker\Big(\H^1(K,A) \rightarrow \prod_{s\in div\Val(K)}\H^1(K_s^{\sh},A)\Big),
$$
where $s$ runs through runs through the divisorial valuations $div\Val(  K)$  of $K$, 
and for such an $s$, $K_s^{\sh}$ is the fraction field of a strict Henselisation of the local ring $\sO_s$ associated to $s$. 

\item 
Assume given an integral regular Noetherian scheme $S$ with function field ${K}(S)=K$. 
Define the \emph{Tate--Shafarevich group} of $A$ with respect to $S$ by
\[
\Sha_S(A) \colonequals \ker\Big(\H^1({K}(S),A) \to \prod_{s \in S^1}\H^1({K}(S)_s^{\sh},A)\Big),
\]
where $S^1$ denotes the set of codimension-$1$ points of $S$ and $  {K}(S)_s^{\sh}$ the fraction field of a strict Henselisation of the local ring $\sO_{S,s}$.	
\end{enumerate}
\end{definition}

We will see in \cref{twovaluations} that the first of the above definitions coincides with the one given in the introduction. 
Note that the second definition a priori depends on $S$. 
Ideally one would like to work with a definition that only depends on the function field of $S$.  
The question to which extent this is possible is one of the topics of
this section.

One of the properties of the Tate--Shafarevich group which will be important later on is its invariance under Weil restriction.

\begin{proposition}\label[proposition]{prop: Sha invariant under Weil}
Let $L/K$ be a finite separable extension of function fields and $A$ an abelian variety over $L$. 
Denote by $\Res_{L/K}(A)$ the Weil restriction of $A$. 
Then there is an isomorphism
$$
\Sha'(A/L) \cong \Sha'(\Res_{L/K}(A)/K).
$$
\end{proposition}
\begin{proof}
We note first, that the statement makes sense, since the Weil restriction of an abelian sheaf is itself an abelian sheaf. 
Let $\alpha\colon \Spec(L) \to \Spec(K)$. 
We have the following isomorphisms
\begin{align*}
\H^1(K,\alpha_\ast A) &\cong \H^1(L, A),\\
\H^1(K^{\sh}_s,\alpha_{s,\ast} A) &\cong \H^1(L^{\sh}_s, A), \text{ for any $s\in div\Val(K)$}.
\end{align*}
Indeed, more generally for a finite morphism $\alpha: Y \rightarrow Z$ 
and an étale sheaf $\sF$ on $Y$, the Leray spectral sequence
$$
\H^i(Z, \R^j\alpha_*\sF) \Longrightarrow \H^{i+j}(Y, \sF)
$$
degenerates on the $E_2$-page because push-forwards by finite morphisms 
are acyclic for the \'etale cohomology~\cite[Cor.\,II\,3.6]{Milne_1980}. 
Consequently, we obtain isomorphisms 
$\H^i(Z, \alpha_*\sF) \cong \H^i(Y, \sF)$ for all $i \geqslant 0$.

Using the surjective map $div\Val(L) \rightarrow div\Val(K)$ with finite fibres 
and the definition of the Tate--Shafarevich group, we obtain a commutative diagram
$$
\xymatrix{
0 \ar[r] & \Sha'(A/L)\ar[d] \ar[r] & \H^1(L, A) \ar@{=}[d]^\sim \ar[r] & 
\prod_{s\in div\Val(L)} \H^1(L^{\sh}_s, A) \ar@{->>}[d] \\
0\ar[r] & \Sha'(\alpha_\ast A/K) \ar[r] & \H^1(K,\alpha_\ast A) \ar[r] & 
\prod_{s\in div\Val(K)} \H^1(K^{\sh}_s,\alpha_{s,\ast} A) 
}
$$
and a simple diagram chase implies that $\Sha'(A/L)\cong \Sha'(\alpha_\ast A/K)$. 
Now the claim follows, 
since  the Weil restriction of an abelian variety is exactly the sheaf push-forward for the associated abelian sheaf, that is
$\Res_{L/K}(A)= \alpha_\ast A$ and 
$\Res_{\Phi_s/K_s}(A)= \alpha_{s,\ast} A$ for any $s\in div\Val(K)$. 
\end{proof}

\subsection{Fibrations of Brauer groups and $\ell$-torsion}

We first provide some technical preliminaries.

\begin{proposition} \label[proposition]{picard-brauer}
	Let $\sY$ be an irreducible regular scheme with function field $K=K(\sY)$.
	Let $\rho\colon \sX \rightarrow \sY$ be a smooth proper morphism such that its  generic fibre $X/K$ is geometrically connected. 
	Denote by $j\colon\Spec( K)\rightarrow  \sY$ the inclusion of the generic point of $\sY$. 
	
	\begin{enumerate}
		\item The natural morphism of étale sheaves on $\sY$
		$$
		\R^1\rho_*\bbG_m\longrightarrow j_*j^*\R^1\rho_*\bbG_m
		$$
		is an isomorphism.
		
		\item For any prime $\ell$ invertible on $\sY$
		the natural map
		$$
		\R^2\rho_*\bbG_m(\ell)\longrightarrow j_*j^*\R^2\rho_*\bbG_m(\ell)
		$$
		is an isomorphism.
	\end{enumerate}
\end{proposition}
\begin{proof}
	This is just~\cite[Lemma~3.1]{qinbrauer0}, but for a higher-dimensional base.
\end{proof}

\begin{theorem}\label[theorem]{big}
	Let $\rho\colon\sX\rightarrow \sY$ be a dominant morphism between smooth geometrically connected varieties over a finitely generated field $k$ of characteristic $p\geqslant 0$. Let $K=K(\sY)$ be the function field of $\sY$. Assume that the generic fibre $X$ of $\rho$ is smooth projective geometrically connected over $  K$. Set $  K^\prime \colonequals  Kk^{\sep}$.  Let $\ell\neq p$ be a prime. 
The morphism
$$
V_\ell\Br(\sX_{k^{\sep}})^{G_k}\rightarrow V_\ell \Br(X_{  K^\sep})^{G_{  K}}
$$
is surjective, and its kernel fits into a canonical short exact sequence of the form
	$$
	0\longrightarrow V_\ell\Br(\sY_{k^{\sep}})^{G_k}\longrightarrow \ker \big(V_\ell\Br(\sX_{k^{\sep}})^{G_k}\rightarrow V_\ell \Br(X_{  K^\sep})^{G_{  K}}\big) \longrightarrow V_\ell\Sha_{K'}(\Pic^0_{X/  K})^{G_k}\longrightarrow 0,
	$$
where $V_\ell\Sha_{K'}(\Pic^0_{X/  K})^{G_k}:=V_\ell\Sha_{\sY_{k^{\sep}}}(\Pic^0_{X/  K})^{G_k}$  is by \cref{prop: sha} independent of models. 
\end{theorem}
\begin{proof}
	By shrinking $\sY$ and using~\cref{bir}, we may assume that $\rho$ is projective and smooth. 
	
	Firstly, we will show that the natural map
	$$V_\ell\Br(\sX_{k^{\sep}})\longrightarrow V_\ell\Br(X_{  K^\sep})^{G_{  K^\prime}}$$
	is surjective. Consider the Leray spectral sequence
	$$E^{p,q}_2=\H^p(\sY_{k^{\sep}},\R^q\rho_*\bbQ_\ell(1))\Rightarrow \H^{p+q}(\sX_{k^{\sep}},\bbQ_\ell(1)). $$
	By Deligne's Lefschetz criteria (cf.\,\cite{Deligne_1968}), the above spectral sequence degenerates at $E_2$. Thus, we get a canonical surjective map
	$$\H^2(\sX_{k^{\sep}},\bbQ_\ell(1))\longrightarrow \H^0(\sY_{k^{\sep}},\R^2\rho_*\bbQ_\ell(1)).
	$$
	Since $\R^2\rho_*\bbQ_\ell(1)$ is lisse, we have
	$$\H^0(\sY_{k^{\sep}},\R^2\rho_*\bbQ_\ell(1))=\H^2(X_{  K^\sep},\bbQ_\ell(1))^{G_{  K^\prime}}.$$
	It follows that the canonical map
	$$\H^2(\sX_{k^{\sep}},\bbQ_\ell(1))\longrightarrow \H^2(X_{  K^\sep},\bbQ_\ell(1))^{G_{  K^\prime}}$$
	is surjective. Consider the following commutative diagram with exact rows:
	\begin{displaymath}
		\xymatrix{
			\H^2(\sX_{k^{\sep}},\bbQ_\ell(1))\ar[r]\ar[d] & V_\ell\Br(\sX_{k^{\sep}})\ar[d]\ar[r] &0
			\\
			\H^2(X_{  K^\sep},\bbQ_\ell(1))^{G_{  K\prime}}\ar[r] & V_\ell\Br(X_{  K^\sep})^{G_{  K^\prime}}	\ar[r]& 0
		}
	\end{displaymath}
	Since the first column is surjective, the map 
	$$V_\ell\Br(\sX_{k^{\sep}})\longrightarrow V_\ell\Br(X_{  K^\sep})^{G_{  K^\prime}}$$
	from the sequence in the statement of the theorem is also surjective.
	
	Secondly, we will show the surjectivity of the displayed morphism of the theorem. By the Leray spectral sequence
	$$\H^p(\sY_{k^{\sep}},\R^q\rho_*\bbG_m)\Rightarrow \H^{p+q}(\sX_{k^{\sep}},\bbG_m),
	$$
	we get a long exact sequence
	$$ \Br(\sY_{k^{\sep}})\longrightarrow \ker\big(\Br(\sX_{k^{\sep}})\longrightarrow \H^0(\sY_{k^{\sep}},\R^2\rho_*\bbG_m)\big) \longrightarrow \H^1(\sY_{k^{\sep}},\R^1\rho_*\bbG_m)\longrightarrow \H^3(\sY_{k^{\sep}},\bbG_m).
	$$
	Note that without loss of generality, we can always replace $k$ by a finite Galois extension. There exists a finite Galois extension $L/  K$ such that $X(L)$ is not empty. Choose a $  K$-morphism $\Spec(L)\rightarrow X$. We may assume that $L$ is the function field of {a} smooth variety $\cZ$ over $k$. By shrinking $\sY$ and $\cZ$, we can assume that $\Spec( L) \rightarrow X$ extends to a morphism $\cZ \rightarrow \sX$ such that {the composition 
	$\cZ \rightarrow \sX \rightarrow \sY$, which we denote by $\rho'$,} is finite flat. 
	Then we have a commutative diagram induced by the Leray spectral sequences for $\rho$ and $\rho^\prime$
	\begin{displaymath}
		\xymatrix{
			\H^1(\sY_{k^{\sep}},\R^1\rho_*\bbG_m)\ar[r]\ar[d]& \H^3(\sY_{k^{\sep}},\bbG_m)\ar[d]
			\\
			\H^1(\sY_{k^{\sep}},\R^1\rho^\prime_*\bbG_m)\ar[r]& \H^3(\sY_{k^{\sep}},\rho^\prime_*\bbG_m)}
	\end{displaymath}
	Since $\rho^\prime$ is finite, we have $\R^i\rho^\prime_*\bbG_m=0$ for $i>0$ and  $\H^3(\sY{k^{\sep}},\rho^\prime_*\bbG_m)=\H^3(\cZ_{k^{\sep}},\bbG_m)$. 
	{Moreover, again because of the finiteness of $\rho'$,}
	there is a canonical norm map $\rho_*^\prime\bbG_m\rightarrow \bbG_m$ which induces a norm map
	$$N\colon \H^i(\cZ_{k^{\sep}},\bbG_m) \longrightarrow \H^i(\sY_{k^{\sep}},\bbG_m)
	$$
	such that the composition of $N$ with the pull back map $\H^i(\sY_{k^{\sep}},\bbG_m) \rightarrow \H^i(\cZ_{k^{\sep}},\bbG_m)$ is equal to the multiplication by $\deg(\rho^\prime)$. It follows  that the second column in the above diagram has kernel killed by $\deg(\rho^\prime)$. Therefore, the first row has image killed by $\deg(\rho^\prime)$ and the sequence $0\rightarrow V_\ell\Br(\sY_{k^{\sep}})\rightarrow V_\ell\Br(\sX_{k^{\sep}})$ is split as $G_k$-representations. Thus, we get a split exact sequence of $G_k$-representations
	$$0\longrightarrow  V_\ell\Br(\sY_{k^{\sep}})\longrightarrow  \ker(V_\ell\Br(\sX_{k^{\sep}})\rightarrow  V_\ell \H^0(\sY_{k^{\sep}},\R^2\rho_*\bbG_m)) \longrightarrow  V_\ell \H^1(\sY_{k^{\sep}},\R^1\rho_*\bbG_m)\longrightarrow  0.$$
	In fact, the norm map $N$ induces a projection  
	$$\ker(V_\ell\Br(\sX_{k^{\sep}})\rightarrow V_\ell \H^0(\sY_{k^{\sep}},\R^2\rho_*\bbG_m))\longrightarrow V_\ell\Br(\sY_{k^{\sep}}).$$
	So we get an isomorphism
	$$
	\ker(V_\ell\Br(\sX_{k^{\sep}})\rightarrow V_\ell \H^0(\sY_{k^{\sep}},\R^2\rho_*\bbG_m))\cong V_\ell\Br(\sY_{k^{\sep}})\oplus V_\ell \H^1(\sY{k^{\sep}},\R^1\rho_*\bbG_m),
	$$
	where the direct sum decomposition depends on the choice of $\rho^\prime$. Next, we will show that there is a canonical isomorphism
	$$V_\ell \H^1(\sY_{k^{\sep}},\R^1\rho_*\bbG_m)\cong V_\ell\Sha_{\sY_{k^{\sep}}}(\Pic^0_{X/  K}).$$
	Let $j\colon\Spec(  K)\rightarrow \sY$ be the generic point. By \cref{picard-brauer}\,(i), the natural map
	$$\R^1\rho_*\bbG_m\longrightarrow j_*j^*\R^1\rho_*\bbG_m$$
	is an isomorphism. It follows that
	$$  \H^1(\sY_{k^{\sep}},\R^1\rho_*\bbG_m)\cong  \H^1(\sY_{k^{\sep}},j_*j^*\R^1\rho_*\bbG_m).$$
	There is an exact sequence
	$$ 0\longrightarrow \H^1(\sY_{k^{\sep}},j_*j^*\R^1\rho_*\bbG_m)\longrightarrow \H^1(  K^\prime,\Pic_{X/  K})\longrightarrow \prod _{s\in \sY_{k^{\sep}}}\H^1(  K^{\sh}_s,\Pic_{X/  K}).
	$$
	Denote by $\Sha_{\sY_{k^{\sep}}}(\Pic_{X/  K})$ the kernel of the third arrow. It suffices to show that
	$$V_\ell\Sha_{\sY_{k^{\sep}}}(\Pic^0_{X/  K})\cong V_\ell\Sha_{\sY_{k^{\sep}}}(\Pic_{X/  K}).$$
	By the canonical exact sequence
	$$ 0\longrightarrow \Pic^0_{X/  K}\longrightarrow \Pic_{X/  K}\longrightarrow \NS(X_{  K^\sep}) \longrightarrow 0,$$
	we get a long exact sequence
	$$\H^0(  K^\prime,\NS(X_{  K^\sep}))\longrightarrow \H^1(  K^\prime,\Pic^0_{X/  K})\longrightarrow \H^1(  K^\prime,\Pic_{X/  K}) \longrightarrow \H^1(  K^\prime,\NS(X_{  K^\sep})).
	$$
	Let $L/  K$ be a finite Galois extension such that $\Pic(X_L)\rightarrow \NS(X_{  K^\sep})$ is surjective. One can show that $\H^1(  K^\prime,\NS(X_{  K^\sep}))$ is killed by $[L:  K]\cdot |\NS(X_{  K^\sep})_{\tor}|$ and the cokernel of $\H^0(  K^\prime,\Pic_{X/  K})\rightarrow \H^0(  K^\prime,\NS(X_{  K^\sep}))$ is killed by $[L:  K]$. Thus, the kernel and the cokernel of 
	$$  \H^1(  K^\prime,\Pic^0_{X/  K}) \longrightarrow \H^1(  K^\prime,\Pic_{X/  K})
	$$
	are killed by $[L:  K]\cdot |\NS(X_{  K^\sep})_{\tor}|$. By the same argument, the kernel and cokernel of
	$$\prod _{s\in \sY_{k^{\sep}}}\H^1(  K^{\sh}_s,\Pic^0_{X/  K})\longrightarrow \prod _{s\in \sY_{k^{\sep}}}\H^1(  K^{\sh}_s,\Pic_{X/  K})$$
	are also killed by $[L:  K] \cdot |\NS(X_{  K^\sep})_{\tor}|$. Taking $V_\ell$, they become isomorphisms. Therefore, we have a canonical isomorphism
	$$
	V_\ell\Sha_{\sY_{k^{\sep}}}(\Pic^0_{X/  K})\cong V_\ell\Sha_{\sY_{k^{\sep}}}(\Pic_{X/  K}),$$
	{and hence a canonical isomorphism }
	$$
	V_\ell \H^1(\sY_{k^{\sep}}, \R^1\rho_*\bbG_m)\cong V_\ell\Sha_{\sY_{k^{\sep}}}(\Pic^0_{X/  K}).
	$$
	By~\cref{picard-brauer}\,(ii), the natural map
	$$\R^2\rho_*\bbG_m(\ell)\longrightarrow j_*j^*\R^2\rho_*\bbG_m(\ell)$$
	is an isomorphism. It follows directly that
	$$\H^0(\sY_{k^{\sep}},\R^2\rho_*\bbG_m)(\ell)\cong \Br(X_{  K^\sep})^{G_{  K^\prime}}(\ell).$$
	Combining the above exact sequences and isomorphisms, we get an exact sequence
	$$0\longrightarrow V_\ell\Br(\sY_{k^{\sep}})\oplus V_\ell \Sha_{\sY_{k^{\sep}}}(\Pic^0_{X/  K}/)
	\longrightarrow V_\ell\Br(\sX_{k^{\sep}}) \longrightarrow V_\ell\Br(X_{  K^\sep})^{G_{  K^\prime}}\longrightarrow 0,$$
	where the second arrow depends on the choice of $\rho^\prime$. Taking $G_k$-invariants, we have proved the surjectivity of the morphism as desired.
	
	It remains to show that the natural map
	$$V_\ell\Br(\sX_{k^{\sep}})^{G_k} \longrightarrow V_\ell\Br(X_{  K^\sep})^{G_{  K}}$$
	in the second exact sequence of the statement
	is surjective. To prove this, we will use a pull-back trick (cf.\,\cite{ColliotTheleneSkorobogatov_2013} or \cite{Yua2}). Let $W$ be a smooth projective geometrically connected curve over $  K$ contained in $X$ which is a complete intersection of hyperplane sections. By shrinking $\sY$, we can assume that $W$ admits a smooth projective model $\cW\rightarrow \sY$ where $\cW$ is the Zariski closure of $W$ in $\sX$. By further shrinking $\sY$, we can assume that the map $\cZ\rightarrow \sX$ chosen before factors through $\cZ\rightarrow \cW$. Thus, we get a commutative diagram with exact rows
	\begin{displaymath}
		\xymatrix{
			0\ar[r] & V_\ell\Br(\sY_{k^{\sep}})\oplus V_\ell \Sha_{\sY_{k^{\sep}}}(\Pic^0_{X/  K})
			\ar[r]\ar[d] & V_\ell\Br(\sX_{k^{\sep}}) \ar[r] \ar[d]& V_\ell\Br(X_{  K^\sep})^{G_{  K^\prime}}\ar[d] \ar[r] & 0\\
			0\ar[r] & V_\ell\Br(\sY_{k^{\sep}})\oplus V_\ell \Sha_{\sY_{k^{\sep}}}(\Pic^0_{W/  K})
			\ar[r] & V_\ell\Br(\cW_{k^{\sep}}) \ar[r] & V_\ell\Br(W_{  K^\sep})^{G_{  K^\prime}} \ar[r] &0
		}
	\end{displaymath}
	Taking $G_k$ invariants, we get 
	\begin{displaymath}
		\xymatrix{
			V_\ell\Br(\sX_{k^{\sep}})^{G_k}\ar[r]^a\ar[d] & V_\ell\Br(X_{  K^\sep})^{G_{  K}}\ar[d] \ar[r] & \H^1(G_k,V_\ell\Br(\sY_{k^{\sep}}))\oplus V_\ell \Sha_{\sY_{k^{\sep}}}(\Pic^0_{X/  K}) \ar[d]^c\\
			V_\ell\Br(\cW_{k^{\sep}})^{G_k} \ar[r] & V_\ell\Br(W_{  K^\sep})^{G_{  K}}\ar[r] & \H^1(G_k,V_\ell\Br(\sY_{k^{\sep}}))\oplus V_\ell \Sha_{\sY_{k^{\sep}}}(\Pic^0_{W/  K} )
		}
	\end{displaymath}
	Since $V_\ell\Br(W_{  K^\sep})^{G_{  K}}=0$ because $W_{  K^\sep}$ is a curve over a separably closed field, the second arrow in the bottom vanishes. Thus, to show that $a$ is surjective, it suffices to show that $c$ is injective. This actually follows from the fact
	$$ V_\ell \Sha_{\sY_{k^{\sep}}}(\Pic^0_{X/  K})\longrightarrow V_\ell \Sha_{\sY_{k^{\sep}}}(\Pic^0_{W/  K})$$
	is split as $G_k$-representations. By our choice of $W$ and the Lefschetz hyperplane section theorem, the pullback map 
	$$\H^1(X_{  K^\sep},\bbQ_\ell)\longrightarrow \H^1(W_{  K^\sep},\bbQ_\ell)$$
	is injective. This implies that the induced map $\Pic^0_{X/  K}\rightarrow \Pic^0_{W/  K}$ has a finite kernel. Therefore, there exists an abelian variety $A/  K$  and an isogeny $\Pic^0_{X/  K}\times A\rightarrow \Pic^0_{W/  K}$. It follows that
	$$V_\ell \Sha_{\sY_{k^{\sep}}}(\Pic^0_{X/  K})\oplus V_\ell \Sha_{\sY_{k^{\sep}}}(A)\cong V_\ell \Sha_{\sY_{k^{\sep}}}(\Pic^0_{W/  K}).
	$$
	This proves the splitness. Thus, the natural map
	$$
	V_\ell\Br(\sX_{k^{\sep}})^{G_k} \longrightarrow  V_\ell\Br(X_{  K^\sep})^{G_{  K}}
	$$
	is surjective.
\end{proof}

In the theorem above, we assume that the generic fibre is projective over $K$. In the following corollary, we will deduce the surjectivity of $V_\ell\Br(\sX_{k^{\sep}})^{G_k}\rightarrow V_\ell\Br(X_{K^\sep})^{G_K}$ from the above theorem for any fibration with a smooth generic fibre.

\begin{corollary}\label{hardsur}
	Let $k$ be a finite field or a number field. Let $\rho\colon\sX\rightarrow \sY$ be a dominant morphism between smooth geometrically connected varieties over $k$. Let $  K=K(\sY)$ be the function field of $\sY$ and $X/K$ be the generic fibre of $\rho$. 
	Let $\ell \neq \Char(k)$ be a prime. 
Assuming that $X$ is smooth over $ K$, then the natural map 
	$$V_\ell\Br(\sX_{k^{\sep}})^{G_k}\longrightarrow V_\ell\Br(X_{  K^\sep})^{G_{  K}}$$
	is surjective.
\end{corollary}
\begin{proof}
	Firstly, note that without loss of generality, we can always extend $k$ to a finite extension.
	
	Secondly, we can assume that $X$ is geometrically irreducible over $  K$: Let $L$ be the algebraic closure of $  K$ in the function field of $X$. Then $L/  K$ is a finite separable extension and $X$ is smooth and geometrically irreducible over $\Spec( L)$. By spreading out $X\rightarrow\Spec( L)$, we get a morphism with a smooth geometrically connected generic fibre. Since $\Br(X\otimes_{L}L^\sep)^{G_L}= \Br(X_{  K^\sep})^{G_{  K}}$, it suffices to prove the claim for this morphism. 
	
	Thirdly, if $X/ K$ is birational equivalent to a smooth projective variety $X^\prime$ over $ K$,  the claim will follow from Theorem \ref{big} since $V_\ell\Br(X_{ K^\sep})^{G_{K}}$ is a birational invariant. 
	
	Next, we will use de Jong's alteration theorem to reduce the question to this case. Let $X^\prime_{\Bar{ K}}\rightarrow X_{\bar{  K}}$ be an alteration such that $X^\prime_{\bar{ K}}$ admits a smooth projective compactification. We can assume that $X^\prime$ and the alteration are defined over a finite normal extension $L$ of $ K$. By shrinking $X$, we can assume that the alteration $X^\prime \rightarrow X_{L}$ is finite flat. By spreading out $X^\prime \rightarrow \Spec( L)$, we get a commutative diagram
	\begin{displaymath}
		\xymatrix{
			\sX^\prime \ar[r] \ar[d]^{f} &\sY^\prime\ar[d] \\
			\sX \ar[r] & \sY }
	\end{displaymath}
	where $\sY^\prime$ has function field $L$. By shrinking $\sY$, we can assume that $f$ is finite flat.  
	{Regarding $f$ as a }
	$\sY$-morphism and base{changing} to $\Spec( K^\sep)$, we get a {finite flat} $X^\prime _{  K^\sep}\rightarrow X_{  K^\sep}$.
	It induces a norm map $\Br(X^\prime_{  K^\sep})\rightarrow \Br(X_{ K^\sep})$ which is compatible with the norm map $\Br(\sX^\prime_{k^{\sep}})\rightarrow \Br(\sX_{k})$. Therefore we get a commutative diagram
	\begin{displaymath}
		\xymatrix{
			V_\ell\Br(\sX^\prime_{k^{\sep}})^{G_k} \ar[r] \ar[d] &V_\ell\Br(X^\prime_{ K^\sep})^{G_{ K}}\ar[d] \\
			V_\ell\Br(\sX_{k^{\sep}})^{G_k} \ar[r]  & V_\ell\Br(X_{ K^\sep})^{G_{ K}} }
	\end{displaymath}
	Since vertical maps are surjective, it suffices to show that the first row is surjective. Since the claim holds for $\sX^\prime \rightarrow \sY^\prime$, the natural map
	$$V_\ell\Br(\sX^\prime_{k^{\sep}})^{G_k}\longrightarrow V_\ell\Br(X^\prime_{L^\sep})^{G_L}$$
	is surjective. Thus, it suffices to show that there is a natural isomorphism 
	$$\Br(X^\prime_{ K^\sep})^{G_K}\cong\Br(X^\prime\otimes_{L} L^\sep)^{G_L}.$$
	Fix an algebraic closure $\bar{ K}$ of $ K$ and assume that $L$ and $L^\sep$ are contained in $\bar{ K}$. Since $X^\prime_{ K^\sep}=X^\prime \otimes_L L\otimes_{ K}  K^\sep$, the natural map $ L\otimes_{K} K^\sep\rightarrow L^\sep$ induces a map
	$$\Br(X^\prime \otimes_L L\otimes_{  K}  K^\sep)\longrightarrow \Br(X^\prime\otimes_{L} L^\sep).$$
	{Set $  K_1=L\cap  K^\sep$. By our assumption that $L/  K$ is a normal extension, $  K_1/  K$ is finite Galois.} Write $G=G_{  K}$ and $H=G_{  K_1}$. We have 
	$$L\otimes_  K  K^\sep=L\otimes_{  K_1}  K_1\otimes_{  K}  K^\sep=\prod_{\sigma\colon  K_1\hookrightarrow  K^\sep}L\otimes_{  K_1,\sigma}  K^\sep
	$$
	One can show that $L\otimes_{  K_1,\sigma}  K^\sep$ is a  closure of $L$ and admits a $H$-action (acting on $  K^\sep$). We denote it by $L^{s}_{\sigma}$. Then
	$$\Br(X^\prime \otimes_L L\otimes_{  K}  K^\sep)\cong \prod_{\sigma\colon  K_1\hookrightarrow  K^\sep}\Br(X^\prime _{L^{s}_{\sigma}}).
	$$
	Taking $H$-invariant, we get
	$$\Br(X^\prime \otimes_L L\otimes_{  K}  K^\sep)^H\cong \prod_{\sigma\colon  K_1\hookrightarrow  K^\sep}\Br(X^\prime _{L^{s}_{\sigma}})^H.
	$$
	Now $G/H$ acts as permutations on the right side, so the $G/H$-invariant of right side can be identified with the factor with the index $\sigma$ equal to the inclusion map $  K_1\subset  K^\sep$. The factor can be identified with $\Br(X^\prime_{L^\sep})^{G_L}$ through the natural isomorphism $L\otimes_{  K_1}  K^\sep\cong L^\sep$. It follows that there is a natural isomorphism
	$$\Br(X^\prime_{  K^\sep})^{G_  K}\cong\Br(X^\prime\otimes_{L} L^\sep)^{G_L}.$$
	This completes the proof.
\end{proof}

\subsection{The $\ell$-adic Tate module of the Tate--Shafarevich group}

The goal of this section is to study more closely 
the Tate--Shafarevich group in the  $\ell$-adic case. 
For this we need some general technical preparations.

\begin{lemma}\label[lemma]{size}
Let $G$ be a smooth connected commutative algebraic group over a finitely generated field $K$. 
Let $\ell\neq \Char(K)$ be a prime. 
Then for any $n\in\bbN$ the size of 
$$
\Hom(G[\ell^n],\bbQ_\ell/\bbZ_\ell)^{G_K}
$$
is bounded by a constant independent of $n$ and $\ell$.
\end{lemma}

\begin{proof}
Note that without loss of generality, 
we can replace $K$ by a finite extension. 
We will do this several times throughout the proof. 
By Chevalley's Theorem, the base change $G_{\overline{K}}$ 
of $G$ to an algebraic closure  is an extension of an abelian variety 
by a linear algebraic group.
After taking a finite field extension of $K$,
we can assume that this extension is defined over $K$. 
Hence we have an exact sequence of group schemes
$$
0\longrightarrow H\longrightarrow G\longrightarrow A\longrightarrow 0.
$$
Since $H$ is commutative, by {possibly further} extending $K$, we can assume that $H\cong \bbG_m^s\times \bbG_a^t$ over $K$.\\
Taking $\ell^n$-torsion and since $H$ is $\ell^n$-divisible, we have a short exact sequence of finite flat group schemes
$$
0\longrightarrow H[\ell^n]\longrightarrow G[\ell^n] \longrightarrow A[\ell^n]\longrightarrow 0.
$$
For later use, note that taking the dual, we obtain
$$
0\longrightarrow A[\ell^n]^\vee \longrightarrow G[\ell^n]^\vee \longrightarrow H[\ell^n]^\vee\longrightarrow 0.
$$
Since $H[\ell^n]=\bbG_m^s[\ell^n]$, 
it suffices to prove the claim for abelian varieties and the multiplicative group scheme. 
Let $\sY$ be an irreducible regular scheme of finite type over 
$\Spec\bbZ[{\ell^{-1}}]$ with function field $K$. 
Let $G$ be an abelian variety (resp.\ equal to $\bbG_m$) over $K$. 
By shrinking $\sY$, we can assume that $G$ extends to an abelian scheme $G_\sY$ (resp.\ $\bbG_{m,\sY}$). 
Let $s\in \sY$ be a closed point with finite residue field $  k(s)$.
Write $\sF$ for the \'etale sheaf represented by $G_\sY[\ell^n]$ on $\sY$. 
As $G_\sY[\ell^n]$ is finite \'etale over $\sY$, 
$\sF$ is a locally constant sheaf of $\bbZ/\ell^n$-modules on $\sY$. 
For the generic point $\eta_\sY$ of $\sY$ we have $\sF_{\bar{\eta}_\sY}=G[\ell^n]$. 
However, given that $\sY$ is connected, 
we have $\sF_{\overline{\eta}}\cong \sF_{\overline{s}}$, 
and hence $G[\ell^n]\cong G_s[\ell^n]$ as  $G_{  k(s)}$-modules. 
Through this isomorphism, there is an inclusion
$$
(G[\ell^n]^\vee)^{G_K} \subseteq (G_s[\ell^n]^\vee)^{G_{  k(s)}}.
$$
Therefore, it suffices to prove the claim for the case that $K$ is a finite field. 
In this case, 
$$
|(G[\ell^n]^\vee)^{G_K}|=|G[\ell^n]_{G_K}|=|G[\ell^n]^{G_K}|.
$$
Since $G[\ell^n]^{G_K}\subseteq G(K)$, 
the size of $(G[\ell^n]^\vee)^{G_K}=\Hom(G[\ell^n],\bbQ_\ell/\bbZ_\ell)^{G_K}$ 
is bounded by $|G(K)|$. 
To finish the proof, we note that the choice of $\sY/\bbZ[\ell^{-1}]$ 
works for all but finitely many $\ell\neq\Char(K)$. 
Thus we only have to repeat the above arguments and constructions 
for finitely many choices of $\sY$ and obtain finitely many constants 
of which we can simply take the maximum to obtain a constant as a bound 
for the size of $\Hom(G[\ell^n],\bbQ_\ell/\bbZ_\ell)^{G_K}$ that is independent of $n$ and $\ell\neq \Char(K)$.
\end{proof}

\begin{lemma}\label[lemma]{finite}
Let $R$ be a Henselian discrete valuation ring with residue field $k$ and quotient field $K$. 
Assume that $k$ is a finitely generated field. 
Let $R^{\sh}$ denote a strict Henselisation of $R$,  
$K^{\sh}$  its quotient field, and $k^{\sep}$ a separable closure of $k$. 
Let $A$ be an abelian variety over $K$ and $\ell$ be a prime different from $\Char(k)$. 
Set $G\colonequals \Gal(k^{\sep}/k)=\Gal(K^{\sh}/K)$. 
Then the size of 
$$
\H^1(K^{\sh}, A[\ell^n])^{G}
$$
is bounded by a constant independent of $n$. 
Moreover, $\H^1(K^{\sh}, A)^{G}(\ell)$ is a finite group.
\end{lemma}
\begin{proof}
Let $I$ denote the inertia group $\Gal(K^\sep/K^{\sh})$ and $I_1$ the wild inertia subgroup. 
Since $\H^1(I_1, A[\ell^n])= 0$ 
by \cite[Chap.\,I, Lem.\,2.18]{Milne_2006}, 
we have 
$$
\H^1(K^{\sh}, A[\ell^n])=(A[\ell^n](-1))_{I}
$$
by purity.
By the Weil pairing, its dual 
$\Hom((A[\ell^n](-1))_{I},\bbQ_\ell/\bbZ_\ell)=\Hom(A[\ell^n](-1),\bbQ_\ell/\bbZ_\ell)^I$ 
is canonically isomorphic to $A^t[\ell^n]^I$, 
where $A^t$ denotes the dual abelian variety of $A$. 
Let $\sA^t$ be a N\'eron model of $A^t/K$ and 
$s\colon\Spec(k) \rightarrow \sY=\Spec(R)$ the closed point. 
Then 
$$
A^t[\ell^n]^I=\sA^t_s[\ell^n].
$$
There is an exact sequence
$$
0\longrightarrow (\sA^t_s)^0 \longrightarrow \sA^t_s \longrightarrow \Phi(A^t)\longrightarrow 0,
$$
where $\Phi(A^t)$ denotes the Néron component group of $A^t$ 
(i.e., the component group of the Néron model $\sA^t$ of $A^t$),
which is a finite \'etale scheme over $k$. 
Taking first the dual and then $G$-invariants, 
we get an exact sequence
$$
0 \longrightarrow (\Phi(A^t)^\vee)^G \longrightarrow ((\sA^t_s[\ell^n])^\vee)^G \longrightarrow (((\sA^t_s)^0[\ell^n])^\vee)^G.
$$
Note that the size of $(\Phi(A^t)^\vee)^G$ is bounded by a constant independent of $n$ since $\Phi(A^t)$ is finite. 
Since
$\H^1(K^{\sh},A[\ell^n])^G \cong ((\sA^t_s[\ell^n])^\vee)^G$, 
it suffices to show that the size of $(((\sA^t_s)^0[\ell^n])^\vee)^G$ 
is bounded by a constant independent of $n$ as well, 
which follows directly from Lemma \ref{size}. 
This proves the first claim.

Next let us show that $\H^1(K^{\sh}, A)^{G}(\ell)$ is finite. 
By Kummer theory, we obtain an exact sequence
$$
0\longrightarrow A(K^{\sh})\otimes_{\bbZ}\bbQ_\ell/\bbZ_\ell \longrightarrow \H^1(K^{\sh}, A(\ell)) \longrightarrow \H^1(K^{\sh}, A)(\ell)\longrightarrow 0.
$$
Let $\sA$ be a N\'eron model of $A$ over $\sY$ and $\sA^0$ its identity component.
Again we have an exact sequence
$$
0\longrightarrow \sA^0(\sY)\longrightarrow \sA(\sY)\longrightarrow \Phi(A)\longrightarrow 0,
$$
where $\Phi(A)$ is the Néron component group of $A$. 
Since $\sA^0(\sY)$ is $\ell$-divisible and $\Phi(A)$ is finite, 
by tensoring with $\bbQ_\ell/\bbZ_\ell$, we obtain
$$
A(K^{\sh})\otimes_{\bbZ}\bbQ_\ell/\bbZ_\ell=0,
$$
where we used that $A(K^{\sh})=\sA(\sY)$.
It follows that 
$$
\H^1(K^{\sh}, A(\ell)) \cong \H^1(K^{\sh}, A)(\ell).
$$ 
Since 
$$\H^1(K^{\sh}, A(\ell))^G=\varinjlim_{n}\H^1(K^{\sh}, A[\ell^n])^{G},$$
is finite by the first claim, 
it follows that $\H^1(K^{\sh}, A)^G(\ell)$ is also finite.
\end{proof}
\noindent
We can now prove the main result of this subsection.

\begin{proposition} \label[proposition]{prop: sha}
Let $k$ be a finitely generated field.  
Denote by $k^{\sep}$ a separable closure of $k$ and by $G_k$ its absolute Galois group. 
Let $S$ be a smooth geometrically connected variety over $k$ 
with function field $ K$ and let $A$ be an abelian variety over $  K$. 
For any prime $\ell\neq\Char(k)$ we have the following statements:
\begin{enumerate}
\item 
Let $\sY\subseteq S$ be an open dense subscheme. 
Then the inclusion
$$
\Sha_{S_{k^{\sep}}}(A)^{G_k}(\ell) \longrightarrow \Sha_{\sY_{k^{\sep}}}(A)^{G_k}(\ell)
$$ 
has a finite cokernel.

\item If $A$ extends to an abelian scheme $\sA\longrightarrow  \sY$ 
on an open dense subscheme $\sY\subseteq S$ 
there is a canonical isomorphism
$$
\Sha_{\sY_{k^{\sep}}}(A)(\ell) \cong \H^1(\sY_{k^{\sep}},\sA)(\ell).
$$

\item
Write $  K^\prime$ for $  Kk^{\sep}$.
Then we have 
$$
V_\ell\Sha_{S_{k^{\sep}}}(A)^{G_k}=V_\ell \H^1(  K^\prime, A)^{G_k},
$$
which only depends on $A/  K$. 
We will denote it by $V_\ell\Sha_{  K^\prime}(A)^{G_k}$. 

\item
If $k$ is a finite field, we have
$$
V_\ell\Sha_{S}(A)=V_\ell \H^1(  K, A) \cong V_\ell\Sha_{  K^\prime}(A)^{G_k}.
$$
In this case, we will write $V_\ell\Sha'(A/K)$ for $V_\ell\Sha_{S}(A)$.
Note that this notation is consistent with~\cref{def: Shafarevich-Tate for abelian varieties}.
\end{enumerate}
\end{proposition}

\begin{proof}
It is easy to see that
$$
\Sha_{S_{k^{\sep}}}(A)^{G_k}=\ker\Big(\H^1(  K^\prime, A)^{G_k}\longrightarrow \prod_{s \in S^1} \H^1(K^{\sh}_s, A)^{G_{K^h_s}}\Big).
$$ 
By Lemma \ref{finite}, $\H^1(  K^{\sh}_s, A)^{G_{  K^h_s}}(\ell)$ is finite. Since $\sY$ omits only finitely many codimension $1$ points of $S$, it follows that the inclusion
$$
\Sha_{S_{k^{\sep}}}(A)^{G_k}(\ell) \longrightarrow \Sha_{\sY_{k^{\sep}}}(A)^{G_k}(\ell)
$$ 
has a finite cokernel and 
$V_\ell\Sha_{S_{k^{\sep}}}(A)^{G_k}=V_\ell \H^1(  K^\prime, A)^{G_k}$. 
In the case that $k$ is a finite field, by the same argument, the inclusion
$$
\Sha_{S}(A)(\ell)\longrightarrow \Sha_{\sY}(A)(\ell)
$${big}
has a finite cokernel and $V_\ell\Sha_{S}(A)=V_\ell \H^1(  K, A)$.

If $A/  K$ extends to an abelian scheme $\sA\rightarrow \sY$, we will show that
	\begin{equation}\label{puri}
		\H^1(\sY_{k^{\sep}}, \sA[\ell^n])=\ker\Big(\H^1(  K^\prime, A[\ell^n])\longrightarrow \prod_{s\in \sY_{k^{\sep}}^1}
		\H^1(  K_{s}^{\sh}, A[\ell^n])\Big).
	\end{equation}
	This will imply that the composition of 
	$$\H^1(\sY_{k^{\sep}},\sA)(\ell)\longrightarrow \H^1(\sY_{k^{\sep}}, j_*j^*\sA)\longrightarrow\Sha_{\sY_{k^{\sep}}}(A)(\ell)$$
	is surjective, where $j\colon\Spec(  K^\prime) \rightarrow \sY_{k^{\sep}}$ is the generic point. Since $\sA[\ell^n]$ is a locally constant sheaf, we have 
	$$\sA[\ell^n]\cong j_*j^*\sA[\ell^n].$$
	Thus, we get an exact sequence
	$$0\longrightarrow \H^1(\sY_{k^{\sep}},\sA[\ell^n])\longrightarrow \H^1(  K^\prime, A[\ell^n]) \longrightarrow \prod_{s \in \sY_{k^{\sep}}}\H^1(  K_s^{\sh}, A),$$
	and therefore
	$$\H^1(\sY_{k^{\sep}}, \sA[\ell^n])\subseteq\ker\Big(\H^1(  K^\prime, A[\ell^n])\longrightarrow \prod_{s\in \sY_{k^{\sep}}^1}
	\H^1(  K^{\sh}_s, A[\ell^n])\Big).$$
	Next, we will show that this inclusion is actually an equality.
	Let $x\in \H^1(  K^\prime, A[\ell^n])$.
	Since
	$$\H^1(  K^\prime, A[\ell^n])=\underset{V\subseteq \sY_{k^{\sep}}}{\varinjlim} \H^1(V,\sA[\ell^n]),$$
	there exists an open subscheme $V\subseteq \sY_{k^{\sep}}$ such that $x\in \H^1(V,\sA[\ell^n])$. Let $Y$ denote the reduced closed scheme $\sY_{k^{\sep}}-V$.  Since $\sA[\ell^n]$ is a locally constant sheaf, by the semi-purity \cite[\S 8]{Fujiwara_2002}, removing a closed subscheme of codimension $\geqslant 2$ from $\sY_{k^{\sep}}$ will not change $\H^1(\sY_{k^{\sep}}, \sA[\ell^n])$. So by shrinking $Y$, we can assume that $Y$ is regular and of codimension $1$ in $\sY_{k^{\sep}}$, then $\H^2_{Y}(\sY_{k^{\sep}}, \sA[\ell^n]) \cong \H^0(Y,\sA[\ell^n](-1))$ by purity \cite[\S 0]{Fujiwara_2002}. Let $s$ be a generic point of $Y$, by purity, we have $\H^1(  K^{\sh}_s,A[\ell^n])\cong (\sA[\ell^n](-1))_{\bar{s}}$. Thus, we get an exact sequence
	$$0\longrightarrow \H^1(\sY_{k^{\sep}}, \sA[\ell^n])\longrightarrow \H^1(V,\sA[\ell^n]) \longrightarrow \prod_{s\in Y^0} \H^1(  K^{\sh}_s, A[\ell^n]).
	$$
	Therefore, if $x$ lies in 
	$$\ker\Big(\H^1(  K^\prime, A[\ell^n])\rightarrow \prod_{s\in \sY_{k^{\sep}}^1}\H^1(  K^{\sh}_s, A[\ell^n])\Big),$$
	it also lies in $\H^1(\sY_{k^{\sep}}, \sA[\ell^n])$. This proves (\ref{puri}).
	By \cite[Thm. 3.3]{Keller_2016}, we have $\sA\cong j_*j^*\sA$. Thus 
	$$
	\H^1(\sY_{k^{\sep}}, \sA)=\ker (\H^1(  K^\prime, A)\rightarrow \prod_{s\in \sY_{k^{\sep}}^1}\H^1(  K^{\sh}_s, A)).
	$$
	It follows that $\H^1(\sY_{k^{\sep}}, \sA)(\ell)\subseteq \Sha_{\sY_{k^{\sep}}}(A)(\ell)$. Since the natural map $\H^1(\sY_{k^{\sep}}, \sA(\ell))\rightarrow\Sha_{\sY_{k^{\sep}}}(A)(\ell)$ factors through $\H^1(\sY_{k^{\sep}}, \sA)(\ell)$ and $\H^1(\sY_{k^{\sep}}, \sA(\ell))\rightarrow \Sha_{\sY_{k^{\sep}}}(A)(\ell)$ is surjective, thus
	$$\H^1(\sY_{k^{\sep}}, \sA)(\ell)=\Sha_{\sY_{k^{\sep}}}(A)(\ell).$$
	By the same argument, one can show that
	$$\H^1(\sY, \sA)(\ell)=\Sha_{\sY}(A)(\ell).$$
	
	It remains to prove that $V_\ell\Sha'(A/K)\cong V_\ell\Sha_{  K^\prime}(A)^{G_k}$ when $k$ is a finite field. Firstly, we will show that the natural map
	$$V_\ell\Sha'(A/K)\longrightarrow V_\ell\Sha_{  K^\prime}(A)^{G_k}$$
	is surjective. Since there exists an abelian variety $B/  K$ such that $A\times B$ is isogenous to  $\Pic^0_{X/  K}$ for some smooth projective geometrically connected curve $X/  K$, we can assume that $A=\Pic^0_{X/  K}$. One can spread out $A/  K$ to a smooth projective relative curve $\rho\colon\sX\rightarrow S$ where $S$ is smooth and geometrically connected over $k$. In the following, we will show that there is commutative diagram with exact second row
	\begin{displaymath}
		\xymatrix{V_\ell\Br(\sX)\ar[r]\ar[d]& V_\ell\Sha'(\Pic^0_{X/  K}/K)\ar[d]\\
			V_\ell\Br(\sX_{k^{\sep}})^{G_k}\ar[r]& V_\ell\Sha_{  K^\prime}(\Pic^0_{X/  K})^{G_k}\ar[r]& 0}
	\end{displaymath}
	Then the surjectivity of the second column will follow from the surjectivity of the first column.
	By the Leray spectral sequence
	$$E^{p,q}_2=\H^p(S,\R^q\rho_*\bbG_m)\Rightarrow \H^{p+q}(\sX,\bbG_m)
	$$
	we get an exact sequence
	$$\H^2(S,\mathbb{G}_m)\longrightarrow \ker(\H^2(\sX,\mathbb{G}_m)
	\rightarrow \H^0(S,\R^2\rho_*\mathbb{G}_m))\longrightarrow \H^1(S,\R^1\rho_*\mathbb{G}_m).
	$$
	By Lemma \ref{picard-brauer}\,(i) above, $\H^0(S,\R^2\rho_*\mathbb{G}_m)(\ell)\cong \Br(X_{  K^\sep})^{G_{  K}}(\ell)$. Since $X$ is a smooth projective curve over $  K$, we have $\Br(X_{  K^\sep})=0$. Thus, we get a natural map
	$$
	V_\ell\Br(\sX)\longrightarrow  V_\ell \H^1(S, \R^1\rho_*\bbG_m).
	$$
	Let $j\colon\Spec(  K)\rightarrow S$ be the generic point. By Lemma \ref{picard-brauer}\,(i), we have 
	$$\R^1\rho_*\bbG_m\cong j_*j^*\R^1\rho_*\bbG_m.$$
	Thus, we have
	$$\H^1(S, \R^1\rho_*\bbG_m)(\ell)\cong \H^1(S, j_*\Pic_{X/  K})(\ell)\cong \H^1(S,j_*\Pic^0_{X/  K})(\ell)\cong \Sha_{S}(\Pic^0_{X/  K})(\ell).$$
	Therefore,
	$$V_\ell \H^1(S, \R^1\rho_*\bbG_m)\cong V_\ell\Sha_{S}(\Pic^0_{X/  K})$$
	By compositions, we get a natural map
	$$
	V_\ell\Br(\sX)\longrightarrow V_\ell\Sha'(\Pic^0_{X/  K}/K).
	$$
	By Theorem \ref{big}, there is a surjective map compatible with the above map
	
	$$
	V_\ell\Br(\sX_{k^{\sep}})^{G_k}\longrightarrow V_\ell\Sha_{  K^\prime}(\Pic^0_{X/  K})^{G_k}.
	$$
	It suffices to show that $V_\ell\Br(\sX)\rightarrow V_\ell\Br(\sX_{k^{\sep}})^{G_k}$ is  surjective. By Lemma \ref{ssplit}, the natural map 
	$$\H^2(\sX_{k^{\sep}},\bbQ_\ell(1))^{G_k}\longrightarrow V_\ell\Br(\sX_{k^{\sep}})^{G_k}$$
	is surjective. Since $\H^2(\sX,\bbQ_\ell(1))\rightarrow V_\ell\Br(\sX)$ is surjective by definition, it suffices to show that
	$$\H^2(\sX,\bbQ_\ell(1))\longrightarrow \H^2(\sX_{k^{\sep}},\bbQ_\ell(1))^{G_k}$$
	is surjective. Consider the spectral sequence
	$$\H^p(G_k,\H^q(\sX_{k^{\sep}},\mu_{\ell^n}))\Rightarrow \H^{p+q}(\sX,\mu_{\ell^n})$$
	Since $\H^q(G_k,-)$ vanishes when $q\geqslant 2$, we get a surjection
	$$\H^2(\sX,\mu_{\ell^n})\longrightarrow \H^2(\sX_{k^{\sep}},\mu_{\ell^n})^{G_k}.$$
	Taking {the limit over $n$}, we get 
	$$\H^2(\sX,\bbQ_\ell(1))\longrightarrow \H^2(\sX_{k^{\sep}},\bbQ_\ell(1))^{G_k}\longrightarrow 0.$$
	Next, we will show that the natural map
	$$V_\ell\Sha'(A/K)\longrightarrow V_\ell\Sha_{  K^\prime}(A)^{G_k}$$
	is injective. Consider the exact sequence
	$$0\longrightarrow \H^1(G_k,A(  K^\prime))\longrightarrow \H^1(  K,A)\longrightarrow \H^1(  K^\prime, A).$$
	Let $\Tr_{  K/k}(A)$ be the $  K/k$ trace of $A$. By the Lang--N\'eron theorem, $A(  K^\prime)/\Tr_{  K/k}(A)(k^{\sep})$ is finitely generated. Since $\H^1(G_k,\Tr_{  K/k}(A))=0$ by Lang's theorem, we have that $\H^1(G_k,A(  K^\prime))$ is finite. It follows that the natural map
	$$V_\ell \H^1(  K,A) \longrightarrow V_\ell \H^1(  K^\prime, A)$$
	is injective. So $V_\ell\Sha'(\Pic^0_{X/  K}/K)\longrightarrow V_\ell\Sha_{  K^\prime}(\Pic^0_{X/  K})^{G_k}$ is also injective. This completes the proof.	
\end{proof}

\begin{corollary}\label{kktrace}
  Let $k$ be a finitely generated field. Let $K$ be the function field of a smooth geometrically connected variety over $k$. Let $A$ be an abelian variety over $  K$. Let $B=\Tr_{  K/k}(A)$ denote the $  K/k$ trace of $A$. Let $\ell\neq \Char(k)$ be a prime. Then there is a canonical exact sequence
	$$
	0\longrightarrow V_\ell\Sha'(B/k) \longrightarrow V_\ell\Sha'(A/K) \longrightarrow V_\ell\Sha_{  K^\prime}(A)^{G_k}\longrightarrow 0.
	$$
\end{corollary}
\begin{proof}
	Let's first prove the statement for the case $k$ is of characteristic $p>0$. Consider the exact sequence
	$$0\longrightarrow \H^1(G_k, A(  K^\prime))\longrightarrow \H^1(  K,A)\longrightarrow \H^1(  K^\prime, A),$$
	we will show that $V_\ell \H^1(G_k, A(  K^\prime))$ is naturally isomorphic to $V_\ell \H^1(k,B)$. This will imply the left exactness of the sequence in the claim. By the Lang--N\'eron theorem, $A(  K^\prime)/B(k^{\sep})$ is finitely generated, so $\H^1(G_k, A(  K^\prime)/B(k^{\sep}))$ is of finite exponent. By the exact sequence
	$$0\longrightarrow B(k^{\sep})\longrightarrow A(  K^\prime)\longrightarrow A(  K^\prime)/B(k^{\sep})\longrightarrow 0,$$
	we get 
	$$V_\ell  \H^1(G_k, A(  K^\prime))=V_\ell \H^1(G_k, B(k^{\sep}))=V_\ell\Sha'(B/k),$$
	where the last equality follows from \cref{prop: sha}. It follows that
	$$
	0\longrightarrow V_\ell\Sha'(B/k) \longrightarrow V_\ell\Sha'(A) \longrightarrow V_\ell\Sha_{  K^\prime}(A)^{G_k}
	$$
	is exact.
	
	To show the surjectivity of the last arrow, we can assume $A=\Pic^0_{X/  K}$ where $X$ is a smooth projective geometrically connected curve over $  K$. Spreading out $X/  K$, we get a smooth projective relative curve $\rho\colon\sX \rightarrow S$. Since $\Br(X_{  K^\sep})=0$, by Theorem \ref{big}, we get a surjection
	$$ \Br(\sX_{k^{\sep}})^{G_k}\longrightarrow V_\ell\Sha_{  K^\prime}(A)^{G_k}.$$
	Let $k_0$ be the algebraic closure in $k$ of the prime field of $k$. Then there exists a smooth geometrically connected variety $\cZ$ over $k_0$ with function field $k$. By shrinking $\cZ$, we can spread out $\rho$ to a smooth projective morphism $\rho_{\cZ}\colon\sX_{\cZ}\rightarrow S_{\cZ}$. By shrinking $\cZ$ further, we can assume that $S_{\cZ}$ is smooth over $\cZ$. Using Theorem \ref{big} for $\rho_{\cZ}$, we get a surjection 
	$$V_\ell\Br(\sX_{\cZ_{k_0^\sep}})^{G_{k_0}}\longrightarrow 
	V_\ell\Sha'(A/Kk_0^\sep)^{G_{k_0}}.$$
	It remains to show that the natural map
	$$ V_\ell\Br(\sX_{\cZ_{k_0^\sep}})^{G_{k_0}} \longrightarrow \Br(\sX_{k^{\sep}})^{G_k}$$
	is surjective. Since $\sX$ can be identified with the generic fibre of $\sX_{\cZ}\rightarrow \cZ$, then the claim follows from Corollary \ref{hardsur}. This proves the positive characteristic case.

 Now assume that $\Char(k)=0$. Let's show the surjectivity of the last map first. This is same as the positive case. Since $k$ is of characteristic $0$, by Hironaka's theorem of resolution of singularity, we can assume that $\sX$ and $S$ are both smooth projective geometrically connected varieties over $k$. Then the surjectivity follows from the fact that $\Br_{\nr}(K(\sX){/k}) \rightarrow  \Br(\sX_{k^\sep})^{G_k}$ is surjective (cf.\,Theorem~\hyperlink{bigtheorem}{D}). This proves the surjectivity of the last map. For any discrete valuation $v$ of $K$. Fix an extension of $v$ on $\bar{K}$. the strictly Henselian field at $v$ of $K$ contains the strictly Henselian field of $k$ at $v|k$. So there is a natural map $\H^1(Kk^{\sh}_v,A) \rightarrow  \H^1(K^{\sh}_v, A)$. The composition of the natural maps
 $$\H^1(k_v^{\sh}, B) \rightarrow  \H^1(K^\prime/Kk^{\sh}_v, A(K^\prime))\rightarrow  \H^1(Kk^{\sh}_v,A) \rightarrow  \H^1(K^{\sh}_v, A)$$
 is compatible with $\H^1(k, B)\rightarrow  \H^1(K,A)$. This shows that 
 $$
 V_\ell\Sha'(B/k) \subseteq \ker( V_\ell\Sha'(A/K) \longrightarrow V_\ell\Sha_{  K^\prime}(A)^{G_k})
 $$
 To show that they are actually same, we compute their dimensions as $\bbQ_\ell$-linear spaces. Let $\rho\colon\sX\rightarrow S$.
be a projective surjective $k$-morphism between smooth projective geometrically connected $k$-varieties with a generic fibre $X$ identified with $A^t/K$. Now we apply the first part of Theorem~\hyperlink{bigthm}{E} and Theorem~\hyperlink{bigtheorem}{D} to $\rho$.  We get
\begin{align*}
\dim(V_\ell\Sha_{K^\prime}(A)^{G_k})&=\dim (V_\ell\Br(\sX_{k^\sep})^{G_k})-\dim (V_\ell\Br(X_{K^\sep})^{G_K})-\dim (V_\ell\Br(S_{k^\sep})^{G_k}),\\
\dim(V_\ell\Sha'(A/K))&=\dim (V_\ell{\Br_{\nr}}(K(X){/k}))-\dim (V_\ell\Br(X_{K^\sep})^{G_K})-\dim (V_\ell{\Br_{\nr}}(K/k)),\\
\dim (V_\ell\Br(\sX_{k^\sep})^{G_k})&=\dim (V_\ell{\Br_{\nr}}(K(\sX){/k}))-\dim (V_\ell\Sha'(\Pic^0_{\sX/k}/k))-\dim (V_\ell{\Br_{\nr}}(k)),\\
\dim (V_\ell\Br(S_{k^\sep})^{G_k})&=\dim (V_\ell{\Br_{\nr}}(K/k))- \dim (V_\ell\Sha'(\Pic^0_{S/k}/k))-\dim (V_\ell
{\Br_{\nr}}(k)).
\end{align*}
Note that $K(\sX)=K(X)$. If we plug the last two equations into the first one and subtract the second one, we get:
$$\dim(V_\ell\Sha'(A/K))=\dim(V_\ell\Sha_{K^\prime}(A)^{G_k})+\dim (V_\ell\Sha'(\Pic^0_{\sX/k}))-\dim (V_\ell\Sha'(\Pic^0_{S/k}/k))
$$
It remains to prove
$$
\dim(V_\ell\Sha'(B/k))=
\dim(V_\ell\Sha'(\Pic^0_{\sX/k}/k))-\dim (V_\ell\Sha'(\Pic^0_{S/k}/k)).
$$
This follows from the next lemma.
\end{proof}
\begin{lemma}
Let $k$ be a field of characteristic $0$. Let $\rho\colon\sX\rightarrow S$
be a projective surjective $k$-morphism between smooth projective geometrically connected $k$-varieties with a generic fibre $X$. Assume that $X$ is an abelian variety. Let $A$ denote $\Pic^0_{X/K}$ and $B$ denote $\Tr_{  K/k}(A)$. Then there is an exact sequence of abelian varieties over $k$
$$0\rightarrow  \Pic^0_{S/k} \dashrightarrow \Pic^0_{\sX/k} \rightarrow  B \rightarrow  0,$$
where the dash arrow means that the map is an isogenous to the identity component of the kernel of the next arrow.
\end{lemma}
\begin{proof}
This can be proved by the same method as the proof of~\cite[Prop.\,4.4, p.\,302]{Ulm}. By the definition of Picard functor, there is a natural $K$-morphism of abelian varieties 
$$\Pic^0_{\sX/k}\times_k K \rightarrow  \Pic^0_{X/K}.$$
By definition of the $K/k$-trace, this factor through $B_K$ and there is a $k$-morphism $\Pic^0_{\sX/k} \rightarrow  B$. By looking at $\bar{k}$-points, this map has a cokernel finitely generated over $\bbZ$ since the Neron--Severi group of $\sX_{\bar{k}}$ is finitely generated and $\Pic(\sX_{\bar{k}}) \rightarrow  \Pic(X_{K\bar{k}})$ is surjective. It remains to show that the kernel of $\Pic^0(\sX_{\bar{k}}) \rightarrow  \Pic(X_{K\bar{k}})$ contains $\Pic^0(S_{\bar{k}})$
as a subgroup with finite index. In fact, it is enough to prove that 
$$ \Pic(S_{\bar{k}})\rightarrow  \ker(\Pic(\sX_{\bar{k}}) \rightarrow  \Pic(X_{K\bar{k}}))$$
is injective and has a cokernel finitely generated over $\bbZ$. The injectivity follows from the fact that $\rho$ have a section after removing a close subset of codimension at least 2. Let $U\subseteq S$ be an open dense subset such that $\rho$ is smooth over $U$. By the proof of \cref{bira2}, 
$$
0 \rightarrow  \Pic(U_{\bar{k}}) \rightarrow  \Pic(\rho^{-1}(U_{\bar{k}})) \rightarrow  \Pic(X_{K\bar{k}})) \rightarrow  0
$$
is exact. This shows that the pullback of $\Pic(S_{\bar{k}})$ and finitely many divisors contained in $\rho^{-1}(S_{\bar{k}}-U_{\bar{k}})$ generate $\ker(\Pic(\sX_{\bar{k}}) \rightarrow  \Pic(X_{K\bar{k}}))$. This shows that 
$$ \Pic(S_{\bar{k}})\rightarrow  \ker(\Pic(\sX_{\bar{k}}) \rightarrow  \Pic(X_{K\bar{k}}))$$
has a cokernel finitely generated over $\bbZ$. Thus,
$$ 
\Pic^0(S_{\bar{k}})\rightarrow  \ker(\Pic^0(\sX_{\bar{k}}) \rightarrow  \Pic^0(X_{K\bar{k}}))
$$
also has a cokernel finitely generated over $\bbZ$. Since $\ker(\Pic^0(\sX_{\bar{k}}) \rightarrow  \Pic^0(X_{K\bar{k}}))$ is same as the set of $\bar{k}$-points of $\ker(\Pic^0_{\sX/k} \rightarrow  B)$ which is an extension of finite group by an abelian subvariety $C$, the natural map $\Pic^0_{S/k} \rightarrow  C $ must be surjective. It is actually an isomorphism since $k$ has characteristic $0$.
\end{proof}

\subsection{Fibrations of Brauer groups and $p$-torsion}

In this section, we are interested in a $p$-adic analogue of \cref{big}.
We consider the situation that $k=\bbF_q$, $q=p^r$, is a finite field, 
and denote by $K_0$ the fraction field of $W(k)$.

We first want to justify the definition of the Tate--Shafarevich group from \cref{def: Shafarevich-Tate for abelian varieties} in this context. 
Step by step we will argue that we can find a ``nice enough'' 
model of a finitely generated field $K$ over $k$. 

The first statement concerns the existence of models for a finitely generated field $K/k$ whose (unramified) Brauer groups compare well enough: 

\begin{lemma}\label[lemma]{lem: S for K}
For any finitely generated field $K$ over a field $k$, there exists a smooth irreducible (hence integral) variety $S/k$ with function field $K$ such that
the natural injection
$$
{\Br^{\nr}}(K/k) \rightarrow \Br(S)
$$
has a cokernel of finite exponent.
\end{lemma}
\begin{proof}
There exists a projective integral variety $S'/k$ with $K(S') \cong K$. 
Furthermore, by de Jong's alteration theorem, 
there exists a finite separable field extension $L/K$ and a smooth projective variety $X/k$ with $  K(X) \cong L$ \cite{deJong_1996}.  
According to Raynaud--Gruson's flattening theorem \cite{RaynaudGruson_1971}, 
we may assume that there exists a birational morphism $X' \to X$ 
and a flat morphism $X' \rightarrow S'$. 
Taking normalisations, we may further assume that $X'$ and $S'$ are normal 
and $X' \to S'$  is finite. 
Since for every point $x \in X'$ of codimension $1$, its image under $X' \to S'$ is also of codimension $1$, there is a well defined morphism $X'_\reg \to S'_\reg$ and hence a commutative diagram
\[
\begin{tikzcd}
	\Br(S'_\reg) \ar[r]\ar[d,hookrightarrow] & \Br((X')_\reg) \ar[d,hookrightarrow]\\
	\Br(K) \ar[r] & \Br(L).
\end{tikzcd}
\]
Because any codimension-$1$ point of $X$ is the image of a codimension-$1$ point in $X'$ under the birational map $X' \to X$ and because $X'$ is regular in codimension $1$, 
one has $\Br(X) = \Br(X'_\reg) = {\Br^{\nr}}(L/k)$. 
Hence the above commutative diagram reads
\[
\begin{tikzcd}
	\Br(S'_\reg) \ar[r]\ar[d,hookrightarrow] & {\Br^{\nr}}(L/k) \ar[d,hookrightarrow]\\
	\Br(K) \ar[r] & \Br(L).
\end{tikzcd}
\]
The composite with the corestriction (which exists since $L/K$ is finite separable)
\[
\begin{tikzcd}
	{\Br^{\nr}}(L/k) \ar[r]\ar[d,hookrightarrow] & {\Br^{\nr}}(K/k) \ar[d,hookrightarrow]\\
	\Br(L) \ar[r] & \Br(K)
\end{tikzcd}
\]
is multiplication by $n = [L:K] < \infty$. 
In other words, $[n]\colon \Br(K) \to \Br(K)$ maps $\Br((S')_\reg)$ into ${\Br^{\nr}}(K/k)$. 
Hence the cokernel $\Br(S'_\reg)/{\Br^{\nr}}(K/k)$ is killed by $n$. 
Taking $S = S'_\reg$ gives the desired smooth integral $k$-variety.
\end{proof}

We first consider the case of a relative curve $X/K$,  we need to relate the unramified Brauer group and the Tate--Shafarevich group. 

\begin{lemma}\label[lemma]{lem: Br to Sha}
Let $X/K$ be a smooth projective geometrically connected curve over a finitely generated field $K/k$. 
The natural map 
$$
\Br(X) \rightarrow \H^1(K,\Pic_{X/K})
$$
sends
${\Br^{\nr}}(  K(X)/k)$ into $\Sha'(\Pic_{X/K}/K)$.
\end{lemma}

\begin{proof}
Let $v$ be a divisorial valuation on $  K(X)/k$. 
By Raynaud--Gruson's flattening theorem \cite{RaynaudGruson_1971}, 
$v|_K$ is also a divisorial valuation on $K$. 
Consider the induced morphism $\Spec\cO_{v|_{  K(X)}} \to \Spec\cO_{v|_K}$. 
By resolution of singularities of surfaces, 
we can assume that there exists a smooth projective relative curve 
$\rho\colon \frX \to \Spec\cO_{v|_K}$ such that $\frX$ is regular 
and $v$ extends to $\frX$. 
Hence we get a morphism 
$\Br(\frX) \to \H^1(\Spec\cO_{v|_K},\R^1\rho_*\bbG_m)$ (from $\ker(E^2 \rightarrow E_2^{0,2})\to E_2^{1,1}$ of the Leray spectral sequence for $\rho_*\bbG_m$ and the fact $E^{0,2}_2=0$ by Artin's theorem).
With $S_v \colonequals \Spec\cO_{v|_K}$ consider the commutative diagram
\[
\begin{tikzcd}
	\H^1(S_v,\R^1\rho_*\bbG_m) \ar[r]\ar[d] & \H^1(S_v^\sh,\R^1\rho_*\bbG_m)\ar[d]\\
	\H^1(K,\R^1\rho_*\bbG_m) \ar[r] & \H^1(K_v^\sh,\R^1\rho_*\bbG_m).
\end{tikzcd}
\]
This gives a map $\Br(\frX) \to \ker(\H^1(K,\Pic_{X/K}) \to \H^1(K_v^\sh,\Pic_{X/K}))$. 
The image of $\Br(\frX)$ in the upper left entry goes to $0$ on the right because $S_v^\sh$ is a strictly Henselian local ring, hence has trivial \'etale cohomology with degree greater than $1$.
Taking the intersection over all $v$ results in a map 
\[
	\bigcap\Br(\frX) \to \ker\Big(\H^1(K,\Pic_{X/K}) \to \prod_{v \in div\Val(K)}\H^1(K_v^\sh,\Pic_{X/K})\Big).
\]
where the target is by \cref{def: Shafarevich-Tate for abelian varieties} equal to  $\Sha'(\Pic_{X/K})/K$.
Hence we obtain the desired morphism ${\Br^{\nr}}(  K(X)/k) \to \Sha(\Pic_{X/K}/K)$.
\end{proof}

\begin{lemma} \label[lemma]{claim 2}
	Let $p \in X(K)$ and $p\colon \Spec(K) \to X$ the corresponding morphism. This induces $p^*\colon \Br(X) \to \Br(K)$. For all $x \in S^1$, $p$ extends uniquely to $p_x\colon \Spec(\sO_{S,x}) \to \sX_{\Spec(\sO_{S,x})}$ and $p_x^*\colon \Br(\sX_{\Spec(\sO_{S,x})}) \to \Br(\Spec(\sO_{S,x}))$ is compatible with $p^*$.
\end{lemma}
\begin{proof}
	This follows from functoriality of the (unramified) Brauer group.
\end{proof}

\begin{lemma}\label[lemma]{claim 3}
	Let $R$ be an excellent DVR, set $S=\Spec(R)$ and $K=\mathrm{Frac}(R)$. 
Let $\rho\colon \sX\rightarrow S$, with $\sX$ regular, be a projective flat morphism with smooth projective geometrically connected generic fibre $X\rightarrow \Spec(K)$. Assume that $X(K)$ is not empty.
	\begin{enumerate}
		\item Let $j\colon \Spec(K) \to S$ be the inclusion of the generic point of $S$. Then then the adjunction $\id \to j_*j^*$ yields a sequence
		\[
		0 \to \sF \to \R^1\rho_*{\bbG_m} \to j_*j^*\R^1\rho_*{\bbG_m} \to 0
		\]
		of \'etale sheaves on $S$.
		The sheaf $\sF$ is supported on the closed point of $S$ and $\sF = 0$ if $\rho$ is smooth. 
		\item 
		The induced morphism
		\[
		\H^1(S, \R^1\rho_*{\bbG_m}) \to \ker\big(\H^1(K,\Pic(X_{K^\sep})) \to \H^1(K^\sh,\Pic(X_{K^\sep}))\big).
		\]
		is an isomorphism up to finite exponent (and is an isomorphism if $\rho$ is smooth). 
	\end{enumerate}
\end{lemma}
\begin{proof}
(i) For the existence of the short exact sequence, we only have to show surjectivity. Consider the stalk at the closed point $x$. we have
  $$(\R^1\rho_*{\bbG_m})_{\bar{x}}=\Pic(\sX_{S^\sh}),$$
  $$(j_*j^*\R^1\rho_*{\bbG_m})_{\bar{x}}=\Pic_{X/K}(K^{sh}).$$ 
  We have $\Pic_{X/K}(K^{sh})=\Pic(X_{K^\sh})$ since $X(K)\neq \emptyset$. This proves (i).
		
(ii) Let $S^\h$ be the Henselisation and $S^\sh$ be the strict Henselisation of $S$ at $x$. Consider
		\[
		\begin{tikzcd}
			\H^1(x,\sF) \ar[d,equal]\ar[r] & \H^1(S,\R^1\rho_*{\bbG_m}) \ar[d]\ar[r,"a"] & \H^1(S,j_*j^*\R^1\rho_*{\bbG_m}) \ar[d]\ar[r] & \H^2(x,\sF) \ar[d,equal]\\ 
			\H^1(x,\sF) \ar[r] & \H^1(S^\h,\R^1\rho_*{\bbG_m}) \ar[r,"b"] & \H^1(S^\h,j_*j^*\R^1\rho_*{\bbG_m}) \ar[r] & \H^2(x,\sF).
		\end{tikzcd}
		\]
		Since $\sF$ is supported on $x$, hence $\H^i(S,\sF) = \H^i(x,\sF)$.
		
		Now $\H^1(S^\h,\R^1\rho_*{\bbG_m}) \cong \H^1(k(x),\Pic(\sX_{S^\sh}))$ and $\H^1(S^\h,j_*j^*\R^1\rho_*{\bbG_m}) \cong \H^1(k(x),\Pic(X_{K^\sh}))$. Consider the exact sequence of $G_{k(x)}$-modules
		\[
		0 \to \sF|_{\overline{x}} \to \Pic(\sX_{S^\sh}) \to \Pic(X_{K^\sh}) \to 0.
		\]
		By the lemma below, this sequence is almost split as continuous $G_{k(x)}$-modules in the sense of \cite[Lem. 3.19]{QinBrauer2}. 
		Thus by taking $\H^1(k(x),-)$, we get the second row and that $b$ is an isomorphism up to finite exponent. 
		The same is true for  $a$ since $\coker{a} \subseteq \coker{b}$ and $\H^1(k(x),\sF|_{\overline{x}})$ is of finite exponent. 
		By the Leray spectral sequence for the morphism $j\colon \Spec(K) \rightarrow S$ and the sheaf $j^*\R^1\rho_*{\bbG_m}$ on $\Spec(K)$, there is an exact sequence
  $$0\rightarrow \H^1(S,j_*j^*\R^1\rho_*{\bbG_m}) \rightarrow \H^1(\Spec(K), j^*\R^1\rho_*{\bbG_m}) \rightarrow\H^0(S, \R^1j_*j^*\R^1\rho_*{\bbG_m}).
  $$
Since $(\R^1j_*j^*\R^1\rho_*{\bbG_m})_{\bar{x}}=\H^1(K^\sh,\Pic(X_{K^\sep}))$, hence 
\[\H^1(S,j_*j^*\R^1\rho_*{\bbG_m})=	\ker\big(\H^1(K,\Pic(X_{K^\sep})) \rightarrow\H^1(K^\sh,\Pic(X_{K^\sep}))\big). \qedhere\]
\end{proof}
\begin{lemma}
Let $\rho\colon \sX\longrightarrow S=\Spec(R)$ be a proper flat morphism with a smooth and geometrically connected generic fibre $X$, where $R$ is a Henselian DVR with residue field $k$ and function field $K$. Assume that $\sX$ is regular. Then the exact sequence of continuous $G_k$-modules
$$0\rightarrow D_{\sX} \rightarrow \Pic(\sX_{R^{\sh}}) \rightarrow \Pic(X_{K^{\sh}})\rightarrow 0
$$
is almost split ( cf. \cite[Lem. 3.19]{QinBrauer2} for the definition), where $D_\sX$ is defined as the kernel of the map on Picard groups.
\end{lemma}
\begin{proof}
By de Jong's alteration theorem \cite{deJong_1996}, there exists a strictly semi-stable morphism $\sX^\prime \rightarrow \Spec(R^\prime)$ with a geometrically connected generic fibre and an alteration $f\colon \sX^\prime \rightarrow \sX$ with $R \rightarrow R^\prime$ a finite flat morphism of discrete valuation rings. By the proof of~\cite[Prop.\,3.20]{QinBrauer2} (where it was shown that the middle column of the commutative diagram in ~\cite[p. 25]{QinBrauer2} has a kernel of finite exponent 
and the second row is almost split ),
the claim holds for $\sX^\prime \rightarrow \Spec(R^\prime)$. The morphism
$\sX^\prime_{R^{\prime \sh}} \rightarrow \sX_{R^{\sh}}$ induces a pullback $\Pic(\sX^\prime_{R^{\prime \sh}}) \rightarrow \Pic(\sX_{R^{\sh}})$ that has a kernel of finite exponent. In fact, Raynaud--Gruson's flattening theorem gives a commutative diagrams with birational columns and a flat first row
\begin{displaymath}
\xymatrix{\widetilde{\sX^\prime_{R^{\prime \sh}}} \ar[r]\ar[d]& \widetilde{\sX_{R^{\sh}}}\ar[d]\\
\sX^\prime_{R^{\prime \sh}} \ar[r]&\sX_{R^{\sh}}
}
\end{displaymath}
The first row is proper and generically finite, hence it is finite flat. So the pullback on Picard groups induced by the first row has a kernel of finite exponent. Since the vertical maps induce injections on Picard groups (birational proper morphism admits a section after removing a codimension $\geqslant 2$ closed subset from the target), the second row induces a map with a kernel of finite exponent on Picard groups. Thus, the exact sequence above is almost split
as continuous $G_l$-modules, where $l$ is the residue field of $R^\prime$.
Let $P$ be a $G_l$-equivariant almost projection $\Pic(\sX_{R^{\sh}}) \rightarrow D_\sX$. Then $\sum_{\sigma \in G_K/G_l}P^{\sigma}$ is an almost $G_k$-equivariant projection\footnote{i.e., projection in the isogeny category of abelian groups with continuous $G_k$-action}.
This proves that the exact sequence is split up to finite exponent.
\end{proof}

\begin{lemma} \label[lemma]{claim 4}
Let $R$ be an excellent DVR with fraction field $K$. Let $X/K$ be a smooth projective geometrically integral curve. Assume that $X(K) \neq \emptyset$. Let $\rho\colon \sX \to S$ be a projective regular model of $X/K$. Then there are exact sequences
	\[
	0 \to \Br(K) \to \Br(X) \to \H^1(K,\Pic(X_{K^\sep})) \to 0
	\]
	and
	\[
	0 \to \Br(S) \to \Br(\sX) \to \H^1(S,\R^1\rho_*{\bbG_m}) \to 0
	\]
\end{lemma}
\begin{proof}
	Consider the Leray spectral sequence for $\rho\colon \sX \to S$
	\[
		E_2^{pq} = \H^p(S,\R^q\rho_*{\bbG_m}) \Longrightarrow \H^{p+q}(\sX,{\bbG_m}).
	\]
	It gives rise to an exact sequence $E_2^{2,0} \to \ker(E^2 \to E_2^{0,2}) \to E_2^{1,1} \to E_2^{3,0}$
	\begin{equation} \label{eq:7 term exact}
		\Br(S) \to \ker\big(\Br(\sX) \to \H^0(S,\R^2\rho_*{\bbG_m})\big) \to \H^1(S,\R^1\rho_*{\bbG_m}) \to \H^3(S,{\bbG_m}).
	\end{equation}
	Note that the left morphism is split because there is a section and $\R^2\rho_*{\bbG_m} = 0$ by Artin's theorem~\cite[Cor.\,3.2]{Grothendieck_1968-III} since $S_x$ is the spectrum of a discrete valuation ring and $\rho$ is flat and projective and $S$ Japanese. Let $p\colon S \to \sX$ be a section to $\rho\colon \sX \to S$, i.e., $\rho \circ p = \id_{S}$. It induces a commutative square
	\[
	\begin{tikzcd}
		\H^1(S,\R^1\rho_*{\bbG_m}) \ar[r] \ar[d] & \H^3(S,{\bbG_m}) \ar[d]\\
		\H^1(S,\R^1(\id_S)_*{\bbG_m}) \ar[r] & \H^3(S,{\bbG_m})
	\end{tikzcd}
	\]
	with the lower left entry equal to $0$. Hence the morphism $\H^1(S,\R^1\rho_*{\bbG_m}) \to \H^3(S,{\bbG_m})$ is $0$, so one gets the second short exact sequence of the claim from~\eqref{eq:7 term exact}. The first one follows from an analogous argument for $X\rightarrow \Spec(K)$.
\end{proof}

\begin{lemma}\label[lemma]{model}
For a smooth projective geometrically connected curve $X$ over a finitely generated field $K/k$, 
there exists a smooth integral variety $S/k$ with function field $K$ a projective flat morphism $\rho\colon\sX^\prime \rightarrow S$ with $\sX^\prime$ integral and the generic fibre identified with $X/K$ such that the cokernels of the natural injections
\begin{align*}
{\Br^{\nr}}(K/k) \rightarrow \Br(S),\\
{\Br^{\nr}}(K(X)/k) \rightarrow \Br(\sX^\prime_{\reg})
\end{align*}
are of finite exponent. 
\end{lemma}

\begin{proof}
We write $\sX$ for $\sX^\prime.$ Applying~\cref{lem: S for K}, we may obtain a rational map $\sX \dashrightarrow S$ between smooth integral $k$-varieties, 
such that $K$ is the function field of $S$ and the inclusions
\begin{align*}
{\Br^{\nr}}( K(X)/k) \subseteq \Br(\sX)\\
{\Br^{\nr}}(K/k) \subseteq \Br(S)
\end{align*}
have cokernels of finite exponent. Namely, \cref{lem: S for K} gives the existence of models $\sX$ and $S$, and the existence of the rational map follows by extending the structure morphism $X \to \Spec(K)$.
Let $\overline{S}$ be a compactification of $S$.
By the valuative criterion for properness, the induced rational map $\sX\dashrightarrow \overline{S}$ extends to all points of codimension~$1$, 
so we may assume that we have an honest morphism $\sX\rightarrow \overline{S}$.
Note that codimension $>1$ points do not contribute to the unramified Brauer group.
Take a compactification $\overline{\sX}\rightarrow \overline{S}$ of $\sX$ over $\overline{S}$, which exists by Nagata's theorem. 
Applying Raynaud--Gruson's flattening theorem \cite{RaynaudGruson_1971} and taking normalisation, 
we assume without loss of generality that  
$\overline{\sX}\rightarrow \overline{S}$ is flat. Replacing $\overline{\sX}$ by its normalisation, we can assume that  all fibres of $\overline{\sX}\rightarrow \overline{S}$ has dimension $\leq 1$.
Then the generic fibre of  $\overline{\sX} \to \overline{S}$ 
is a normal curve over $K$ with function field $ K(X)/K$, 
hence it equals $X/K$. We write now $\overline{\sX}\rightarrow \overline{S}$ as $\rho\colon \sX\rightarrow S$.
This has generic fibre $X/K$ and by~\cref{lem: S for K}, the inclusions
\begin{align*}
{\Br^{\nr}}(  K(X)/k) \subseteq \Br(\sX_{\reg})\\
{\Br^{\nr}}(K/k)\subseteq \Br(S_{\reg})
\end{align*}
have a cokernel of finite exponent. 
Moreover, every fibre of $\rho$ has dimension $\leqslant 1$, thus for the codimension-1 points we have 
$\rho(\sX^1) \subset S^1\cup\{\eta_S\}$.
So,
$$
\bigcap_{s\in S^1}\Br((\sX_{\sO_{S,s}})_\reg) = \Br(\sX_{\reg}),
$$
where $\sX_{\sO_{S,s}}$ denotes the base change of $\rho$ to the local ring at $s$. Thus, by removing a closed subset with a codimension at least 2 from $S$, we get a projective-flat morphism (this does not affect Brauer groups since $\rho(\sX^1) \subset S^1\cup\{\eta_S\}$).
\end{proof}

We introduce and compare two kinds of unramified Brauer groups with respect to different valuations.

\begin{lemma} \label[lemma]{exact sequences with Brauer and Sha}
	Let $X/K$ be a smooth projective geometrically integral curve and $S$ and $\sX^\prime$ as in~\cref{model}. Assume that $X(K) \neq \emptyset$. For any $x\in S^1$, let $\sX_{S_x}$ be a desingularisation of $\sX^\prime_{S_x}$. Let $\rho\colon \sX_{S_x}\rightarrow S_x\colonequals\Spec(\cO_{S,x})$ be the induced morphism. Then from~\cref{claim 4}, there is a commutative diagram with exact rows
	\[
	\begin{tikzcd}
		0 \ar[r] & \Br(K) \ar[r] \ar[d] & \Br(X) \ar[r] \ar[d]& \H^1(K,\Pic(X_{K^\sep})) \ar[d,"\varphi"]\ar[r] & 0\\
		0 \ar[r] & \displaystyle\bigoplus_{x \in S^1}\frac{\Br(K)}{\Br(\sO_{S,x})} \ar[r]& \displaystyle\bigoplus_{x \in S^1}\frac{\Br(X)}{\Br(\sX_{S_x})} \ar[r]& \displaystyle\bigoplus_{x \in S^1}\frac{\H^1(K,\Pic(X_{K^\sep}))}{\H^1(S_x,\R^1\rho_*{\bbG_m})} \ar[r] & 0.
	\end{tikzcd}
	\]	
	The diagram induces short exact sequences (the last arrow has a cokernel of finite exponent)
	\[
	0 \to \Br^\nr(K/k) \to \Br^\nr(X/k) \to \ker(\varphi) \dashrightarrow 0
	\]
	and an exact sequence
	\[
	0 \to \ker(\varphi) \to \Sha_S(\Pic_{X/K}) \to \bigoplus_{x\in S^1}\frac{\ker(\H^1(K,\Pic_{X/K}) \to \H^1(K_x^\sh,\Pic_{X/K}))}{\H^1(S_x,\R^1\rho_*{\bbG_m})},
	\]
	where
	\[
	\Sha_S(\Pic_{X/K}) \colonequals \ker\Big(\H^1(K,\Pic_{X/K}) \to \bigoplus_{x \in S^1}\H^1(K_x^\sh,\Pic_{X/K})\Big),
	\]
 and the last group in the second exact sequence is of finite exponent.
\end{lemma}
\begin{proof}
	By~\cref{claim 4},
	\[
	\begin{tikzcd}
		0 \ar[r] & \Br(\sO_{S,x}) \ar[d]\ar[r] & \Br(\sX_{S_x}) \ar[d]\ar[r] & \H^1(S_x,\R^1\rho_*{\bbG_m}) \ar[d]\ar[r] & 0\\
		0 \ar[r] & \Br(K) \ar[r] & \Br(X) \ar[r] & \H^1(K,\Pic_{X/K}) \ar[r] & 0
	\end{tikzcd}
	\]
	is a commutative diagram with exact rows for $x \in S^1$. By~\cref{claim 2}, $p \in X(K)$ induces $\Br(X) \to \Br(K)$ such that the two rows split compatibly. The snake lemma gives a short exact sequence
 \[
	0 \to \Br(S) \to \bigcap_{s\in S^1}\Br(\sX_{S_x}) \to \ker(\varphi) \to 0
	\]
By~\cref{model}, the inclusions
$$\Br^\nr(K/k)\rightarrow\Br(S),$$
$${\Br^\nr}(K(X)/k) \rightarrow \bigcap_{s\in S^1}\Br(\sX_{S_x})$$
have cokernels of finite exponent since $$\bigcap_{s\in S^1}\Br(\sX_{S_x})\subseteq \bigcap_{s\in S^1}\Br((\sX^\prime_{S_x})_{\reg})=\Br(\sX^\prime_{\reg}).$$
This implies that the induced natural map
$${\Br^\nr}(K(X)/k)\rightarrow \ker(\varphi)$$
has a cokernel of finite exponent and has a kernel which is equal to $\ker(\Br^{\nr}(X/k))\rightarrow \H^1(K,\Pic_{X/K}))$. Since $\Br(O_{K,v}) \rightarrow\Br(\sX_{O_{K,v}})$ is split for any $v\in div\Val(K)$, $\ker(\Br^{\nr}(X/k))\rightarrow \H^1(K,\Pic_{X/K}))$ is equal to
$$\Br(K)\cap\Br^{\nr}(X/k)=\Br(K)\cap \bigcap_{v\in div\Val(K)}\Br(\sX_{O_{K,v}})=\bigcap_{v\in div\Val(K)}\Br(O_{K,v})=\Br^\nr(K/k).$$
The second short exact sequence follows from~\cref{claim 3}. Since for all but finite many $x\in S^1$, $\sX^\prime$ is smooth and projective over $S_x$, \cref{claim 3} implies the last claim.
\end{proof}

\begin{corollary} \label[corollary]{main theorem with a point}
Let the conditions be the same as in the previous lemma. Then there is an exact sequence
	\[
	0 \to {\Br^\nr}(K/k) \to {\Br^\nr}(K(X)/k) \to\Sha_S(\Pic^0_{X/K})
	\]
in which the last arrow has a cokernel of finite exponent.
\end{corollary}
\begin{proof}
	By~\cref{exact sequences with Brauer and Sha}, it suffices to show that $\Sha_S(\Pic^0_{X/K}) \xrightarrow{\sim} \Sha_S(\Pic_{X/K})$. Since $X(K) \neq \emptyset$, $\Pic^0(X_{K^\sep}) \oplus \bbZ \cong \Pic(X_{K^\sep})$ as $G_K$-modules. The claim follows from $\H^1(K,\bbZ) = 0$.
\end{proof}

We can now obtain the following key result for relative curves: 

\begin{proposition} \label[proposition]{ShaX cokernel of finite exponent}
Let $K/k$ be a finitely generated field and $X/K$ a smooth projective geometrically connected curve.
There exists a smooth integral variety $S/k$ 
with function field $K$ such that the injection
$$
\Sha'(\Pic^0_{X/ K}/K)\rightarrow \Sha_S(\Pic^0_{X/  K})
$$
has a cokernel of finite exponent.
\end{proposition}

\begin{proof}
First, we show that if the claim is true for $X_L/L$ for some finite Galois extension $L/K$, then the claim is true for $X/K$. Let $T/k$ be a smooth integral variety over $k$ with function field $ K(T)=L$ such that 
$$
\Sha'(\Pic^0_{X_L/L}/L)\rightarrow \Sha_T(\Pic^0_{X_L/ L})
$$
has a cokernel of finite exponent.
Let $S_0$ be projective normal integral variety over $k$ with function field $K$. Then the $k$-rational map $T\dashrightarrow S_0$ can be defined at all $t\in T^1$. Thus, by shrinking $T$, we may assume that $T\dashrightarrow S_0$ is a morphism. By Raynaud--Gruson’s flattening theorem, there are birational morphisms $S_1 \rightarrow S_0$ and $T^\prime\rightarrow T$, and a flat morphism $T^\prime\rightarrow S_1$ with $S_1$ normal. Since $\T^1\subseteq T^\prime_{\mathrm{reg}}$, we may replace $T$ by $T^\prime_{\mathrm{reg}}$. Let $S$ denote the regular locus of $S_1$. So $\T^1$ is mapped into $S^1$. We will show that $S$ satisfies the claim. 
The Galois group $G \colonequals \Gal(L/K)$ acts on the set $div\Val(L)$ of divisorial valuations on $L/k$ and is compatible with the restriction $div\Val(L) \rightarrow div\Val(K)$.
Consider the commutative diagram
$$
\xymatrix{
\Sha'(\Pic^0_{X/K}/K)\ar[r]\ar[d] &\Sha_S(\Pic^0_{X/K}) \ar[d]
\\
\Sha'(\Pic^0_{X_L/L}/L)^{G}\ar[r]&\Sha_T(\Pic^0_{X_L/L}) 
}
$$
and the natural map $H^1(K,\Pic^0_{X/K}) \rightarrow H^1(L,\Pic^0_{X_L/L})^{G}$.
Let $a\in \Sha_S(\Pic^0_{X/K})$ and $b$ is its image $H^1(L,\Pic^0_{X_L/L})^{G}$. By the assumption, the cokernel of the second row is killed by some positive integer $N$. So $Nb \in \Sha'(\Pic^0_{X_L/L}/L) \cap H^1(L,\Pic^0_{X_L/L})^{G}$. 
Since the first column has a kernel and a cokernel of finite exponent. For this, note that for the Galois cohomology of $L/K$ and $L \otimes_K K_v/K_v$, we have $\res$ and $\cores$, and $\res \circ \cores$ is the norm on $\H^1(L)$, which restricts to multiplication by $[L:K]$ on $\H^1(L)^G$. The other composition $\cores \circ \res$ is also multiplication by $[L:K]$. Now consider
\[
\begin{tikzcd}
	0 \ar[r] & \Sha'(\Pic^0_{X/K}/K) \ar[d]\ar[r] & \H^1(K,\Pic^0_{X/K}) \ar[d]\ar[r] & \prod_v\H^1(K_v,\Pic^0_{X/K}) \ar[d]\\
	0 \ar[r] & \Sha'(\Pic^0_{X_L/L}/L)^G \ar[r] & \H^1(L,\Pic^0_{X_L/L})^G \ar[r] & \prod_v\H^1(L_v,\Pic^0_{X_L/L})^G
\end{tikzcd}
\]
By the above, the kernel and cokernel of the two last vertical arrows are killed by $[L:K]$, hence also the kernel on the Tate--Shafarevich groups. For the cokernel of the Tate--Shafarevich groups, take $a \in \Sha'(\Pic^0_{X_L/L}/L)^G$. Then $[L:K]a \in \H^1(L,\Pic^0_{X_L/L})$ is the restriction of $b \in \H^1(K,\Pic^0_{X/K})$. One has $[L:K]b_v$ goes to $0$ on the bottom right, hence $[L:K]^2b$ goes to $0$ on the top right, hence comes from $\Sha'(\Pic^0_{X/K}/K)$. By diagram chasing, the first row has a cokernel of finite exponent.

Second, by the argument above, we can assume without loss of generality that $X(K)$ is not empty.
The Leray spectral sequence $\H^p(K,\R^q\rho_*\bbG_m) \Rightarrow \H^{p+q}(X,\bbG_m)$ for $\rho\colon X \to \Spec(K)$ gives a short exact sequence
$$
0 \rightarrow \Br(K) \rightarrow \Br(X) \rightarrow \H^1(K, \Pic_{X/K}) \rightarrow 0
$$
Thus by \cref{lem: Br to Sha} and the fact $\Sha'(\Pic^0_{X/K}/K)=\Sha'(\Pic_{X/K}/K)$, we also have an exact sequence
$$
0 \rightarrow \Br^{\nr}(K/k) \rightarrow \Br^{\nr}( K(X)/k) \rightarrow \Sha'(\Pic^0_{X/K}/K).
$$
By \cref{main theorem with a point}, there exists a smooth integral variety $S$ with function field $K$ and an exact sequence
$$
0 \rightarrow \Br^{\nr}(K/k) \rightarrow \Br^{\nr}( K(X)/k) \rightarrow \Sha_S(\Pic^0_{X/K})
$$
such that the last morphism has a cokernel of finite exponent. Since the last map factor through $\Sha'(\Pic^0_{X/K}/K)\rightarrow \Sha_S(\Pic^0_{X/K})$, this finishes the proof.
\end{proof}
Finally we obtain the general version of the above result: 

\begin{corollary}\label[corollary]{cor: sha}
Let $K/k$ be a finitely generated infinite field and $A/K$ an abelian variety.
There exists a smooth irreducible (hence integral) variety $S/k$ 
with function field $K$ such that the injection
$$
\Sha'(A/K)\rightarrow \Sha_S(A)
$$
has a cokernel of finite exponent.
\end{corollary}

\begin{proof}
Since $K$ is infinite, by \cite[Thm.\,10.1]{Milne_1986_JV} there exists a curve $X/K$ and a surjection $\Pic^0_{X/K} \to A$. 
By Poincar\'e's complete reducibility theorem \cite[Prop.\,12.1]{Milne_1986_AV}, 
it is split up to isogeny, i.e., there exists an abelian variety $B/K$ such 
that there is an isogeny $\Pic^0_{X/K} \xrightarrow{sim} A \times B$. 
But both $\Sha'(-/K)$ and $\Sha_S(-)$ commute with finite products and 
having cokernel of finite exponent is preserved by isogenies.
It therefore suffices to prove the statement for $\Pic^0_{X/K}$. 

As $X$ is a curve, we have a short exact sequence of the form
$$0 \to \Pic^0_{X/K} \to \Pic_{X/K} \to \bbZ \to 0.$$
In étale cohomology it induces an exact sequence of the form
$$
0 \rightarrow  \bbZ/\deg(\Pic_{X/K}(K)) \rightarrow  \H^1(K, \Pic^0_{X/K}) \rightarrow \H^1(K, \Pic_{X/K}) \rightarrow \H^1(K,\bbZ).
$$
Note that $\H^1(K,\bbZ)$ is finite. 
Indeed, after a finite field extension $L/K$ the Galois action on $\bbZ$ is trivial, and hence $\H^1(L,\bbZ)\cong \H^1(G_L,\bbZ)^{G_K/G_L}\cong\Hom_{cont}(G_L,\bbZ)=0$. 
Consider the inflation-restriction exact sequence
$$
0 \rightarrow \H^1(G_K/G_L,\bbA^{G_L}) \rightarrow \H^1(G_K,\bbZ) \rightarrow \H^1(G_L,\bbZ)^{G_K/G_L}
$$
arising as the five-term exact sequence of the Hochschild--Serre spectral sequence $\H^i(G_K/G_L, \H^j(G_L,\bbZ)) \Rightarrow \H^{i+j}(G_K,\bbZ) $.
By what we said above, the last term of the exact sequence is trivial. 
Moreover, the first term $\H^1(G_K/G_L,\bbA^{G_L})=\H^1(\Gal(L/K),\bbA^{G_L})$ is finite, as $L/K$ is finite cyclic. 
Consequently, the middle term of the exact sequence is finite as desired. 
This means that $ \H^1(K, \Pic^0_{X/K})$ and $ \H^1(K, \Pic_{X/K})$ are isomorphic modulo the Serre subcategory of finite abelian groups.  

We can make the same argument locally at each divisorial valuation of $K$, 
meaning that for each $s\in div\Val(K)$, we have an exact sequence
$$
0 \rightarrow  \bbZ/\deg(\Pic_{X/K}(K^{\sh}_s)) \rightarrow  \H^1(K^\sh_s, \Pic^0_{X/K}) \rightarrow \H^1(K^\sh_s, \Pic_{X/K}) \rightarrow \H^1(K^\sh_s,\bbZ)
$$
where the first and the last group are finite. 
We observe that the order of $\bbZ/\deg(\Pic_{X/K}(K^{\sh}_s))$ is bounded by $|\bbZ/\deg(\Pic_{X/K}(K))|$ and  that of $\H^1(K^\sh_s,\bbZ)$ is bounded by $|\H^1(K,\bbZ)|$. 
It follows that the natural morphisms of products
\begin{align*}
\prod_{s\in div\Val(K)} \H^1(K^\sh_s, \Pic^0_{X/K}) \rightarrow \prod_{s\in div\Val(K)}\H^1(K^\sh_s, \Pic_{X/K}),\\
\prod_{s\in S^1} \H^1(K^\sh_s, \Pic^0_{X/K}) \rightarrow \prod_{s\in S^1}\H^1(K^\sh_s, \Pic_{X/K})
\end{align*}
are isomorphisms modulo the Serre subcategory of finite abelian groups, too. 

By \cref{def: Shafarevich-Tate for abelian varieties}, it follows that 
the same is true for the natural morphisms
\begin{align*}
\Sha'(\Pic^0_{X/K}/K) \rightarrow \Sha'(\Pic_{X/K}/K),\\
\Sha_S(\Pic^0_{X/K}) \rightarrow \Sha_S(\Pic_{X/K}). 
\end{align*}
As a consequence, the fact that the injection $\Sha'(\Pic_{X/K}/K)\hookrightarrow \Sha_S(\Pic_{X/K})$ has cokernel of finite exponent by \cref{ShaX cokernel of finite exponent}
implies that the same holds for the injection $\Sha'(\Pic^0_{X/K}/K)\hookrightarrow \Sha_S(\Pic^0_{X/K})$. 
\end{proof}

We refine the above arguments now to obtain the main theorem of this section.

\begin{lemma}\label[lemma]{injectivelem}
Let $X/K^\prime$ be a smooth projective geometrically connected variety over a field $K'$. Let $R$ be a strictly Henselian DVR with fraction field $K$. Assume that $K/K^\prime$ is a separable algebraic field extension and $\rho\colon \sX\rightarrow S=\Spec(R)$ a proper flat morphism with the generic fibre identified with $X_K/K$.
Then the kernel of 
$$
\Br(\sX)\rightarrow \Br(X_{K^{\sep}})
$$
is of finite exponent which only depends on $X/K^\prime$.
\end{lemma}
\begin{proof}
If $\dim(\sX)\leqslant 2$, Artin's theorem implies that both groups vanish. In general, we will use a pull-back trick developed by Colliot-Thél\`ene and Skorobogatov to reduce it to cases of relative dimension $1$. Let $L^\prime/K^\prime$ be a finite Galois extension such that $\Pic(X_{L^\prime})\rightarrow \NS(X_{K^\sep})$ is surjective. 
By \cite[\S\,3.3]{Yua2}, there are finitely many smooth projective geometrically connected curves $Z_i/L^\prime \subset X_{L^\prime}/L^\prime$, an abelian variety $A/L^\prime$ and a $G_{L^\prime}$-equivariant map
$$\Pic(X_{K^\sep})\times A((L^\prime)^\sep)\rightarrow\bigoplus_i\Pic(Z_{i,(L^\prime)^\sep})$$
with finite kernel and cokernel. By extending $L^\prime$, we may assume that $Z_i(L^\prime)$ is not empty for all $i$. Let $L$ denote $KL^\prime$. For simplicity, we still use $Z_i$ to denote their base change to $L$ over $L^\prime$. Think $Z_i$ as a scheme over $K$. It follows that the natural map
$$\H^1(K,\Pic_{X_K/K})\rightarrow \bigoplus_i \H^1(K,\Pic_{Z_i/K})$$
has a kernel killed by a positive integer that only depends on $X/K^\prime$ since the kernel is mapped into the kernel of 
$$\H^1(L,\Pic_{X_K/K})\rightarrow \bigoplus_i \H^1(L,\Pic_{Z_i/L})$$ 
under the natural map $\H^1(K,\Pic_{X_K/K}) \rightarrow \H^1(L,\Pic_{X_K/K})$ which has a kernel killed by $[L^\prime:K^\prime]$. By taking the Zariski closures of $Z_i$ in $\sX_{R_L}$ and then desingularizing it, we may assume that there are $\pi_i\colon\sZ_i\longrightarrow S$ with generic fibre $Z_i/K$ and $S$-morphisms $\sZ_i\longrightarrow \sX$.
Consider the  commutative diagram \\
\begin{displaymath}
\xymatrix{ 
&\Br(K)\ar[r]\ar[d] & \ker(\Br(X_K)\rightarrow  \Br(X_{K^\sep})^{G_K}) \ar[d]\ar[r] &\H^1(K,\Pic_{X_K/K}) \ar[d]
\\
0\ar[r]&\bigoplus_i\Br(L)\ar[r] & \bigoplus_i\ker( \Br(Z_i)\rightarrow \Br(Z_{i,K^\sep})^{G_K}) \ar[r]&\bigoplus_i \H^1(K,\Pic_{Z_i/K}) }
\end{displaymath}
The exact rows come from the spectral sequences\\
$$
\H^p(K,\H^q(X_{K^\sep},\mathbb{G}_m))\Rightarrow \H^{p+q}(X_K,\mathbb{G}_m)
$$
and its analogue for $Z_i/K$. Since $\Br(\sZ_i)=0$, the intersection $\Br(\sX)\cap \ker(\Br(X_K)\rightarrow  \Br(X_{K^\sep})^{G_K})$ will be mapped to zero in $\Br(Z_i)$. Thus, a diagram chase tells us that it is contained in the kernel of the second column which is killed by a positive number that only depends on $X/K^\prime$.
\end{proof}
\begin{lemma}\label[lemma]{keylemma}
Let $K/k$ be a finitely generated field and $X/K$ a smooth projective geometrically connected variety.
There exists a smooth integral variety $S$ with function field $K$ together with a projective $k$-morphism $\rho\colon\sX\rightarrow S$ of integral varieties over $k$ which extends the $k$-morphism $\Spec(K(X))\rightarrow \Spec(K)$ such that the natural injections
\begin{align*}
\Br_{\nr}(  K/k) &\rightarrow \Br(S)\\
\Br_{\nr}(  K(X)/k) & \rightarrow \bigcap_{s\in S^1} \Br((\sX_{\sO_{S,s}})_{\reg})\\
\Sha'(\Pic^0_{X/  K,\red}/K) &\rightarrow \Sha_S(\Pic^0_{X/  K,\red})
\end{align*}
have cokernels of finite exponent.
\end{lemma}
\begin{proof}
By  \cref{lem: S for K} and \cref{cor: sha}, there exists a $k$-rational map $\sX \dashrightarrow S$ between smooth integral $k$-varieties which extends the $k$-morphism $\Spec(K(X))\rightarrow \Spec(K)$  such that the inclusions
\begin{align*}
\Br_{\nr}(  K(X)/k)  & \rightarrow \Br(\sX)\\
\Br_{\nr}(K/k)  & \rightarrow \Br(S)\\
\Sha'(\Pic^0_{X/  K,\red}/K)  & \rightarrow \Sha_S(\Pic^0_{X/  K,\red})
\end{align*}
have cokernels of finite exponent. 
Let $\overline{S}$ be a projective compactification of $S$.
The induced rational map $\sX\dashrightarrow \overline{S}$ extends to all points of codimension 1, so we may assume that we have an honest morphism $\sX\rightarrow \overline{S}$. Take a compactification $\overline{\sX}\rightarrow \overline{S}$ of $\sX$ over $\overline{S}$, which exists by Nagata's theorem. 
Applying Raynaud--Gruson's flattening theorem \cite{RaynaudGruson_1971} and taking normalisation, we get $\overline{\sX}_1\rightarrow \overline{S}_1$ such that $\overline{\sX}_1$ and $\overline{S_1}$ are normal, and $(\overline{\sX}_1)^1$ is mapped into  $(\overline{S}_1)^1\cup\{\eta_S\}$. Since the morphisms
 $\overline{\sX}_1\rightarrow \overline{\sX}$ and $\overline{S}_1\rightarrow\overline{S}$ are birational morphisms, any codimension-1 point of $\overline{\sX}$ (resp. $\overline{S}$) has at least one preimage of codimension 1 in $\overline{\sX}_1$ (resp. $\overline{S}_1$). From now, denote $\overline{\sX}_1$ (resp. $\overline{S}_1$) by $\overline{\sX}$ (resp. $\overline{S}$). So the following inclusion maps
\begin{align*}
\Br_{\nr}(  K(X)/k)  & \rightarrow \Br(\overline{\sX}_\reg)\\
\Br_{\nr}(K/k)  & \rightarrow \Br(\overline{S}_\reg)\\
\Sha'(\Pic^0_{X/  K,\red}/K)  & \rightarrow \Sha_{\overline{S}_\reg}(\Pic^0_{X/  K,\red})
\end{align*}
have cokernels of finite exponents.
Since all codimension-$1$ points of $\overline{\sX}$ are contained in $\sX$ and are mapped into $S^1\cup\{\eta_S\}$ by the map $\sX\rightarrow S$, we have
$$
\bigcap_{s\in S^1}\Br((\sX_{\sO_{S,s}})_\reg) = \Br(\sX_{\reg})=\Br(\overline{\sX}_{\reg}).
$$
Thus, the map $\sX\rightarrow S$ satisfies the conditions in the lemma.
\end{proof}

\begin{lemma} \label[lemma]{prop1}
Let $R$ be an excellent DVR, set $S=\Spec(R)$ and $K=\mathrm{Frac}(R)$. 
Let $\rho\colon \sX\rightarrow S$, with $\sX$ regular, be a projective flat morphism with smooth projective geometrically connected generic fibre $X\rightarrow \Spec(K)$. Assume that $X(K)$ is not empty.
Then there exists an exact
sequence 
$$
0 \rightarrow \Br(S) \rightarrow \ker\big(\Br(\sX)\rightarrow \Br(X_{K^{\sep}})\big) \rightarrow \Sha_S(\Pic_{X/K}),
$$
where the last morphism has a cokernel of finite exponent.
By abuse of notation (which is justified because $R$ is a local ring), 
we write here
$$\Sha_S(\Pic_{X/K}) \colonequals \ker\big(\H^1(K,\Pic_{X/K})\rightarrow \H^1(K^{\sh},\Pic_{X/K})\big).$$
\end{lemma}
\begin{proof}
Consider the Leray spectral sequence for $\rho\colon \sX \to S$
	\[
		E_2^{pq} = \H^p(S,\R^q\rho_*{\bbG_m}) \Longrightarrow \H^{p+q}(\sX,{\bbG_m}).
	\]
	It gives rise to an exact sequence $E_2^{2,0} \to \ker(E^2 \to E_2^{0,2}) \to E_2^{1,1} \to E_2^{3,0}$
	\begin{equation}
		\Br(S) \to \ker\big(\Br(\sX) \to \H^0(S,\R^2\rho_*{\bbG_m})\big) \to \H^1(S,\R^1\rho_*{\bbG_m}) \to \H^3(S,{\bbG_m}).
	\end{equation}
	Note that the left morphism is split because there is a section. Let $p\colon S \to \sX$ be a section to $\rho\colon \sX \to S$, i.e., $\rho \circ p = \id_{S}$. It induces a commutative square
	\[
	\begin{tikzcd}
		\H^1(S,\R^1\rho_*{\bbG_m}) \ar[r] \ar[d] & \H^3(S,{\bbG_m}) \ar[d]\\
		\H^1(S,\R^1(\id_{S})_*{\bbG_m}) \ar[r] & \H^3(S,{\bbG_m})
	\end{tikzcd}
	\]
	with the lower left entry equal to $0$. Hence the morphism $\H^1(S,\R^1\rho_*{\bbG_m}) \to \H^3(S,{\bbG_m})$ is $0$, so one gets a short exact sequence
 \begin{equation}
		0\to\Br(S) \to \ker\big(\Br(\sX) \to \H^0(S,\R^2\rho_*{\bbG_m})\big) \to \H^1(S,\R^1\rho_*{\bbG_m}) \to 0.
	\end{equation}
 Since $\H^0(S,\R^2\rho_*{\bbG_m})\subseteq \H^0(S^{\sh},\R^2\rho_*{\bbG_m})=\Br(\sX_{S^{\sh}})$ and $\Br(\sX_{S^\sh})\rightarrow \Br(X_{K^\sep})$ has a kernel of finite exponent by \cref{injectivelem}, the inclusion
 $$\ker\big(\Br(\sX) \rightarrow \H^0(S,\R^2\rho_*{\bbG_m})\big) \rightarrow \ker\big(\Br(\sX) \rightarrow\Br(X_{K^\sep})\big)  $$
 has a cokernel of finite exponent by~\cref{claim 3}\,(ii), the natural map 
 $$
 \H^1(S,\R^1\rho_*{\bbG_m}) \rightarrow \Sha_S(\Pic_{X/K})
 $$
 has a kernel and a cokernel of finite exponent and is an isomorphism if $\rho$ is smooth. This proves the lemma.
 \end{proof}
 
\begin{remark}
Let $K$ be the function field of a smooth integral variety $S$ over $\mathbb{F}_p$. 
Then for all but finitely many $s\in S^1$, $X/K$ extends to a smooth projective morphism $\sX \rightarrow \Spec(R)$, where $R$ is the local ring at $s$. 
If one now varies the lemma above over such $R$, by \cref{injectivelem}, the exponents involved in the statement are uniformly bounded by a positive number independent of $s\in S^1$.
\end{remark}

\begin{lemma}\label[lemma]{H1Pic almost injective}
Let $L/K$ be a finite extension over a field of characteristic $p > 0$, $Y/L$ and $X/K$ smooth geometrically connected projective varieties such that $Y\rightarrow X_L$ is a generically finite $L$-morphism. 
Then the natural map 
$$
\H^1(K,\Pic_{X/K}) \rightarrow \H^1(L,\Pic_{Y/L})
$$
has a kernel of finite exponent. 
\end{lemma}
\begin{proof}
The morphism $\H^1(K,\Pic_{X/K}) \rightarrow \H^1(L,\Pic_{Y/L})$ factors as $\beta \circ \alpha$ where $\alpha\colon \H^1(K,\Pic_{X/K}) \rightarrow \H^1(L,\Pic_{X/K})$ and $\beta\colon \H^1(L,\Pic_{X/K}) \to \H^1(L,\Pic_{Y/L})$.
Thus the claim follows if we prove that both  $\alpha$ and $\beta$ have kernel of finite exponent.
It suffices to prove the statements with $\Pic$ replaced by $\Pic^0$.
Indeed, their quotient, the Néron--Severi group, becomes the trivial Galois module over a finite extension, over which its first cohomology group is  
of finite exponent. 
We first show the statement for $\beta$.
As $Y \to X_L$ is generically finite and proper, hence in particular surjective, it follows from a restriction-corestriction argument that $\H^1(X_{L^{\sep}},\bbQ_\ell) \to \H^1(Y_{L^{\sep}},\bbQ_\ell)$ is injective for all $\ell \neq p$.
But these are the duals of the rationalised $\ell$-adic Tate modules, 
and thus $\Pic^0_{X/L} \to \Pic^0_{Y/L}$ has finite kernel. 
Consequently $\beta$ has kernel of finite exponent.
 
Next we show the statement for $\alpha$.
If $L/K$ is separable, it follows again from a restriction-corestriction argument that the kernel of $\alpha$ is of finite exponent. 
Thus let us assume without loss of generality that $L/K$ is purely inseparable of degree $p$, hence so is $L^{\sep}/K^{\sep}$. Consider the morphisms 
\begin{align*}
[p]\colon \Pic^0_{X/K} \rightarrow \Pic^0_{X/K},&&
F\colon \Pic^0_{X/K} \rightarrow (\Pic^0_{X/K})^\sigma,&&
 V\colon (\Pic^0_{X/K})^\sigma \rightarrow  \Pic^0_{X/K},
\end{align*} 
where the first one is multiplication by $p$, and the second and third satisfy the equation $[p] = V \circ F$. Here, $-^\sigma$ denotes the Frobenius twist.
Since $F(L^{\sep}) \subseteq K^{\sep}$, $F(\Pic^0_{X/K}(L^{\sep})) \subseteq (\Pic^0_{X/K})^{\sigma}(K^{\sep})$, one has $[p](\Pic^0_{X/K}(L^{\sep})) \subseteq V((\Pic^0_{X/K})^\sigma(K^{\sep})) = \Pic^0_{X/K}(K^{\sep})$. Hence $\Pic^0_{X/K}(L^{\sep})/\Pic^0_{X/K}(K^{\sep})$ has exponent dividing $p$. Since $\Gal(K^{\sep}/K) = \Gal(L^{\sep}/L)$, it follows that $\alpha\colon \H^1(K,\Pic^0_{X/K}(K^{\sep})) \to \H^1(L,\Pic^0_{X/K}(L^{\sep}))$ has kernel of finite exponent.
\end{proof}

\begin{lemma}
Let $K/k$ be a finitely generated field and  $X/K$ be smooth projective geometrically connected variety. 
Assume that $X(K)$ is not empty. For a divisorial valuation $s$ on $K$, 
let $R=\sO_{K,s}$ be the discrete valuation ring at $s$ and set $S=\Spec(R)$. 
There exists an exact sequence 
$$
0 \rightarrow \Br(R) \rightarrow \ker\big(\Br_{\nr}(  K(X)/R) \rightarrow \Br(X_{K^{\sep}})\big) \rightarrow \Sha_S(\Pic_{X/K})
$$
such that the last map has a cokernel of finite exponent
which only depends on $X/K$ and is independent of $s$.
Here by abuse of notation we write $\Br_{\nr}( K(X)/R) \colonequals \cap_{\sX}\Br(\sX)$, taking the intersection over all separated regular models (not necessarily proper) of finite type over $R$.
\end{lemma}
\begin{proof}
By~\cite[Thm.\,8.2]{deJong_1996}, there is a finite extension $L/K$
and a flat projective morphism $\frX \to \Spec(R_L)$ 
with $\frX$ regular and $\frX_L$ smooth over $L$ 
and an $R$-morphism $\frX \to \widetilde{\frX}$ which is an alteration, 
where $\widetilde{\frX}$ is a projective model of $X$ over $R$. 
Replacing $L$ by its algebraic closure, we can assume that $\frX_L$ is geometrically connected. 
Note that $L/K$ may be inseparable. By~\cref{prop1}, the claim is true for $\frX_L \to \Spec(L)$.

Consider the commutative diagram
\[
\begin{tikzcd}
	\ker(\Br_{\nr}(K(X)/R) \to \Br(X_{K^{\sep}})) \ar[d,"b"]\ar[r,"a"] & \H^1(K,\Pic_{X/K}) \ar[d,"d"]\\
	\ker(\Br(\frX) \to \Br(X_{L^{\sep}})) \ar[r,"c"] & \H^1(K,\Pic_{X/K}).
\end{tikzcd}\]
We have $\ker(a) \subseteq \ker(c) \cap \Br(K) = \Br(R_L) \cap \Br(K) {\leftarrow} \Br(R)$ where the last morphism is an isomorphism up to finite exponent by a restriction-corestriction argument. 
By~\cref{H1Pic almost injective}, for any $w \mid v$, $\H^1(K_v^{\sh},\Pic_{X/K}) \to \H^1(L_w^{\sh},\Pic_{\frX_L/L})$ has a kernel of finite exponent. 
Hence the same is true for  $\im(a) \subseteq \Sha_{S_v}(\Pic_{X/K})$ 
where $\Sha_{S_v}(\Pic_{X/K})$ is defined as in~\cref{prop1}.

On the other hand, for all $y \in \ker(\Br(X) \to \Br(X_{K^{\sep}}))$ such that $x = a(y) \in \Sha_{S_v}(\Pic_{X/K})$, $c(b(y)) = d(x) \in \bigcap_{w \mid v}\Sha_{S_w}(\Pic_{\frX_L/K})$. Hence there exists $n_v > 0$ (independent of $x$) such that $n_vb(y) \in \bigcap_{w \mid v}\ker(\Br(\frX_{R_{L,w}}) \to \Br(\frX_{L^{\sep}})) = \ker(\Br(\frX) \to \Br(\frX_{L^{\sep}}))$. Hence there is $m_v > 0$ independent of $y$ such that $m_vy \in \ker(\Br_{\nr}( K(X)/R) \to \Br(X_{K^{\sep}}))$.

Note that $\Br(\frX) \cap \Br(X)$ and $\Br_{\nr}( K(X)/R)$ can be identified up to finite exponent 
by a restriction-corestriction argument. 
Since $ K(\frX)/ K(X)$ is a finite extension, the map from $\Sha_{S_v}(\Pic_{X/K}) \cap \Br(X)$ to $\ker(\Br(\frX) \to \Br(X_{K^{\sep}}))$ has a kernel of finite exponent.
\end{proof}

\begin{remark}
The Hochschild--Serre spectral sequence \[E_2^{p,q} = \H^p(K,\H^q(X_{K^\sep},{\bbG_m})) \Rightarrow \H^{p+q}(X,{\bbG_m})\] in \'etale cohomology gives an exact sequence
	\[
	\Br(K) \to \ker\big(\Br(X) \to \Br(X_{K^\sep})^{G_K}\big) \to \H^1(K,\Pic(X_{K^\sep})) \to \H^3(K,{\bbG_m}).
	\]
Let $L/K$ be a finite Galois extension such that $X(L)$ is not empty. Then there is an exact sequence
\[
	0\to\Br(L) \to \ker\big(\Br(X_L) \to \Br(X_{K^\sep})^{G_L}\big) \to \H^1(L,\Pic(X_{K^\sep})) \to 0.
	\]
Assuming that we have proved 
$$
0\longrightarrow \Br_{\nr}(L/k)\longrightarrow 
\ker\big(\Br_{\nr}(  K(X_L)/k)\rightarrow \Br(X_{L^{\sep}})\big) 
\longrightarrow \Sha'(\Pic_{X_L/L})\dashrightarrow 0,
$$
it is easy to deduce the theorem below for $X/K$ by a diagram chase since the natural maps
\begin{align*}
\Br_{\nr}( K/k)  & \rightarrow \Br_{\nr}(L/k)^{G}\\
\Br_{\nr}(  K(X)/k)  & \rightarrow \Br_{\nr}(  K(X_L)/k)^G\\
\Sha'(\Pic_{X/K}/K)  & \rightarrow \Sha'(\Pic_{X/L}/L)^G
\end{align*}
have kernels and cokernels of finite exponent, where $G$ is the Galois group of $L/K$. Thus, we can always assume that $X(K)$ is not empty (or even better, can always take a finite Galois extension of $K$)
\end{remark}

\begin{theorem}\label[theorem]{thm: fibrations p-torsion}
Let $K$ be a finitely generated field over $k$ and $X/  K$  a smooth projective geometrically connected variety. 
Then there exists an exact sequence 
$$
\Br_{\nr}(  K/k)\longrightarrow 
\ker\big(\Br_{\nr}(  K(X)/k)\rightarrow \Br(X_{ K^{\sep}})^{G_{ K}}\big) 
\rightarrow \Sha(\Pic^0_{X/ K,\red}/K),
$$
where the cokernel of the last map and the kernel of the first map have finite exponent.
\end{theorem}
\begin{proof}
By the remark above, by taking a finite Galois extension $L$ of $K$, we can assume that $X(K)$ is not empty and $\Pic(X) \rightarrow \NS(X_{K^\sep})$ is surjective. So there exists a finitely generated abelian group $N \subseteq \Pic(X)$ such that the inclusion $\Pic^0(X_{K^\sep})\oplus N \subseteq \Pic(X_{K^\sep}) $ has a finite cokernel. Thus, the natural map  $\Sha'(\Pic^0_{X/K,\red}/K)\rightarrow \Sha'(\Pic_{X/K}/K)$
has a kernel and a cokernel of finite exponent. So it suffices to show the existence of an exact sequence
$$
0\longrightarrow \Br_{\nr}(K/k)\longrightarrow 
\ker\big(\Br_{\nr}(  K(X)/k)\rightarrow \Br(X_{K^{\sep}})\big) 
\longrightarrow \Sha'(\Pic_{X/K}/K)
$$
where the last map has cokernel of finite exponent.
Let $P$ be a point in $X(K)$. Then $P$ induces $\tilde{P}\colon \Br( K(X)) \rightarrow \Br(K)$. This implies that $\Br(K)\rightarrow\Br( K(X))$ is split. Let $\sX$ and $S$ be as in the \cref{keylemma}. Then we have 
$$\bigcap_{s\in S^1}\Br_{\nr}(  K(X)/\sO_{S,s}) \subseteq \bigcap_{s\in S^1} \Br((\sX_{\sO_{S,s}})_{\reg}).$$
So $\Br_{\nr}(K/k) \subseteq\cap_{s\in S^1}\Br_{\nr}(  K(X)/\sO_{S,s})$ has a cokernel of finite exponent. Consider the  exact sequence
	\[
	0\to \Br(K) \to \ker\big(\Br(X) \to \Br(X_{K^\sep})^{G_K}\big) \to \H^1(K,\Pic(X_{K^\sep})) \to 0
	\]
 Let $a \in \Sha'(\Pic_{X/K}/K) $. There exists $b \in\ker\big(\Br(X) \to \Br(X_{K^\sep})^{G_K}\big) $ such that $\tilde{P}(b)=0$ and $b$ is mapped to $a$. Note that there is an integer $m>0$ such that for any $s\in S^1$, $m\tilde{P}(\Br_{\nr}(  K(X)/R_s) )\subseteq\Br(R_s)$, where $R_s=\sO_{S,s}$. By the lemma above, there exists a positive integer $n>0$, which is independent of $a$, and $b_s\in \ker\big(\Br_{\nr}(  K(X)/R) \rightarrow \Br(X_{K^{\sep}})\big)$ such that $b_s=na$ in $\H^1(K,\Pic(X_{K^\sep}))$ for all $s\in S^1$. Thus $mb_s-\tilde{P}(mb_s)$ is mapped to $mna$ for any $s \in S^1$. Thus $mnb=mb_s-\tilde{P}(mb_s)$ for all $s\in S^1$. This proves that $mnb\in \cap_{s\in S^1}\Br_{\nr}(  K(X)/\sO_{S,s}) $. Thus, there exists an integer $e>0$, independent of $a$, such that $eb \in \ker\big(\Br_{\nr}(  K(X)/k)\rightarrow \Br(X_{K^{\sep}})\big)$. Thus, for any $a \in \Sha'(\Pic_{X/K}/K)$, $ea$ has a preimage in $\ker\big(\Br_{\nr}(  K(X)/k)\rightarrow \Br(X_{K^{\sep}})\big)$. This proves the theorem.
\end{proof}

\begin{remark}\label{remarkforchar0}
	The theorem is also true for $K$ being a finitely generated field of characteristic~$0$. In this case, divisorial valuations for a finitely generated field are defined as valuations associated to prime divisors on all regular schemes of finite type over $\bbZ$ with the function field identified with the given finitely generated field.
	As a result, the finiteness of the exponent of the Tate--Shafarevich groups for all abelian varieties over finitely generated fields is equivalent to the finiteness of Brauer groups for all arithmetic projective 3-folds (cf.\,\cite[Thm.\,1.5]{QinBrauer2}).
\end{remark}

\begin{corollary}\label[corollary]{cor: fibrations p-torsion}
Let $S/k$ be a smooth (projective) geometrically irreducible variety with function field $K$ and  $X/K$ be a smooth projective geometrically connected variety. 
Then there exists an exact sequence
$$
0\longrightarrow V_p\Br_{\nr}(  K/k)\longrightarrow 
\ker\big(V_p\Br_{\nr}(  K(X)/k)\rightarrow V_p\Br(X_{k^{\sep}})^{G_k}\big) 
\longrightarrow V_p\Sha'(\Pic^0_{X/K,\red}/K)\longrightarrow 0.
$$
\end{corollary}

The next statement consolidates the definition $\Sha(A/K)$ from the introduction with \cref{def: Shafarevich-Tate for abelian varieties}.

\begin{proposition}\label[proposition]{twovaluations}
	Let $K$ be a finitely generated field.
	\begin{enumerate}
	\item
	The inclusion 
	$$\Br_{\nr}(K)\rightarrow \Br^{\nr}(K)$$
	has a cokernel of finite exponent.
	\item
	Let $A$ be an abelian variety over $K$. 
	Then the inclusion $\Sha(A/K)\rightarrow\Sha'(A/K)$ has a cokernel of finite exponent.
	\end{enumerate}
\end{proposition}
\begin{proof}
By de Jong's alteration theorem, there is a finite extension $L/K$ such that $L$ is the function field of a projective integral regular scheme over $\bbZ$. So $\Br_{\nr}(L)=\Br^{\nr}(L)$. 
	Then the claim follows from the fact that both corestriction maps $\Br_{\nr}(L) \rightarrow\Br_{\nr}(K)$ and $\Br^{\nr}(L) \rightarrow\Br^{\nr}(K)$ have cokernels of finite exponent and are compatible with the corestriction $\Br(L)\rightarrow\Br(K)$. 
	This shows (i).

To show (ii), we note that it suffices to show the claim for the Picard variety $\Pic^0_{X/K}$ of a smooth projective geometrically connected curve $X/K$. The \cref{lem: Br to Sha} still holds if we replace $\Br^{\nr}(  K(X)/k)$ (resp.\ $\Sha'(\Pic_{X/K}/K)$) by $\Br_{\nr}(  K(X))$ (resp.\ $\Sha(\Pic_{X/K}/K)$). Since 
	$$ \Br^{\nr}(  K(X)/k)\rightarrow\Sha'(\Pic_{X/K}/K)$$
	has a cokernel of finite exponent and $\Br_{\nr}(  K(X)/k) \rightarrow \Br^{\nr}(  K(X)/k)$  has a cokernel of finite exponent by~(i), the natural map
	$$ \Br_{\nr}(  K(X))\rightarrow\Sha'(\Pic_{X/K}/K)$$
	has a cokernel of finite exponent, too. Hence  the inclusion $\Sha(\Pic_{X/K}/K) \rightarrow \Sha'(\Pic_{X/K}/K)$ has a cokernel of finite exponent. Let $L$ be a line bundle on $X$ with $\deg(L) >0$. The natural map
	$$\Pic^0(X_{K^\sep})\oplus \bbZ[L] \rightarrow \Pic(X_{K^\sep})$$
	has a kernel and a cokernel of finite exponent. Thus, the natural maps 
	\begin{align*}
		\Sha(\Pic^0_{X/K}/K)\rightarrow\Sha(\Pic_{X/K}/K)\\
		\Sha'(\Pic^0_{X/K}/K)\rightarrow\Sha'(\Pic_{X/K}/K)
	\end{align*}
	have kernels and cokernels of finite exponents. This proves that 
	$\Sha(\Pic^0_{X/K}/K)\rightarrow\Sha'(\Pic^0_{X/K}/K)$
	has a cokernel of finite exponent.
\end{proof}

We finally can give the proof of Theorem~\hyperlink{bigtheorem}{D}, which  we recall here: 
\begin{theorem}
Let $X$ be a smooth projective geometrically connected variety over a finitely generated field $K$. Then there exists a canonical exact sequence
$$
0\longrightarrow V_\ell\Br_{\nr}(K/k)\longrightarrow 
\ker\big(V_\ell\Br_{\nr}(K(X))\rightarrow V_\ell\Br(X_{K^{\sep}})^{G_K}\big) 
\longrightarrow V_\ell\Sha(\Pic^0_{X/  K,\red})\longrightarrow 0.
$$
Moreover, if $\ell \neq \Char(K)$, the natural map
$$\Br_{\nr}(K(X))(\ell)\rightarrow \Br(X_{K^{\sep}})^{G_K}(\ell)$$
has a finite cokernel and is surjective for sufficiently large $\ell$.
\end{theorem}
\begin{proof}
By the proposition above, we can replace the unramified Brauer groups and the Tate--Shafarevich groups defined by using divisorial valuations by the ones defined by using all discrete valuations. So the first part of  Theorem~\hyperlink{bigtheorem}{D} follows from \cref{cor: fibrations p-torsion} and \cref{remarkforchar0}. It remains to show that the natural map
$$\Br_{\nr}(K(X){/k})\rightarrow \Br(X_{K^{\sep}})^{G_K}$$
has a cokernel of finite exponent. Let $k$ be the algebraic closure of $\bbQ$ in $K$. By Hironaka's theorem of resolution of singularity, we can assume that
$X/K$ can be spread out to $\rho\colon\sX\rightarrow S$, a projective surjective $k$-morphism between smooth projective geometrically connected $k$-varieties. Now applying Theorem~\hyperlink{bigthm}{E} to $\rho$, we get 
$$\Br(\sX_{k^{\sep}})^{G_k}\rightarrow \Br(X_{K^{\sep}})^{G_K}
$$
has a cokernel of finite exponent (Theorem~\hyperlink{bigthm}{E} implies that the map is surjective on the divisible parts of two groups, the claim follows from the fact that both groups are of finite exponent). It suffices to show that
$$\Br_{\nr}(K(X){/k})\rightarrow \Br(\sX_{k^{\sep}})^{G_k}$$
has a cokernel of finite exponent. This reduces the question to $K$ being a number field. If $X/K$ has a projective regular model over $\sO_K$, this was proved in \cite{QinBrauer2}. The general case can be proved by an alteration argument. By de Jong's alteration theorem, there is finite separable extension $L$ of $K(X)$ such that $K(X)$ admits a proper regular model over a finite extension of $\sO_K$. Since $\Br_{\nr}(L{/k}) \rightarrow  (\Br_{\nr}(L\otimes_K K^\sep{/k}))^{G_K}$ has a cokernel of finite exponent. Then the same claim is true for $K(X)/K$ by restriction-corestriction arguments. If $K$ is not algebraic closed in $L$. Let $K_1$ denote the algebraic closure of $K$ in $L$. By the proof of \cref{hardsur}, $(\Br_{\nr}(L\otimes_K K^\sep{/k}))^{G_K}=(\Br_{\nr}(L\otimes_{K_1} K_1^\sep{/k}))^{G_{K_1}}$. This finishes the proof.
\end{proof}

\section{The $\ell$-adic obstruction of the BSD rank conjecture} \label{sec:elladic obstruction of BSD rank}

\subsection{The BSD rank Conjecture and the $\ell$-torsion of the Tate--Shafarevich group}

Let $K$ be a finitely generated field of characteristic $p>0$  (which means finitely generated over its prime field $k=\bbF_p$).
We are considering abelian varieties $A/K$.

We start with the main result of this section, which requires several statements proved below it.
\begin{theorem}\label[theorem]{thm: bsdconj}
For an abelian variety $A$ over a finitely generated field $K$ of characteristic $p>0$ the following statements are equivalent for any prime $\ell \neq p$:
\begin{enumerate}
\item The Conjecture~\hyperlink{bsd}{B+S-D} for $A/K$ holds. 
\item The subgroup $\Sha(A/K)(\ell)$ of the Tate--Shafarevich group is finite. 
\end{enumerate}
\end{theorem}
\begin{proof}
	(i)$\implies$(ii): Write $K^\prime = Kk^\sep$. Assuming the BSD rank Conjecture for $A$, by~\cref{prop: rBSD inequality} and~\cref{splitexact} below,  we have $V_\ell\Sha_{K^\prime}(A)^{G_k}=0$. By~\cref{prop: sha}, 
	$V_\ell\Sha(A/K)=V_\ell\Sha_{K^\prime}(A)^{G_k}$, so $V_\ell\Sha(A/K)=0$. Thus $\Sha(A/K)(\ell)$ is finite because it is a cofinitely generated $\bbZ_\ell$-module by the finiteness of the $\ell$-Selmer group.
	
	(ii)$\implies$(i): Assuming that $\Sha(A/K)(\ell)$ is finite, one has $V_\ell\Sha_{K^\prime}(A)^{G_k}=0$. Then~\cref{splitexact} below implies that
	$$\H^0(S,  \sA)\otimes_\bbZ \bbQ_\ell\cong \H^1(S_{\bar{k}}, V_\ell\sA)^{G_k}$$
	and that $\H^0(S,\sA)\otimes_\bbZ \bbQ_\ell$ is a direct summand of  $\H^1(S_{\bar{k}}, V_\ell\sA)$ as $G_k$-representations. It follows that $\H^1(S_{\bar{k}}, V_\ell\sA)^{G_k}$ is direct summand of $\H^1(S_{\bar{k}}, V_\ell\sA)$ as $G_k$-representations. This implies that $\H^1(S_{\bar{k}}, V_\ell\sA)^{G_k}$ is equal to the generalised $1$-eigenspace of $\H^1(S_{\bar{k}}, V_\ell\sA)$. Then the claim follows from~\cref{prop: rBSD inequality}.
\end{proof}
\begin{lemma}\label[lemma]{splitexact}
	Let $A$ be an abelian variety defined over a finitely generated field $K$ of characteristic $p>0$. Assume that $A$ extends to an abelian scheme $\rho\colon\sA\rightarrow S$, where $S$ is a smooth geometrically connected variety over a finite field $k$ with function field $K$. Let $\ell\neq p$ be a prime. Then there is a split exact ``$\ell^\infty$-descent'' sequence of $G_k$-representations
	$$ 0\longrightarrow \H^0(S_{\bar{k}},  \sA)\otimes_\bbZ \bbQ_\ell\longrightarrow \H^1(S_{\bar{k}}, V_\ell\sA) \longrightarrow V_\ell \H^1(S_{\bar{k}}, \sA)\longrightarrow 0.$$
\end{lemma}
\begin{proof}
	Write $K^\prime$ for $K\bar{k}$. It suffices to show that $\H^0(S_{\bar{k}},  \sA)\otimes_\bbZ \bbQ_\ell=A(K^\prime)\otimes_\bbZ \bbQ_\ell$ is a direct summand of $\H^1(S_{\bar{k}}, V_\ell\sA)$ as $G_k$-representations. Note that without loss of generality, we can extend $k$ or shrink $S$. By shrinking $S$, we can assume  that $S$ is an open subvariety of a projective normal variety $\bar{S}$ over $k$. Let $V$ denote the smooth locus of $\bar{S}/k$. Then $\bar{S}-V$ has codimension $\geqslant 2$ in $\bar{S}$. Let $j$ denote the open immersion $S\hookrightarrow V$ and $\sF$ denote $V_\ell\sA$ on $S$. By \cite[Prop.\,10.1.18 (iii)]{Fu_2015}, $\R^qj_*\sF$ are constructible $\bbQ_\ell$-sheaves for all $q\geqslant 0$. Set $D=V-S$. By removing a closed subset $Z$ of codimension $\geqslant 2$ from $V$, we can assume that $D$ is a smooth divisor and $\R^1j_*\sF$ and $j_*\R^1\rho_*\bbQ_\ell$ are lisse on $D$. There is a canonical exact sequence
	$$0\longrightarrow \H^1(V_{\bar{k}},j_*\sF)\longrightarrow \H^1(S_{\bar{k}},\sF)\longrightarrow \H^0(D_{\bar{k}},\R^1j_*\sF).
	$$
	We will show that $\H^0(D_{\bar{k}},\R^1j_*\sF)$ is of weight $\geqslant 1$. This will imply that $\H^1(V_{\bar{k}},j_*\sF)$ and $\H^1(S_{\bar{k}},\sF)$ have the same generalised $1$-eigenspace of the Frobenius action. Then the question will be reduced to show the splitness of
	$$0\longrightarrow \H^0(S_{\bar{k}},  \sA)\otimes_\bbZ \bbQ_\ell \longrightarrow \H^1(V_{\bar{k}},j_*V_\ell\sA),
	$$
	and this will follow from Lemma \ref{zerocomp} and Corollary \ref{split} below.
	
	Let $\eta$ be a generic point of $D$ and $k(\eta)$ be the residue field. Then we have
	$$(\R^1j_*\sF)_{\bar{\eta}}=\H^1(K_\eta^{\sh},V_\ell A).$$
	Set $I=\Gal(K^\sep/K_\eta^{\sh})$. By \cite[Chap.\,I, Lem.\,2.18]{Milne_2006}, there is an isomorphism of $G_{k(\eta)}$-representations
	$$ \H^1(K_\eta^{\sh},V_\ell A)\cong (V_\ell A)_I(-1).$$
	Since $(V_\ell A)_I=(\H^1(A_{K^\sep},\bbQ_\ell)^I)^\vee$, it follows that
	$$\H^1(K_\eta^{\sh},V_\ell A)\cong (\H^1(A_{K^\sep},\bbQ_\ell)^I)^\vee(-1)$$
	as $G_{k(\eta)}$-representations. Let $s\in \overline{\{\eta\}}$ be a closed point. By assumption, $j_*\R^1\rho_*\bbQ_\ell$ is lisse on $D$. So $G_{  K}$ acts on $(j_*\R^1\rho_*\bbQ_\ell)_{\bar{\eta}}=\H^1(A_{K^\sep},\bbQ_\ell)^I$. By \cite[Cor.\,1.8.9]{Deligne_1980}, the $G_{  K}$-action on $\H^1(A_{K^\sep},\bbQ_\ell)^I$ is of weight $\leqslant 1$. Thus, the $G_{  K}$-action on $(\H^1(A_{K^\sep},\bbQ_\ell)^I)^\vee(-1)$ is of weight $\geqslant 1$. It follows that $(\R^1j_*\sF)_{\bar{s}}$ is of weight $\geqslant 1$. Since this holds for all closed points in $D$, it follows that $\H^0(D_{\bar{k}},\R^1j_*\sF)$ is of weight $\geqslant 1$. This completes the proof.
\end{proof}

\subsection{The N\'eron--Tate height pairing and Yoneda pairing}

Let $S$ be a projective normal geometrically connected variety over a finite field $k$ with function field $K$. Let $V\subseteq S$ be a regular open subscheme with $S-V$ of codimension $\geqslant 2$. Let $A$ be an abelian variety over $K$. Let $\sY \subseteq V$ be an open dense subscheme such that $A$ extends to an abelian scheme $\rho\colon\sA\rightarrow \sY$. Let $j\colon \sY\hookrightarrow V$ be the inclusion and $K^\prime$ denote $K\bar{k}$. The aim of this section is to show the splitness of
\begin{equation}\label{aimsplit}
	0\longrightarrow A(K^\prime)\otimes_{\bbZ}\bbQ_\ell \longrightarrow \H^1(V_{k^{\sep}},j_*V_\ell\sA).
\end{equation}
In the case that $\dim S=1$, $V$ is equal to $S$. $j_*\sA$ is the \'etale sheaf represented by a N\'eron model of $A$. The claim was proved by Schneider~\cite{Schneider_1982}. In the case that $A$ extends to an abelian scheme over a smooth projective variety $S$ over $k$, we have $j_*\sA=\sA$. This case was proved by Keller~\cite{Keller_2019}. The idea is to construct compatible $G_k$-equivariant pairings
\[
\begin{tikzcd} 
	\H^0(V,\sA^0)\otimes_{\bbZ}\bbQ \arrow[d,"-\delta"] \arrow[r,"\times", phantom] & \Ext^1_{V_{\fppf}}(\sA^0,\bbG_m)\otimes_{\bbZ}\bbQ \ar[d, "\varprojlim r_n"]\ar[r]  & \bbQ \ar[d]  \\
	\H^1(V_{\bar{k}},V_\ell\sA^0) \arrow[r,"\times", phantom]  & \varprojlim_{n} \Ext^1_{(\ell^n)-V}(\sA^0[\ell^n],\mu_{\ell^n})\otimes_{\bbZ_\ell} \bbQ_\ell \ar[r] & \bbQ_\ell
\end{tikzcd}
\]
where $\sA^0$ is a smooth group scheme over $V$ with connected fibres ($\sA^0$ is taken to be the identity component of a N\'eron model of $A$ for two special cases above) and the top pairing is left non-degenerate. The the splitness of (\ref{aimsplit}) will follow from Lemma \ref{splitrep}. In the following, we will explain the construction of $\sA^0$ and the two pairings.

\subsubsection {N\'eron model}
\begin{lemma}
	Let $A$ be an abelian variety over a finitely generated field $K$ of characteristic $p>0$. Let $S$ be a projective normal geometrically connected variety over a finite field $k$ with function field $K$. Let $S^1$ denote the set of points of codimension $1$. Then there exists a regular open subscheme $V\subseteq S$ with codimension $(S-V)\geqslant 2$ and a smooth commutative $V$-group scheme $\sA$ of finite type satisfying the following assumptions:
	\begin{itemize}
		\item[(i)] $A$ is the generic fibre of $\sA$ and the restriction of $\sA$ to $\Spec\sO_{S,s}$ is a N\'eron model of $A/K$ for all $s \in S^1$.
		\item[(ii)] $\sA$ admits an open subgroup scheme $\sA^0$ which equals the identity component of $\sA$ when restricted to $ \Spec\sO_{S,s}$ for all $s \in S^1$ and there exists a closed reduced subscheme $Y$ of $V$ such that $\sA$ is abelian over $V-Y$ and the quotient of $\sA_Y$ by $\sA^0_Y$ is a finite \'etale group scheme over $Y$.
		\item[(iii)] All fibres of $\sA^0\rightarrow V$ are geometrically connected.
	\end{itemize}
	
\end{lemma}
\begin{proof}
	For the case  $\dim(S)=1$, we can take $V=S$ and take $\sA$ to be a N\'eron model of $A/K$ over $S$ and $\sA^0$ to be its identity component. In general, if $\dim(S)>1$, the N\'eron model of $A$ over $S$ need not exist. The idea is that extending $A$ to a N\'eron model over  
	$\Spec(\sO_{S,s}) $ for all $s\in S^1$, then spread out.  Since a N\'eron model of an abelian variety over a discrete valuation ring is of finite type, a N\'eron model of $A$ over $\Spec(\sO_{S,s})$ can spread out to a smooth group scheme of finite type over an open neighborhood of $s\in S^1$. There exists a regular open dense subscheme $\sY$ of $S$ such that $A$ extends to an abelian scheme $\sA_\sY\rightarrow \sY$. By \cite[Thm.\,3.3]{Keller_2016}, $\sA_\sY \rightarrow \sY$  is a N\'eron model of $A$ over $\sY$.  Since $S^1-\sY$ is finite, we can glue N\'eron models of $A$ over $\Spec( \sO_{S,s})$ for all $s\in S^1-\sY $ with $\sA_\sY\rightarrow \sY$ to get a smooth commutative group scheme $\sA$ of finite type over some open subset $V$ that contains $S^1$ satisfying the condition (i). For each $s \in S^1-\sY$, by removing the closure of the complement of the identity component of $\sA_s$, we get an open subscheme  $\sA^0$ of $\sA$, then shrink $V$, this subscheme becomes an open subgroup scheme since $\sA_s/\sA^0_s$ is finite \'etale over $  K$. Taking $Y=V-\sY$, by shrinking $V$, we will show that the quotient $\sA_Y/\sA_Y^0$ exists and is finite \'etale over $Y$. Let $s\in Y\cap S^1$. There is an exact sequence
	$$ 0\longrightarrow \sA_s^0\longrightarrow \sA_s\longrightarrow \sF_s \longrightarrow 0,$$
	where $\sF_s$ is a finite \'etale group scheme over $K$. By shrinking $Y$, $\sF_s$ extends to a finite \'etale group scheme $\sF$ over $Y$ and the map $\sA_s\rightarrow \sF_s$ extends to a faithful flat morphism of finite type $\sA_Y\rightarrow \sF$. By shrinking $Y$ further, the kernel of $\sA_Y\rightarrow \sF$ can be identified with $\sA^0_Y$. We get an exact sequence of fppf sheaves on $Y$ 
	$$0\longrightarrow\sA_Y^{0}\longrightarrow \sA_Y\longrightarrow \sF \longrightarrow 0.$$
	Therefore, $\sA^0\rightarrow V$ satisfies condition (ii). Since $\sA^0_{s}$ is geometrically connected for all $s\in S^1$, by shrinking $Y$, we can assume that all fibres of $\sA^0_{Y}\rightarrow Y$ are geometrically connected.
\end{proof}
In the lemma below, we will show that the $\sA^0$ constructed as above satisfies
$$V_\ell\sA^0\cong  j_*V_\ell(\sA|_\sY).$$
As a result, there is a canonical isomorphism
$$ \H^1(V_{\bar{k}}, j_*V_\ell(\sA|_\sY))\cong \H^1(V_{\bar{k}},V_\ell\sA^0).$$ 

\begin{lemma}\label{zerocomp}
	Notations as above, let $\sY\subseteq V$ be an open dense subset such that $\sA$ is an abelian scheme over $\sY$. Let $j\colon \sY\rightarrow V$ be the inclusion.  Denote by $T_\ell\sA^0 $ (resp. $T_\ell\sA$) the inverse system of $\ell$-torsion sheaves $(\sA^0[\ell^n])_{n\geqslant 0}$ (resp.\ $(\sA[\ell^n])_{n\geqslant 0}$). By shrinking $V$, there is a canonical isomorphism between A--R $\ell$-adic sheaves (cf.~\cite[Chap.\,10]{Fu_2015})
	$$ T_\ell\sA^0\cong j_*(T_\ell\sA|_{\sY}).$$
\end{lemma}

\begin{proof}
	By the lemma below, $T_\ell \sA^0$ is an A--R $\ell$-adic sheaf. There is an exact sequence
	$$0\longrightarrow \sA^0[\ell^n]\longrightarrow \sA[\ell^n] \longrightarrow \sF[\ell^n] \longrightarrow 0.$$
	Since $\sF$ is finite \'etale, $T_\ell\sF$ is A--R zero. Therefore $T_\ell\sA^0$ is A--R isomorphic to $T_\ell\sA$. Since $T_\ell\sA|_{\sY}$ is a $\ell$-adic sheaf,  $j_*(T_\ell\sA|_{\sY})$ is also A--R $\ell$-adic. Since $\sA$ is a N\'eron model of $A$ when restricted to $ \Spec\sO_{S,s}$ for all $s \in S^1$, we have $\sA\cong j_*(\sA|_{\sY})$ over $ \Spec\sO_{S,s}$. It follows that
	$$\sA[\ell^n]_{\bar{s}}\cong (j_*(\sA[\ell^n]|_{\sY}))_{\bar{s}}.$$
	It follows that the kernel and cokernel of 
	$$T_\ell\sA\longrightarrow j_*(T_\ell\sA|_{\sY})$$
	are supported on a closed subset of codimension $\geqslant 2$. Removing this closed subset from $V$, we have 
	\[ T_\ell\sA^0\cong j_*(T_\ell\sA|_{\sY}). \qedhere\]
\end{proof}

\begin{lemma}
	Let $S$ be a Noetherian scheme and $G$ be a smooth commutative group scheme over $S$ of finite type. Let $\ell$ be a prime invertible on $S$. Then the inverse system of $\ell$-torsion sheaves $T_\ell G=(G[\ell^n])_{n\geqslant 1}$ is A--R $\ell$-adic. If all fibres of $G/S$ are geometrically connected, then 
	$$G\stackrel{\ell^n}{\longrightarrow}G$$
	is surjective as a morphism between \'etale sheaves on $S$.
\end{lemma}

\begin{proof}
	Firstly, note that $G[\ell^n]$ is \'etale and of finite type over $S$. Thus the \'etale sheaf represented by $G[\ell^n]$ is constructible.
	
	Secondly, let $f\colon S_1\rightarrow S$ be a surjective morphism between Noetherian schemes. If the statement holds for $G_{S_1}/S_1$, then it holds for $G/S$. Set $\sF_n=G[\ell^n]$ and $\sF=(\sF_n)$. Since $f^*$ is an exact functor, one can show that if $f^*\sF$ satisfies the criterion of \cite[Prop.\,10.1.1]{Fu_2015}, then so does $\sF$.
	
	Thirdly, we might assume that $S$ is integral. By Noetherian induction, it suffices to prove the statement for an open sense subscheme of $S$. Let $K$ denote the function field of $S$. $G_{\bar{K}}$ is smooth and of finite type over $\bar{K}$. Let $G^0_{\bar{K}}$ be the identity component of $G_{\bar{K}}$, then there is an exact sequence 
	$$0\longrightarrow G^0_{\bar{K}} \longrightarrow G_{\bar{K}} \longrightarrow \pi_0(G)\longrightarrow 0,$$
	where $\pi_0(G)$ is a finite \'etale group scheme. By Chevalley's theorem, $G^0_{\bar{K}}$ is an extension of an abelian variety $A$ by $\bbG_m^r\times \bbG_a^\sep$, i.e.,
	$$0\longrightarrow \bbG_m^r\times \bbG_a^\sep\longrightarrow G^0_{\bar{K}}\longrightarrow A\longrightarrow 0.$$
	The above exact sequences can descend to a finite extension $L$ of $K$. Choose a model $S_1$ for $L$. By shrinking $S_1$ and $S$, we can assume that there exists a flat surjective morphism $S_1\rightarrow S$ of finite type whose function field extension corresponds to $L/K$ and exact sequences
	$$0\longrightarrow G^0_{S_1} \longrightarrow G_{S_1} \longrightarrow \pi_0(G)\longrightarrow 0$$
	and 
	$$0\longrightarrow \bbG_m^r\times \bbG_a^\sep \longrightarrow G^0_{S_1}\longrightarrow \sA_{S_1}\longrightarrow 0,$$
	where $\pi_0(G)$ is finite \'etale over $S_1$ and $\sA_{S_1}$ is an abelian scheme. We get exact sequences of systems of $\ell$-torsion sheaves
	$$0\longrightarrow T_\ell G^0_{S_1} \longrightarrow T_\ell G_{S_1} \longrightarrow T_\ell\pi_0(G)\longrightarrow 0$$
	and 
	$$0\longrightarrow T_\ell(\bbG_m^r\times \bbG_a^\sep) \longrightarrow T_\ell G^0_{S_1}\longrightarrow T_\ell\sA_{S_1}\longrightarrow 0.$$
	Since $T_\ell(\bbG_m^r\times \bbG_a^\sep)=T_\ell \bbG_m^r$ and $T_\ell\sA_{S_1}$ are $\ell$-adic sheaves, by \cite[Prop.\,10.1.7\,(iii)]{Fu_2015}, $T_\ell G_{S_1}$ is A--R $\ell$-adic. It follows that $T_\ell G$ is A--R $\ell$-adic on $S$. This proves the first claim.
	
	For the second claim, we might assume $S=\Spec( R)$ where $R$ is a strict local ring with residue field $k$. It suffices to show that $G(S)$ is $\ell$-divisible. Let $\tau\colon S\rightarrow G$ be a section. $\ell^{-1}(\tau)$ is \'etale over $S$. $\ell^{-1}(\tau)(k)$ is not empty since $G(k)$ is $\ell$-divisible. By Hensel's Lemma, it follows that $\ell^{-1}(\tau)(S)$ is not empty. Thus, $G(S)$ is $\ell$-divisible.
\end{proof}

\subsubsection{Yoneda pairing}

Following notations in \cite[\S\,2]{Schneider_1982}, let $V_{\fppf}$ (resp. $(\ell^n)-V $) denote the category of fppf-sheaves (resp.\ \'etale sheaves of $\bbZ/\ell^n\bbZ$-modules) on $V$.

\begin{lemma}\label{commu}
	Notations as before, assuming that $V$ is regular, then there is a commutative diagram
	\begin{displaymath}
		\xymatrix{ 
			\H^0(V,\sA^0) \ar[d]^{-\delta} &\times &\Ext^1_{V_{\fppf}}(\sA^0,\bbG_m)\ar[d]^{r_n}\ar[r]  & \H^1(V,\bbG_m) \ar[d] \\
			\H^1(V,\sA^0[\ell^n]) &\times  & \Ext^1_{(\ell^n)-V}(\sA^0[\ell^n],\mu_{\ell^n}) \ar[r] & \H^2(V,\mu_{\ell^n}) }
	\end{displaymath}
	where $\delta$ is induced by the exact sequence $0\rightarrow \sA^0[\ell^n]\rightarrow \sA^0 
	\xrightarrow{\ell^n} \sA^0\rightarrow 0$ and $r_n$ is defined by 
	$$ (0\longrightarrow \bbG_m\longrightarrow \sX\longrightarrow \sA^0\longrightarrow 0)\mapsto (0\longrightarrow\mu_{\ell^n}\longrightarrow \sX[\ell^n] \longrightarrow \sA^0[\ell^n]\longrightarrow 0).
	$$
	The natural map 
	\begin{equation}\label{iso}
		\Ext^1_{V_{\fppf}}(\sA^0[\ell^n],\bbG_m) \longrightarrow  \Ext^1_{(\ell^n)-V}(\sA^0[\ell^n],\mu_{\ell^n})
	\end{equation}
	defined by
	$$ (0\longrightarrow \bbG_m\longrightarrow \sX\longrightarrow \sA^0[\ell^n]\longrightarrow 0)\mapsto (0\longrightarrow\mu_{\ell^n}\longrightarrow \sX[\ell^n] \longrightarrow \sA^0[\ell^n]\longrightarrow 0)
	$$
	is an isomorphism and it induces a natural map 
	$$  \Ext^1_{(\ell^n)-V}(\sA^0[\ell^{n+1}],\mu_{\ell^{n+1}})\longrightarrow   \Ext^1_{(\ell^n)-V}(\sA^0[\ell^n],\mu_{\ell^n})
	$$
	compatible with pairings in the above diagram.  As a result, by taking limits, there is a commutative diagram
	
	\begin{displaymath}
		\xymatrix{ 
			\H^0(V,\sA^0) \ar[d]^{-\delta} &\times &\Ext^1_{V_{\fppf}}(\sA^0,\bbG_m)\ar[d]^{\varprojlim r_n}\ar[r]  & \H^1(V,\bbG_m) \ar[d] \\
			\H^1(V,T_\ell\sA^0) &\times  & \varprojlim_n \Ext^1_{(\ell^n)-V}(\sA^0[\ell^n],\mu_{\ell^n}) \ar[r] & \H^2(V,\bbZ_\ell(1)) }
	\end{displaymath}
\end{lemma}

\begin{proof}
	The two pairings in the first diagram are defined as Yoneda product of extensions and the commutativity of the diagram can be checked explicitly.
	We will construct a natural map 
	$$\Ext^1_{(\ell^n)-V}(\sA^0[\ell^n],\mu_{\ell^n}) \longrightarrow \Ext^1_{V_{\fppf}}(\sA^0[\ell^n],\mu_{\ell^n})$$
	which induces a map 
	$$\Ext^1_{(\ell^n)-V}(\sA^0[\ell^n],\mu_{\ell^n}) \longrightarrow \Ext^1_{V_{\fppf}}(\sA^0[\ell^n],\bbG_m). $$
	And it is easy to check the induced map is the inverse of (\ref{iso}) by definition. Let $(0\rightarrow\mu_{\ell^n}\rightarrow \sY \rightarrow \sA^0[\ell^n]\rightarrow 0) \in \Ext^1_{(\ell^n)-V}(\sA^0[\ell^n],\mu_{\ell^n})$. $\sY$ may be regarded as a $\mu_{\ell^n}$-torsor over $\sA^0[\ell^n]$ and so it is representable by an \'etale scheme of finite type over $V$ (cf.~\cite[Chap.\,III, Thm.\,4.3]{Milne_1980}). Therefore the exact sequence can be regarded as an exact sequence of fppf sheaves on $V$ which gives an element in $\Ext^1_{V_{\fppf}}(\sA^0[\ell^n],\mu_{\ell^n})$. This defines a map
	$$\Ext^1_{(\ell^n)-V}(\sA^0[\ell^n],\mu_{\ell^n}) \longrightarrow \Ext^1_{V_{\fppf}}(\sA^0[\ell^n],\mu_{\ell^n}).$$
	By composing with the natural map
	$$\Ext^1_{V_{\fppf}}(\sA^0[\ell^n],\mu_{\ell^n})\longrightarrow \Ext^1_{V_{\fppf}}(\sA^0[\ell^n],\bbG_m),$$
	we get the desired map. Define
	$$\Ext^1_{V_{\fppf}}(\sA^0[\ell^{n+1}],\bbG_m) \longrightarrow \Ext^1_{(\ell^n)-V}(\sA^0[\ell^n],\mu_{\ell^n})$$
	as
	$$ (0\longrightarrow \bbG_m\longrightarrow \sX\longrightarrow \sA^0[\ell^{n+1}]\longrightarrow 0)\mapsto (0\longrightarrow\mu_{\ell^n}\longrightarrow \sX[\ell^n] \longrightarrow \sA^0[\ell^n]\longrightarrow 0).
	$$
	Then the composition 
	$$\Ext^1_{(\ell^n)-V}(\sA^0[\ell^{n+1}],\mu_{\ell^{n+1}}) \longrightarrow \Ext^1_{V_{\fppf}}(\sA^0[\ell^{n+1}],\bbG_m) \longrightarrow \Ext^1_{(\ell^n)-V}(\sA^0[\ell^n],\mu_{\ell^n}) $$
	gives a map compatible with the first diagram.
\end{proof}

\subsubsection{Compatibility between the height pairing and Yoneda pairing}

\begin{theorem}\label{height}
	Notations as before. Let $S$ be a projective normal geometrically connected variety over a finite field $k$. Fix a closed $k$-immersion $\iota\colon S\hookrightarrow \bbP^m_{k}$. Let $V$ be a regular open subscheme of $S$ with codim $(S-V)\geqslant 2$ satisfying conditions in the lemma above. Define a degree map $\Pic(S)\rightarrow  \bbZ$ by sending a prime Weil divisor to its degree as a subvariety in $\bbP^m_{k}$. Composing it with the Yoneda pairing in the lemma above, we get a pairing 
	\begin{equation}\label{keypairing}
		\H^0(V,\sA^0) \times \Ext^1_{V_{\fppf}}(\sA^0,\bbG_m) \longrightarrow \bbZ.
	\end{equation}
	There is a natural inclusion $\H^0(V,\sA^0)\hookrightarrow A(K)$ and a natural map $\Ext^1_{V_{\fppf}}(\sA^0,\bbG_m) \rightarrow  \Ext^1_{K}(A,\bbG_m)=A^t(K)$, where $A^t$ is the dual abelian variety of $A$. Then, up to a normalizing constant, the above pairing is compatible 
	with the N\'eron--Tate height pairing defined in \cite[Cor.\,9.17]{Conrad_2006}
	$$A(K)\times A^t(K) \longrightarrow \bbR.$$
\end{theorem}

\begin{remark}
	This is a generalisation of results in \cite[\S\,3]{Keller_2019} and \cite[\S\,3]{Schneider_1982}. The proof is just an imitation of arguments in \cite[\S\,2]{Bloch_1980}, \cite[\S\,3]{Schneider_1982} and \cite[\S\,3]{Keller_2019}. As a consequence, the pairing (\ref{keypairing}) is left non-degenerate after tensoring $\bbQ$.
\end{remark}

\begin{proof}
	Firstly, we define the adele ring associated to $S$ (cf.~\cite[\S 3.1]{Keller_2019}). Define for $T\subset S^1$ finite the $T$-adele ring of $S$ as the restricted product 
	
	$$\bbA_{K,T}=\prod_{s\in T}K_s\times \prod_{s\in S^1\backslash T}\widehat{\sO}_{S,s},$$
	where $\widehat{\sO}_{S,s}$ is the completion of the local ring $\sO_{S,s}$ and $K_s$ is the quotient field of $\widehat{\sO}_{S,s}$, and the adele ring of $S$ as
	
	$$\bbA_K=\varinjlim_{T\subset S^1}\bbA_{K,T}.$$
	For each $s\in S^1$, we define $\deg(s)$ as the degree of the close subvariety $\overline{\{s\}}$ in $\bbP_k^m$. This gives an absolute value on $K_s$ i.e., $|\cdot|_s=q^{-\deg(s) \cdot v_s(\cdot)}$, where $q$ denotes the cardinality of $k$. Since a principal Weil divisor has degree zero, the absolute value $|\cdot|_s$ satisfies product formula. We call the field  $K$ equipped with the set of absolute values $|\cdot|_s$ a generalised global field (cf.~\cite[Def.\,8.1]{Conrad_2006}).
	
	Secondly, we describe the first pairing explicitly. Let $\tilde{a}=(0\rightarrow \bbG_m\rightarrow \sX \rightarrow \sA^0\rightarrow 0) \in \Ext^1(\sA^0,\bbG_m)$. One can think $\sX$ as a $\bbG_m$-torsor on $\sA^0$ and so it is representable by a smooth commutative $S$-group scheme of finite type (cf.~\cite[Chap.\,III, Thm.\,4.3]{Milne_1980}). By Hilbert's theorem $90$, there are exact sequences
	$$0\longrightarrow \bbG_m(K) \longrightarrow \sX(K) \longrightarrow A(K)\longrightarrow 0$$
	and 
	$$0\longrightarrow \bbG_m(\bbA_K) \longrightarrow \sX(\bbA_K) \longrightarrow \sA^0(\bbA_K)\longrightarrow 0.$$
	There is a natural homomorphism, the \emph{logarithmic modulus map},
	$$l\colon\bbG_m(\bbA_K)\longrightarrow \log q\cdot\bbZ \subseteq \bbR, (a_s)\mapsto
	\sum_{s\in S^1}\log |a_s|_s=-\log q\cdot \sum_{s \in S^1}\deg(s)\cdot v_s(a_s).$$
	By the product formula, $l(\bbG_m(K))=0$. Define $G_m^1$ as the kernel of $l$, and $\sX^1$ as 
	$$ \sX^1=\{ a \in \sX(\bbA_K): \exists n \in \bbZ_{\geqslant 1}, na\in\widetilde{\sX^1}(\bbA_K)\},$$
	the rational saturation of $\widetilde{\sX^1}$ with 
	$$ \widetilde{\sX^1}=G_m^1\cdot\prod_{s\in S^1}\sX(\widehat{\sO}_{S,s})\subseteq \sX(\bbA_K).
	$$
	Then there exists a unique extension $l_{\tilde{a}}\colon\sX(\bbA_K)\rightarrow  \bbR$ of $l$ vanishing on $\sX^1$ (cf.~\cite[Lem.\,3.1.4]{Keller_2019} or \cite[Lem.\,1.8]{Bloch_1980}). It induces by restriction to $\sX(K)$ a homomorphism 
	$$ l_{\tilde{a}}\colon A(K)\longrightarrow \bbR.$$
	Next, we will show that for any $a\in\sA^0(S)\subseteq A(K)$, we have
	$$l_{\tilde{a}}(a)=-\log q\cdot\deg(a\vee \tilde{a}),$$
	where $\vee$ denotes the Yoneda pairing $\H^0(V,\sA^0) \times \Ext^1_{V_{\fppf}}(\sA^0,\bbG_m) \rightarrow  \H^1(V,\bbG_m)$. By definition, $a\vee \tilde{a}$ is defined by the following commutative diagram
	\begin{displaymath}
		\xymatrix{
			a\vee \tilde{a}\colon &	0 \ar[r] & \bbG_m \ar[r] \ar[d]^{\id} & \sY \ar[d]\ar[r] & \bbZ \ar[r]\ar[d]^{a} & 0\\
			\tilde{a}\colon &0 \ar[r] & \bbG_m \ar[r] & \sX\ar[r] & \sA^0 \ar[r] & 0}
	\end{displaymath}
	By composition, one gets an extension
	$$
	l_{a\vee \tilde{a}}\colon\sY(\bbA_K)\longrightarrow \sX(\bbA_K)\stackrel{l_{\tilde{a}}} \longrightarrow \bbR
	$$
	of $l\colon\bbG_m(\bbA_K)\rightarrow  \bbR$ to $\sY(\bbA_K)$, which induces because of $l(\bbG_m(K))=0$ in the exact sequence $a\vee \tilde{a}$ by restriction to $\sY(K)$ a homomorphism 
	$$ l_{a\vee \tilde{a}}\colon\bbZ\stackrel{a} \longrightarrow A(K)\stackrel{l_{\tilde{a}}}\longrightarrow\bbR ,$$
	so one obviously has $l_{\tilde{a}}(a)=l_{a\vee\tilde{a}}(1)$. Since 
	$$l_{\tilde{a}}\Big(\prod_{s\in S^1}\sX(\widehat{\sO}_{S,s})\Big)=0,$$
	hence
	$$\ell_{a\vee \tilde{a}}\Big(\prod_{s\in S^1}\sY(\widehat{\sO}_{S,s})\Big)=0.$$
	Set $e=a\vee \tilde{a}$, $e$ is represented by an exact sequence of fppf sheaves on $V$
	$$0 \longrightarrow  \bbG_m \longrightarrow \sY \longrightarrow \bbZ \longrightarrow 0,$$
	and $l_e$ is an extension of $l$ and vanishes on $\prod_{s\in S^1}\sY(\widehat{\sO}_{S,s})$. By \cite[Lem.\,3.3.2]{Keller_2019} (replacing $X$ by $V$ in his argument) or \cite[Lem.\,12]{Schneider_1982}, 
	$$l_{\tilde{a}}(a)=l_e(1)=-\log q\cdot \deg(e).$$
	
	Thirdly, we will show that the pairing $h(a,\tilde{a})\colonequals l_{\tilde{a}}(a)$ can be written as a sum of the local N\'eron pairings which therefore coincides with the canonical N\'eron--Tate height pairing by \cite[p.\,307, Cor.\,9.5.14]{BombieriGubler_2007}. The proof for our case is exactly same as the argument in \cite[\S\,2]{Bloch_1980} and \cite[\S\,3.4]{Keller_2019}.  Given an extension $(0\rightarrow \bbG_m\rightarrow \sX\rightarrow \sA^0\rightarrow 0) \in \Ext^1_{V_{\fppf}}(\sA^0,\bbG_m)$, its restriction to $\Spec( \sO_{S,s})$ for $s \in S^1$ is still an extension. Since $\sO_{S,s}$ is a Dedekind domain, by\cite[p.\,53, Lem.\,5.1]{MazurMessing_2006}, the push-forward of the sheaf 
	$$\sheafext^1_{(\Spec( \sO_{S,s}))_{\fppf}}(\sA^0,\bbG_m)$$
	to the smooth site of $\Spec( \sO_{S,s})$ is represented by the N\'eron model of $A^t\colonequals \Pic^0_{A/K}$ over 
	$\Spec( \sO_{S,s})$. By \cite[Lem.\,9]{Schneider_1982}, 
	$\sheafhom_{(\Spec( \widehat{\sO}_{S,s}))_{\fppf}}(\sA^0,\bbG_m)=0$, it follows that
	$$
	A^t(K_s)=\H^0(\widehat{\sO}_{S,s},\sheafext^1_{(\Spec(\widehat{\sO}_{S,s}))_{\fppf}}(\sA^0,\bbG_m))=\Ext^1_{(\Spec( \widehat{\sO}_{S,s}))_{\fppf}}(\sA^0,\bbG_m).
	$$
	Therefore, the restriction of $(0\rightarrow \bbG_m\rightarrow \sX\rightarrow \sA^0\rightarrow 0)$ to $\Spec( \widehat{\sO}_{S,s})$ gives an element in $A^t(K_s)$. Then one can define the local N\'eron pairings as \cite[(2.6)]{Bloch_1980} for each $s\in S^1$. By \cite[Thm.\,2.7]{Bloch_1980} or \cite[Thm.\,3.4.2]{Keller_2019}, $h(\cdot ,\cdot)$ is equal to the sum of the local N\'eron pairings. This completes the proof of the theorem.
\end{proof}

\begin{corollary}\label{split}
	Notations as before. Set $d=\dim(V)$ and let $K^\prime$ denote $K\bar{k}$. Then, one can shrink $V$ (with $\codim(S-V) \geqslant 2$) such that 
	$$ \H^0(V_{\bar{k}},\sA^0)\otimes_\bbZ \bbQ_\ell=A(K^\prime)\otimes_\bbZ\bbQ_\ell$$
	and the injective map of $G_k$-representations
	$$\H^0(V_{\bar{k}},\sA^0)\otimes_{\bbZ} \bbQ_\ell \longrightarrow \H^1(V_{\bar{k}},V_\ell\sA^0)$$
	is split.
\end{corollary}
\begin{proof}
	By the Lang--N\'eron theorem \cite[Thm.\,7.1]{Conrad_2006}, the quotient group $A(K^\prime)/\Tr_{K/k}(A)(\bar{k})$ is finitely generated where $\Tr_{K/k}(A)$ is a $K/k$-trace of $A/K$. Since $\Tr_{K/k}(A)(\bar{k})$ is a $\ell$-divisible torsion group, we have that
	$$A(K^\prime)\hat{\otimes}\bbQ_\ell=A(K^\prime)\otimes_\bbZ \bbQ_\ell$$
	has finite dimension over $\bbQ_\ell$. Thus, there exists a finite extension $l/k$ such that $A(Kl)\otimes_\bbZ \bbQ_\ell=A(K^\prime)\otimes_\bbZ\bbQ_\ell$. Replacing $k$ by $l$, we might assume $A(K)\otimes \bbQ_\ell=A(K^\prime)\otimes \bbQ_\ell$. Since $A(K)$ is finitely generated, by shrinking $V$, we can assume $\sA(V)=A(K)$. Since $\sA(V)/\sA^0(V)$ is finite, we can shrink $V$ further such that $\sA^0(V)\otimes_\bbZ\bbQ_\ell=A(K)\otimes_\bbZ\bbQ_\ell$. Assume that $A^t(K)$ is generated by  finitely many elements $\tilde{a}_i$. Since 
	$$A^t(K_s)=\Ext^1_{(\Spec( \sO_{S,s}))_{\fppf}}(\sA^0,\bbG_m)
	$$ 
	for any $s\in S^1$, so for a given $\tilde{a}\in A^t(K)$, it corresponds to an extension $\tilde{a}_s \colon(0\rightarrow \bbG_m\rightarrow \sX\rightarrow \sA^0\rightarrow 0)$ over 
	$\Spec( \sO_{S,s})$ whose restriction to $\Spec( K)$ is independent of $s$.  One can think of $\sX$ as a $\bbG_m$-torsor on $\sA^0$ and so it is representable by a smooth commutative group scheme of finite type over $\Spec( \sO_{S,s})$ (cf. \cite[Chap.\,III, Thm.\,4.3]{Milne_1980}). Therefore, by shrinking $V$, one can glue the $\tilde{a}_s$ to get an extension $(0\rightarrow \bbG_m\rightarrow \sX\rightarrow \sA^0\rightarrow 0)$ over $V$. Thus, by shrinking $V$, we can assume that the natural map
	$$\Ext^1_{V_{\fppf}}(\sA^0,\bbG_m)\rightarrow \Ext^1_{K}( A,\bbG_m)=A^t(K)$$
	is surjective. There is a commutative diagram 
	\[
	\begin{tikzcd} 
		\H^0(V,\sA^0) \arrow[d,"-\delta"] \arrow[r,"\times", phantom] & \Ext^1_{V_{\fppf}}(\sA^0,\bbG_m)\ar[d, "\varprojlim r_n"]\ar[r]  & \H^1(V,\bbG_m) \ar[d] \ar[r, "\deg(\cdot)"] & \bbZ \ar[d]  \\
		\H^1(V_{\bar{k}},V_\ell\sA^0) \arrow[r,"\times", phantom]  & \varprojlim_n \Ext^1_{(\ell^n)-V_{\et}}(\sA^0[\ell^n],\mu_{\ell^n})\otimes_{\bbZ_\ell}\bbQ_\ell \ar[r] & \H^2(V_{\bar{k}},\bbQ_\ell(1))\ar[r,"\triv(\cdot\cap \eta^{d-1})"] & \bbQ_\ell
	\end{tikzcd}
	\]
	where $\eta$ denotes the cycle class of $\sO_{S}(1)$ and $\triv$ denotes the trace map  $\H^{2d}(S_{\bar{k}},\bbQ_\ell(d))\cong \bbQ_\ell$. The commutativity of the second square follows from the definition of degree. By Lemma \ref{commu}, the first square commutes. By Theorem \ref{height}, the top pairing is compatible with the N\'eron--Tate pairing which is non-degenerate after tensoring $\bbQ$ (cf. \cite[Cor.\,9.17]{Conrad_2006}). Assuming the surjectivity of $\Ext^1_{V_{\fppf}}(\sA^0,\bbG_m)\rightarrow \Ext^1_{K}( A,\bbG_m)=A^t(K)$, then the top pairing is left non-degenerate after tensoring $\bbQ_\ell$. By Lemma \ref{splitrep}, the injective map of $G_k$-representations
	$$\H^0(V_{\bar{k}},\sA^0)\otimes_\bbZ \bbQ_\ell \stackrel{\delta}\longrightarrow \H^1(V_{\bar{k}},V_\ell\sA^0)$$
	is split.
\end{proof}

\subsection{Applications}

We finish the section with two applications:
\begin{corollary}\label[corollary]{cor1}
Let $K$ a finitely generated field in characteristic $p>0$. 
Let $X/K$ be a smooth projective geometrically connected variety 
with a smooth model $\sX/\sY$.
	Then the analytic Tate Conjecture~\hyperlink{BS3}{B+S-D+2} for  $\sX$, that is $\T^1_{\an}(\sX)$, is equivalent to 
	\begin{align*}
	&\text{ Conjecture~\hyperlink{BS3}{B+S-D+2} for $\sY$, that is, $\T^1_{\an}(\sY)$ }   \\
	+ & \text{ Conjecture~\hyperlink{bsd}{B+S-D} for $X/K$, }  \\
	+ & \text{ Conjecture~\hyperlink{conj2}{2} for $i=1$.  }
	\end{align*}
\end{corollary}
\begin{proof}
By \cref{mlem}, Conjecture~\hyperlink{BS3}{B+S-D+2} for  $\sX$ is equivalent 
to the vanishing of $V_\ell\Br(\sX_{k^{\sep}})^{G_k}$ 
which is equivalent to the vanishing of $V_\ell\Br(\sY_{k^{\sep}})^{G_k}$,  $V_\ell\Sha(\Pic^0_{X/K}/K)$ and $V_\ell\Br(X_{K^\sep})^{G_K}$ by  \cref{big}. 
By \cref{cor: nsrank}, \cref{thm: bsdconj} and \cref{mlem}, 
this is equivalent to the statement that the Conjectures~\hyperlink{BS3}{B+S-D+2} for $\sY$, \hyperlink{bsd}{B+S-D} for $X/K$ and \hyperlink{conj2}{2} hold, as claimed.
\end{proof}

\begin{remark}
	For the case that $S$ is a smooth projective curve and the structure morphism of the model is proper, this result was firstly proved by Geisser (cf.~\cite[Thm.\,1.1]{Geisser_2021}).
\end{remark}

\begin{corollary}\label[corollary]{cor: conjectures follow from Tate}
	Assume that the Conjecture $\T^1(\sX/k)$ (or equivalently the Tate Conjecture $\T^1(\sX/k,\ell)$ for divisors for some and hence all primes $\ell$) holds for all smooth projective varieties $\sX$ over finite fields of characteristic $p$.
	Let $K$ be a finitely generated field of the same characteristic $p$. 
Let $X/K$ be a smooth projective geometrically connected variety 
with a smooth model $\sX/\sY$.
Then the following conjectures are true:
     \begin{enumerate}
     \item The Tate Conjecture~\hyperlink{conj1}{$\T^1(X/K,\ell)$} for divisors, for any and hence all primes $\ell$,
     \item Conjecture~\hyperlink{conj2}{2} for $i=1$,
     \item Conjecture~\hyperlink{bsd}{B+S-D}, for $X/K$,   
     \item the analytic Tate Conjecture~\hyperlink{BS3}{B+S-D+2} for  $\sX$.
     \end{enumerate}	
\end{corollary}
\begin{proof}
The model $\rho\colon \sX\rightarrow  \sY$ of $X/K$ satisfies the properties of \cref{big}. 
By \cref{mlem} we have the analytic Tate Conjecture~\hyperlink{BS3}{B+S-D+2} for  $\sX$. 
Using \cref{cor1}, we see that Conjecture~\hyperlink{conj2}{2} for $i=1$,  (and therefore also Conjecture~\hyperlink{conj1}{$\T^1(X/K,\ell)$}), 
Conjecture~\hyperlink{bsd}{B+S-D} for $X/K$, and Conjecture~\hyperlink{BS3}{B+S-D+2} for  $\sX$ hold. 
\end{proof}

\begin{remark}
	$\T^1(X,\ell)$ for smooth projective varieties over finite fields has been reduced to $\T^1(X,\ell)$ for smooth projective surfaces over finite fields with $\H^1(X,\cO_X)=0$ by the work of de Jong \cite{deJong_unpublished}, Morrow \cite[Thm.\,4.3]{Morrow_2019} and Yuan \cite[Thm.\,1.6]{Yuan_2021}. 
\end{remark}

%
\section{The $p$-adic obstruction of the BSD rank conjecture} \label{sec:padic obstruction to BSD rank}

In this section, we are interested in a $p$-adic equivalent of \cref{thm: bsdconj}. 
Again, let $k$ be a finite field and  $K$  finitely generated  over $k$.

\begin{theorem} \label[theorem]{BSD in char p}
For an abelian variety $A/ K$ the following statements are equivalent:
\begin{enumerate}
\item The Conjecture~\hyperlink{bsd}{B+S-D} for $A/ K$ holds. 
\item The $p$-torsion subgroup $\Sha(A/K)(p)$ of the Tate--Shafarevich group is of finite exponent. 
\end{enumerate}
\end{theorem}

\begin{proof}
We can reduce to the case that the base field is a rational function field. 
Indeed, $ K$ has a rational subfield $k(\underline{x})\colonequals k(x_1,\ldots,x_{\mathrm{trdeg}(K)})=K(\bbP^{{\mathrm{trdeg}(K)}})$, such that 
the field extension $ K/k(\underline{x})$ is finite separable. 
Now on the one hand, we have by 
\cref{prop: Sha invariant under Weil} an isomorphism
$$
\Sha(A/K)\cong \Sha \big(\Res_{  K/k(\underline{x})}(A)/k(\underline{x})\big).
$$
In particular, $\Sha(A/K)(p)$ is of finite exponent, 
if and only if $\Sha \big(\Res_{  K/k(\underline{x})}(A)/k(\underline{x})\big)(p)$ is of finite exponent. 
On the other hand, by the subsequent \cref{invariance of BSD under Weil restriction}, we have equalities
\begin{align*}
\rk_{\an}(A/  K) & = \rk_{\an}(\Res_{  K/k(\underline{x})}(A)),\\
\rk_{\alg}(A/K) & = \rk_{\alg}(\Res_{K/k(\underline{x})}(A)).
\end{align*}
Hence, Conjecture \hyperlink{bsd}{B+S-D} for $A/K$ holds if and only if it holds for $\Res_{  K/k(\underline{x})}(A)/k(\underline{x}) $.

Thus, let us assume now that $K=K(\underline{x})$ is the function field of $\bbP^{\mathrm{trdeg}(K)}$, 
which means that its (unramified) Brauer group vanishes. 
As a consequence, the exact sequence of \cref{thm: fibrations p-torsion} provides an isomorphism
$$
\ker\big(V_p\Br_{\nr}(K(A)/k)\rightarrow V_p\Br(A_{k^{\sep}})^{G_k}\big) 
\xrightarrow{\sim} V_p\Sha(\Pic^0_{A/  K,\red}/K)
$$
What is more, since $A$ is an abelian variety, $\Br(A_{k^{\sep}})^{G_k}$ 
has finite exponent by \cite[Cor.\,5.3]{dAddezio}. 
Thus we obtain an isomorphism
$$
V_p \Br_{\nr}(K(A)/k) \xrightarrow{\sim} V_p \Sha(\Pic^0_{A/  K,\red}/K). 
$$
In particular,
$V_p\Sha(A/K)(p) \cong V_p\Sha(\Pic^0_{A/  K,\red}/K)(p)$ is trivial if and only if $V_p\Br_{\nr}(K(A)/k)(p)$ is trivial. 
By \cref{thm: obstruction analytic Tate} the latter is the $p$-adic obstruction of the analytic Tate Conjecture~\hyperlink{BS3}{B+S-D+2},  
for a smooth integral  $\sA/k$ with function field $K(A)$, 
which is invariant under birational morphisms according to a result of Tate mentioned in the introduction~\eqref{eq:def z}.

Since by our above simplifications $\Br(S)=0$, the Tate Conjecture $\T^1(S/k)$ holds, and so does by \cref{cor: conjectures follow from Tate} 
the analytic Tate Conjecture~\hyperlink{BS3}{B+S-D+2}. 
Moreover, since $A/  K$ is an abelian variety, the Tate Conjecture~\hyperlink{conj1}{$\T^1(X/K,\ell)$} for divisors is known. 
Thus in this case it follows from \cref{cor1}, that 
the analytic Tate Conjecture~\hyperlink{BS3}{B+S-D+2} is equivalent to 
the Conjecture~\hyperlink{bsd}{B+S-D} for $A/  K$.
\end{proof}

\begin{lemma} \label[lemma]{invariance of BSD under Weil restriction}
Let $L/K$ be a finite separable extension and $A/L$ be an abelian variety with Weil restriction $\Res_{L/K}(A)$ over $K$. 
The following equalities hold:
\begin{align*}
\rk_{\an}(A/L) &= \rk_{\an}(\Res_{L/K}(A)/K),\\
\rk_{\alg}(A/L) &= \rk_{\alg}(\Res_{L/K}(A)/K).
\end{align*}
\end{lemma}

\begin{proof}
Let $\alpha\colon \Spec(L) \to \Spec(K)$.
There is a morphism $\rho\colon \sA \to \sY$ with generic fibre $A/L$ 
together with a finite flat morphism $\alpha\colon \sY \to \sY'$,  
which is a model of $\alpha$.  
Denote the push-forward $\alpha \circ \rho$ of $\rho$ along $\alpha$  by $\rho'$. 
For a prime $\ell \neq p$, consider the Leray spectral sequence 
$$
\R^i\alpha_*\R^j\rho_* \bbQ_\ell \Rightarrow \R^{i+j}\rho'_*\bbQ_\ell.
$$ Since push-forwards by finite morphisms are acyclic for the \'etale cohomology~\cite[Cor.\,II\,3.6]{Milne_1980}, 
the Leray spectral sequence degenerates on the $E_2$-page to give 
isomorphisms 
$$
\alpha_*\sR_\ell^n(\sA/\sY)=\pi_*\R^n\rho_*\bbQ_\ell \cong \R^n\rho'_*\bbQ_\ell=\sR_\ell^n(\sA/\sY').
$$ 
Hence  $\Phi_n(\alpha_*(\sA)/\sY',s) = \Phi_n(\sA/\sY,s)$. 
Note that the sheaf push-forward $\alpha_*A$ is the Weil restriction 
$\Res_{L/K}(A)$. 
Since $\dim(\sY) = \dim(\sY')$, 
we use the above equation for $n=1$ to conclude that 
the analytic ranks of $A/L$ and $\Res_{L/K}(A)/K$ are equal.

For the second equality, 
recall that the Weil restriction of an abelian variety under a 
finite separable field extension is again an abelian variety. 
Moreover, an abelian variety $A$ is isogenous to its dual $\Pic^0(A)$. 
Thus $\rk_{\alg}(A/L)=\rk(A(L))$, because isogenous varieties have the same rank. 
Similarly we have $\rk_{\alg}(\Res_{L/K}(A)/K)= \rk (\Res_{L/K}(A)(K))$. 
But by definition of the Weil restriction, $A(L) \cong \Res_{L/K}(A)(K)$. 
and this finishes the proof.
\end{proof}

\bibliography{p-l-BSD}

\providecommand{\bysame}{\leavevmode\hbox to3em{\hrulefill}\thinspace}
\providecommand{\MR}{\relax\ifhmode\unskip\space\fi MR }
\providecommand{\MRhref}[2]{%
  \href{http://www.ams.org/mathscinet-getitem?mr=#1}{#2}
}
\providecommand{\href}[2]{#2}
\begin{thebibliography}{10}

\bibitem{Bloch_1980}
Spencer Bloch, \emph{A note on height pairings, {T}amagawa numbers, and the {B}irch and {Swinnerton-Dyer} conjecture}, Inventiones mathematicae \textbf{58} (1980), no.~1, 65--76.

\bibitem{Ulm}
Gebhard B{\"o}ckle, David Burns, David Goss, Dinesh Thakur, Fabien Trihan, Douglas Ulmer, and Douglas Ulmer, \emph{Curves and jacobians over function fields}, Arithmetic geometry over global function fields (2014), 281--337.

\bibitem{BombieriGubler_2007}
Enrico Bombieri and Walter Gubler, \emph{Heights in {Diophantine} geometry}, New Mathematical Monographs, vol.~4, Cambridge University Press, 2007.

\bibitem{CadoretHuiTamagawa_2023}
Anna Cadoret, Chun~Yin Hui, and Akio Tamagawa, \emph{{$\bbQ_{\ell}$}-versus {$\bbF_{\ell}$}-coefficients in the {G}rothendieck-{S}erre/{T}ate conjectures}, Israel J. Math. \textbf{257} (2023), no.~1, 71--101.

\bibitem{ColliotTheleneSkorobogatov_2013}
Jean-Luis Colliot-Thélène and {Alexei\,Nikolaievich} Skorobogatov, \emph{Descente galoisienne sur le groupe de {Brauer}}, J. reine angew. Math. \textbf{682} (2013), 141--165.

\bibitem{ColliotTheleneSkorobogatov_2021}
\bysame, \emph{The {B}rauer--{G}rothendieck group.}, 1 ed., Ergebnisse der Mathematik und ihrer Grenzgebiete, 3. Folge, vol.~71, Springer Cham, 2021.

\bibitem{Conrad_2006}
Brian Conrad, \emph{Chow's $k/k$-image and $k/k$-trace, and the {Lang--Néron} theorem}, Enseignement Mathématique \textbf{52} (2006), no.~2, 37--108.

\bibitem{dAddezio}
Marco D'Addezio, \emph{Boundedness of the {$p$}-primary torsion of the {B}rauer group of an abelian variety}, Compos. Math. \textbf{160} (2024), no.~2, 463--480.

\bibitem{deJong_unpublished}
{Aise Johan} {de Jong}, \emph{Tate conjecture for divisors}, Unpublished note.

\bibitem{deJong_1996}
\bysame, \emph{Smoothness, semi-stability and alterations}, Publications mathématiques de l’I.H.É.S. \textbf{83} (1996), 51--93.

\bibitem{Deligne_1968}
Pierre Deligne, \emph{Lefschetz’s theorem and criteria of degeneration of spectral sequences}, Publ.\,math.\,I.H.É.S. \textbf{35 (1)} (1968), 107--126.

\bibitem{Deligne_1980}
\bysame, \emph{La conjecture de {W}eil: {II}}, Publ.\,math.\,I.H.É.S. \textbf{52} (1980), 137--252.

\bibitem{EtesseLestum_1993}
Jean-Yves Etesse and Bernard~Le Stum, \emph{Fonctions {$L$} associées aux {$F$-}isocristaux surconvergents, {I}. {I}nterprétation cohomologique}, Math.\,Ann. \textbf{296} (1993), 557--576.

\bibitem{Fu_2015}
Lei Fu, \emph{Etale cohomology theory}, World Scientific, 2015.

\bibitem{Fujiwara_2002}
Koji Fujiwara, \emph{A proof of the absolute purity conjecture (after {G}abber)}, (Algebraic Geometry 2000, Azumino), Advanced Studies in Pure Mathematics, vol.~36, 2002, pp.~153--183.

\bibitem{fulton2013intersection}
William Fulton, \emph{Intersection theory}, vol.~2, Springer Science \& Business Media, 2013.

\bibitem{Geisser_2021}
{Thomas\,Hermann} Geisser, \emph{Tate's conjecture and the {Tate--Shafarevich} group over global function fields}, Journal of the Institute of Mathematics of Jussieu \textbf{20} (2021), no.~3, 1001--1022.

\bibitem{Goa}
Cristian~D. Gonzalez-Aviles, \emph{Brauer groups and {T}ate-{S}hafarevich groups}, J. Math. Sci. Univ. Tokyo \textbf{10} (2003), no.~2, 391--419.

\bibitem{Grothendieck_1968-III}
Alexandre Grothendieck, \emph{{Le groupe de Brauer.\,III\,: Exemples et complements}}, Dix Exposés sur la Cohomologie des Schémas (Alexandre Grothendieck and Nicolaas~H. Kuiper, eds.), {Advanced Studies Pure Math.}, vol.~3, North Holland Publishing Company, Amsterdam, 1968, pp.~88--188.

\bibitem{Hartshorne}
Robin Hartshorne, \emph{Algebraic geometry}, Graduate Texts in Mathematics, vol. No. 52, Springer-Verlag, New York-Heidelberg, 1977.

\bibitem{hindry2005rang}
Marc Hindry and Am{\'\i}lcar Pacheco, \emph{Sur le rang des jacobiennes sur un corps de fonctions}, Bulletin de la Soci{\'e}t{\'e} Math{\'e}matique de France \textbf{133} (2005), no.~2, 275--295.

\bibitem{Holmes2019}
David Holmes, \emph{Néron models of jacobians over base schemes of dimension greater than 1}, Journal für die reine und angewandte Mathematik (Crelles Journal) \textbf{2019} (2019), no.~747, 109--145.

\bibitem{Illusie_1979}
Luc Illusie, \emph{Complex de {de Rham--Witt} et cohomologie cristalline}, Ann. Sci. Ec. Norm. Sup\'er. (4) \textbf{12} (1979), no.~4, 501--661.

\bibitem{katotrihan2003}
Kazuya Kato and Fabien Trihan, \emph{{On the conjectures of Birch and Swinnerton-Dyer in characteristic $p > 0$}}, Inventiones mathematicae \textbf{153} (2003), 537--592.

\bibitem{Kedlaya_2006}
Kiran~S. Kedlaya, \emph{Fourier transforms and $p$-adic ``{Weil} {II}''}, Compos.\,Math. \textbf{142} (2006), 1426--1450.

\bibitem{Keller_2016}
Timo Keller, \emph{On the {T}ate--{S}hafarevich group of abelian schemes over higher dimensional bases over finite fields}, Manuscripta Mathematica \textbf{150} (2016), 211--245.

\bibitem{Keller_2019}
\bysame, \emph{On an analogue of the conjecture of {Birch} and {Swinnerton-Dyer} for {Abelian} schemes over higher dimensional bases over finite fields}, Documenta Mathematica \textbf{24} (2019), 915--993.

\bibitem{lichtenbaum1983zeta}
Stephen Lichtenbaum, \emph{Zeta-functions of varieties over finite fields at $s = 1$}, Arithmetic and Geometry: Papers Dedicated to IR Shafarevich on the Occasion of His Sixtieth Birthday Volume I Arithmetic (1983), 173--194.

\bibitem{Madapusi_2015}
Keerthi {Madapusi Pera}, \emph{The {T}ate conjecture for {K3}-surfaces in odd characteristic}, Invent.\,Math. \textbf{201} (2015), 625--668.

\bibitem{MadoreOrgogozo_2015}
David~A. Madore and Fabrice Orgogozo, \emph{Calculabilit\'{e} de la cohomologie \'{e}tale modulo {$\ell$}}, Algebra Number Theory \textbf{9} (2015), no.~7, 1647--1739.

\bibitem{Matsusaka_1957}
Teruhisa Matsusaka, \emph{The criteria for algebraic equivalence and the torsion group}, American Journal of Mathematics \textbf{79} (1957), 53--66.

\bibitem{MazurMessing_2006}
Barry Mazur and Walter Messing, \emph{Universal extensions and one dimensional crystalline cohomology}, Lecture Notes in Mathematics, vol. 370, Springer, 2006.

\bibitem{Mil2}
J.~S. Milne, \emph{Comparison of the {B}rauer group with the {T}ate-\v safarevi\v c{} group}, J. Fac. Sci. Univ. Tokyo Sect. IA Math. \textbf{28} (1981), no.~3, 735--743 (1982).

\bibitem{milne1975conjecture}
James~S Milne, \emph{On a conjecture of artin and tate}, Annals of Mathematics \textbf{102} (1975), no.~3, 517--533.

\bibitem{Milne_1980}
James~S. Milne, \emph{{É}tale cohomology}, Princeton Mathematical Series, vol.~33, Princeton University Press, 1980.

\bibitem{Milne_1986_AV}
\bysame, \emph{Abelian varieties}, Arithmetic Geometry (Gary Cornell and Joseph~H. Silverman, eds.), Springer New York, New York, NY, 1986, pp.~103--150.

\bibitem{Milne_1986_JV}
\bysame, \emph{Jacobian varieties}, Arithmetic Geometry (Gary Cornell and Joseph~H. Silverman, eds.), Springer New York, New York, NY, 1986, pp.~167--212.

\bibitem{Milne_2006}
\bysame, \emph{{A}rithmetic duality theorems}, second ed., BookSurge LLC, 2006.

\bibitem{Morrow_2019}
Mathew Morrow, \emph{A variational {T}ate conjecture in crystalline cohomology}, J.\,Eur.\,Math.\,Soc. \textbf{21} (2019), no.~11, 3467--3511.

\bibitem{Petrequin_2003}
Denis Petrequin, \emph{Classes de {Chern} et classes de cycles en cohomologie rigide}, Bull. Soc. Math. France \textbf{131} (2003), 59--121.

\bibitem{Poonen_2017}
Bjorn Poonen, \emph{Rational points on varieties}, Graduate Studies in Mathematics, vol. 186, AMS, 2017.

\bibitem{Pal_2022}
Ambrus Pál, \emph{The $p$-adic monodromy group of abelian varieties over global function fields of characteristic $p$}, Doc.\,Math. \textbf{27} (2022), 1509--1579.

\bibitem{QinBrauer1}
Yanshuai Qin, \emph{On geometric {B}rauer groups and {Tate--Shafarevich} groups}, arXiv: \textbf{2012.01681v2 [math.AG]} (2021), 19 pages.

\bibitem{QinBrauer2}
\bysame, \emph{On the {B}rauer groups of fibrations {II}}, arXiv: \textbf{2103.06910v3 [math.AG]} (2021), 46 pages.

\bibitem{qinbrauer0}
\bysame, \emph{On the {B}rauer groups of fibrations}, Mathematische Zeitschrift \textbf{307} (2024), no.~1, 1--20.

\bibitem{RaynaudGruson_1971}
Michel Raynaud and Laurent Gruson, \emph{Critères de platitude et de projectivité, {T}echniques de platification d'un module}, Invertiones Mathematicae \textbf{13} (1971), 1--89.

\bibitem{rosen1998rank}
Michael Rosen and Joseph~H Silverman, \emph{On the rank of an elliptic surface}, Inventiones mathematicae \textbf{133} (1998), 43--67.

\bibitem{Schneider_1982}
Peter Schneider, \emph{On the conjecture of {Birch} and {Swinnerton-Dyer} about global function fields}, Mathematische Annalen \textbf{260} (1982), no.~4, 495--510.

\bibitem{Serre_1965}
Jean-Pierre Serre, \emph{Zeta and {$L$}-functions}, Arithmetical Algebraic Geometry, Harper \& Row, 1965, Proc.\,of a Conference held at Purdue Univ., Dec.\,5--7 1963, pp.~82--92.

\bibitem{TateBSD1965}
John Tate, \emph{{On the conjectures of Birch and Swinnerton-Dyer and a geometric analog}}, S{\'e}minaire Bourbaki \textbf{9} (1965), no.~306, 415--440.

\bibitem{Tate_1965}
John~T. Tate, \emph{Algebraic cycles and poles of zeta functions}, Arithmetical Algebraic Geometry (Proc.\,Conf.\,Purdue Univ., 1963) (New York), Harper \& Row, 1965, pp.~93--110.

\bibitem{Tate_1966}
\bysame, \emph{Endomorphisms of abelian varieties over finite fields}, Invent.\,Math. \textbf{2} (1966), 134--144.

\bibitem{Tate_1994}
\bysame, \emph{Conjectures on algebraic cycles in $\ell$-adic cohomology}, Motives (Seattle, WA, 1991), 1994, pp.~71--83.

\bibitem{Yuan_2021}
Xinyi Yuan, \emph{Positivity of {H}odge bundles of abelian varieties over some function fields}, Compositio Mathematica \textbf{157} (2021), 1964--2000.

\bibitem{Yua2}
\bysame, \emph{Comparison of arithmetic {B}rauer groups with geometric {B}rauer group}, arXiv: \textbf{2011.12897 [math.AG]} (version 2), 21 pages.

\bibitem{Zarhin_1974a}
Yuri~G. Zarhin, \emph{Isogenies of abelian varieties over fields of finite characteristic}, Mat.\,Sb.\,(N.S.) \textbf{95} (1974), no.~137, 461--470 \& 472.

\bibitem{Zarhin_1974b}
Yuri~G Zarhin, \emph{A remark on endomorphisms of abelian varieties over function fields of finite characteristic}, Izv.\,Akad.\,Nauk SSSR Ser.\,Mat.\, \textbf{38} (1974), 471--474.

\end{thebibliography}
\bibliographystyle{amsplain}

\end{document}